\documentclass[11pt]{article}

\pdfoutput=1
\usepackage{mathrsfs}
\usepackage{fancyhdr}
\usepackage{epsfig,morefloats}
\usepackage[authoryear]{natbib}
\setcitestyle{aysep={,}}
\usepackage{amsmath,enumerate,amsbsy,amsfonts,amssymb,mathabx,amscd,graphicx,algorithm,amsthm,soul}
\usepackage{adjustbox}
\usepackage{enumitem}
\usepackage{color}
\usepackage{lscape}
\usepackage{rotating}
\usepackage{caption}
\usepackage{subfigure}
\usepackage[compact]{titlesec}
\usepackage{lineno}
\usepackage{cancel}
\usepackage{multirow}
\usepackage{setspace}

\definecolor{mygray}{gray}{0.6}

\newcommand{\Date}[1]{\def\@Date{#1}}
\def\today{\number\day~\ifcase\month\or
	January\or February\or March\or April\or May\or June\or
	July\or August\or September\or October\or November\or December\fi~\number\year}

\def\t0{\theta_0}

\def\bsc{\begin{scriptsize}}
	\def\esc{\end{scriptsize}}

\newtheorem{theorem}{Theorem}

\newtheorem{lemma}{Lemma}

\newtheorem{proposition}{Proposition}
\theoremstyle{definition}
\newtheorem{definition}{Definition}
\newtheorem{remark}{Remark}
\newtheorem{condition}{Condition}

\makeatletter
\newcommand{\figcaption}{\def\@captype{figure}\caption}
\newcommand{\tabcaption}{\def\@captype{table}\caption}
\makeatother

\newcommand{\beps}{\boldsymbol \epsilon}
\newcommand{\bLambda}{\boldsymbol \Lambda}
\newcommand{\bOmega}{\boldsymbol \Omega}
\newcommand{\bSigma}{\boldsymbol \Sigma}

\newcommand{\blind}{1}

\newcommand{\bpsi}{\boldsymbol \psi}

\newcommand{\bphi}{\boldsymbol \phi}
\newcommand{\bPhi}{\boldsymbol \Phi}
\newcommand{\bmu}{\boldsymbol \mu}

\newcommand{\btheta}{\boldsymbol \theta}

\newcommand{\bdelta}{\boldsymbol \delta}
\newcommand{\bbeta}{\boldsymbol \beta}
\newcommand{\bnu}{\boldsymbol \nu}
\newcommand{\bfeta}{\boldsymbol \eta}
\newcommand{\bvar}{\boldsymbol \varepsilon}

\newcommand{\bvarepsilon}{\boldsymbol \varepsilon}
\newcommand{\bzeta}{\boldsymbol \zeta}
\def\p{{ \mathrm{p} }}
\newcommand{\ba}{{\mathbf a}}
\newcommand{\be}{{\mathbf e}}
\newcommand{\bbf}{{\mathbf f}}
\newcommand{\bx}{{\mathbf x}}
\newcommand{\bz}{{\mathbf z}}
\newcommand{\bg}{{\mathbf g}}
\newcommand{\by}{{\mathbf y}}

\newcommand{\bv}{{\mathbf v}}

\newcommand{\br}{{\mathbf r}}

\newcommand{\bbb}{{\mathbf b}}

\newcommand{\bd}{{\mathrm d}}

\newcommand{\bD}{{\bf D}}
\newcommand{\bA}{{\bf A}}
\newcommand{\bB}{{\bf B}}

\newcommand{\bG}{{\bf G}}

\newcommand{\bK}{{\bf K}}

\newcommand{\bX}{{\bf X}}
\newcommand{\bY}{{\bf Y}}
\newcommand{\bZ}{{\bf Z}}
\newcommand{\bR}{{\bf R}}
\newcommand{\bU}{{\bf U}}
\newcommand{\bV}{{\bf V}}
\newcommand{\bW}{{\bf W}}

\newcommand{\cM}{{\cal M}}

\newcommand{\cU}{{\cal U}}
\newcommand{\cS}{{\cal S}}

\newcommand{\eZ}{\mathbb{Z}}
\newcommand{\eE}{\mathbb{E}}
\newcommand{\eR}{\mathbb{R}}
\newcommand{\cH}{\mathbb{H}}

\newcommand{\tF}{\text{F}}

\newcommand{\cov}{\text{Cov}}

\newcommand{\tr}{\mbox{tr}}

\newcommand{\T}{\mathrm{\scriptscriptstyle T}}

\DeclareMathOperator*{\esssup}{ess\,sup}

\allowdisplaybreaks

\numberwithin{remark}{section}

\usepackage{geometry}
\usepackage{authblk}

\begin{document}
	
	% \def\spacingset#1{\renewcommand{\baselinestretch}%
		% {#1}\small\normalsize} \spacingset{1}

	%\renewcommand{\baselinestretch}{1.0}
	
	%%%%%%%%%%%%%%%%%%%%%%%%%%%%%%%%%%%%%%%%%%%%%%%%%%%%%%%%%%%%%%%%%%%%%%%%%%%%%%
	
	\if1\blind
	{
		\title{\bf An autocovariance-based learning framework for high-dimensional functional time series\thanks{The authors equally contributed to the paper. We thank the editor, the associate editor and two anonymous referees for their constructive comments and suggestions. Chang and Chen were supported in part by the National Natural Science Foundation of China (grant nos.~71991472, 72125008 and 11871401). Chang was also supported by the Center of Statistical Research at Southwestern University of Finance and Economics. Yao was supported in part by the U.K. Engineering and Physical Sciences Research Council (grant no. EP/V007556/1).
				
				Corresponding author: Xinghao Qiao and Qiwei Yao. Email addresses: changjinyuan@swufe.edu.cn (J. Chang), chenc@swufe.edu.cn (C. Chen), x.qiao@lse.ac.uk (X. Qiao), q.yao@lse.ac.uk (Q. Yao)
			}\hspace{.2cm}\\}
		
		\author[1]{Jinyuan Chang}
		\author[1]{Cheng Chen}
		\author[2]{Xinghao Qiao}
		\author[2]{Qiwei Yao}
		\affil[1]{\it \small Joint Laboratory of Data Science and Business Intelligence, Southwestern University of Finance and Economics, Chengdu, Sichuan 611130, China
		}
		\affil[2]{\it \small Department of Statistics, London School of Economics, London, WC2A 2AE, UK
		}
		
		\date{}
		
		\maketitle
	} \fi
	
	\if0\blind
	{
		\bigskip
		\bigskip
		\bigskip
		\begin{center}
			{\LARGE\bf An autocovariance-based learning framework for high-dimensional functional time series}
		\end{center}
		\medskip
	} \fi

	%\maketitle
	
	\bigskip
	
	\begin{abstract}
		Many scientific and economic applications involve the statistical learning of high-dimensional functional time series, where the number of functional variables is comparable to, or even greater than, the number of serially dependent functional observations. In this paper, we model observed functional time series, which are subject to errors in the sense that each functional datum arises as the sum of two uncorrelated
		components, one dynamic and one white noise. Motivated from the fact that the autocovariance function of observed functional time series automatically filters out the noise term, we propose a three-step framework by first performing autocovariance-based dimension reduction, then formulating a novel autocovariance-based block regularized
		minimum distance estimation to produce block sparse estimates, and based on which obtaining the final functional sparse estimates. We investigate theoretical properties of the proposed estimators, and illustrate the proposed  estimation procedure with the corresponding convergence analysis via three sparse high-dimensional functional time series models. We demonstrate via both simulated and real datasets that our proposed estimators significantly outperform their competitors.
	\end{abstract}

	\bigskip \bigskip
	\noindent%
	{\it Key words:}  Block regularized minimum distance estimation;
	Dimension reduction; Functional time series; High-dimensional data;
	Non-asymptotics; Sparsity.
	
	\bigskip
	%\begin{quote}
	\noindent
	{\it JEL code}: C13, C32, C55
	%{\sl MSC2010 subject classifications}: Primary 62G99; secondary 62F40
	%\end{quote}
	
	%\vfill
	
	%\begin{quote}
	%\noindent
	%{\sl Keywords}: Empirical likelihood; Estimating equations; High-dimensional statistical methods; Misspecification; Moment selection; Penalized likelihood.
	%\end{quote}

	%\thispagestyle{empty}
	%\pagenumbering{gobble}
	
	%\newpage
	%\tableofcontents
	
	%\newpage
	%\pagenumbering{arabic}
	%\setcounter{page}{1}
	
	%\linenumbers

	\newpage
	% \spacingset{1.1} % DON'T change the spacing!
	\onehalfspacing

	\section{Introduction}
	\label{sec.intro}
	Functional time series refers to functional data objects that are
	observed consecutively over time. % and constitutes an active research area.
	Existing research on functional time series has mainly focused on
	extending the univariate or
	low-dimensional multivariate time series methods to the functional domain.
	% with theoretical guarantees  under an asymptotic framework.
	An incomplete list of the relevant references includes
	\cite{Bbosq1}, \cite{bathia2010}, \cite{hormann2010}, \cite{panaretos2013}, \cite{aue2015}, \cite{hormann2015}, \cite{li2020} and \cite{chen2022}. % just to list a few.
	The rapid development of data collection
	technology has made high-dimensional functional time series datasets
	increasingly common. Examples include hourly measured
	concentrations of various pollutants such as PM10 trajectories
	\citep{hormann2015} collected at different measuring stations, daily
	electricity load curves \citep{cho2013} for a large number of
	households, cumulative intraday return trajectories \citep{horvath2014},
	daily return density curves \citep{bathia2010} 
	and functional volatility processes \citep{muller2011} for a 
	collection of stocks.
	%and annul temperature curves \citep{aue2020} at different measuring stations.
	
	We consider in this paper a setting for
	modelling high-dimensional functional time series as follows.
	Let $\bW_t(\cdot) = \{W_{t1}(\cdot), \dots, 
	W_{tp}(\cdot)\}^{\T}$, $t=1, \dots,n$, be the observed $p$-vector of
	functional time series defined on a compact interval $\cU$,
	% $\bW_t(\cdot) = \{W_{t1}(\cdot), \dots, W_{tp}(\cdot)\}^{\T}$ for $t=1, \dots,n$
	% with (auto)covariance functions
	% $\bSigma_h^W(u,v) = \cov\{\bW_t(u),\bW_{t+h}(v)\}$ for any integer $h$
	% and $u,v \in \cU$ (a compact interval),
	where the dimension $p$ is large
	in relation to $n$, and $p$ may be greater than $n$.
	% regime, $p$ can diverge with, or even larger than, $n.$
	Suppose that 
	$\bW_t(\cdot)$ is subject to an error:
	% in the form of
	\begin{equation}
		\label{error.model}
		\bW_{t}(\cdot) =  \bX_{t}(\cdot) + \be_{t}(\cdot)\,, % ~~~~u \in \cU\,,
	\end{equation}
	where $\bX_t(\cdot)=\{X_{t1}(\cdot),\dots, X_{tp}(\cdot)\}^{\T}$ is a
	functional time series of interest, 
	$\be_t(\cdot)=\{e_{t1}(\cdot),\dots, \linebreak e_{tp}(\cdot)\}^{\T}$ is white
	noise in the sense (\ref{WN}) below, and $\{ \bX_t(\cdot) \}_{t=1}^n$
	and $\{ \be_t(\cdot)\}_{t=1}^n$ are uncorrelated. Note that both $ \bX_{t}(\cdot)$ and
	$\be_{t}(\cdot)$ are latent.
	We assume that
	both $\bW_{t}(\cdot)$ and $\bX_t(\cdot)$ are weakly stationary,
	and $\mathbb{E}\{\bW_t(\cdot)\} ={\bf 0}$. 
	For any integer $h$ and $u,v \in \cU$, put
	\begin{equation}
		\label{n12}
		\bSigma_h^W(u,v) = \cov\{\bW_{t-h}(u),\bW_{t}(v)\}\,,~~
		\bSigma_h^X(u,v) = \cov\{\bX_{t-h}(u),\bX_{t}(v)\}\,,~~
		\bSigma_h^e(u,v) = \cov\{\be_{t-h}(u),\be_{t}(v)\}\,. 
	\end{equation}
	%In the same manner as $\bSigma_h^W(u,v),$ we define $\bSigma_h^X(u,v)$ and $\bSigma_h^e(u,v)$
	%by replacing $\bW_t(\cdot)$ with $\bX_t(\cdot)$ and $\be_t(\cdot)$, respectively.
	%We define $\bSigma_h^X(u,v)$ and $\bSigma_h^e(u,v)$ in the same manner as $\bSigma_h^W(u,v)$ with replacing $\bW_t(\cdot)$ by $\bX_t(\cdot)$ and $\be_t(\cdot)$, respectively.
	We call $\be_t(\cdot)$ a  white noise if 
	\begin{equation} \label{WN}
		\mathbb{E}\{\be_t(u)\} = {\bf 0}~\textrm{and}~\bSigma_h^e(u,v)={\bf 0}~\textrm{for any}~u,v\in \cU~{\rm and}~h \neq 0\,.
	\end{equation}
	Furthermore, we assume that $\ba^{\T} \bX_t(\cdot)$ is not white noise
	for any non-zero constant vector $\ba\in\mathbb{R}^p $.
	Under this setting, the linear dynamic structure of $\bW_t(\cdot)$
	is entirely determined by that of
	$\bX_t(\cdot),$ and all white noise
	elements $\bW_t(\cdot)$ are absorbed into $\be_t(\cdot).$ The presence of
	$\be_t(\cdot)$ reflects that signal curves $\bX_t(\cdot)$ are seldom
	completely observed. Instead, they are often only measured, with errors,
	on a grid. These noisy discrete data are smoothed to yield
	`observed' curves $\bW_t(\cdot).$
	See \cite{bathia2010} for the univariate version of
	model~(\ref{error.model}).
	% with fully nonparametric structure on $\bSigma^e_0.$
	When $\bW_1(\cdot), \dots, \bW_n(\cdot)$ are univariate and independent,
	\cite{hall2006} considered the same model under a `low noise' setting
	assuming that $\be_t(\cdot)$ goes to 0 as $n$ grows to $\infty.$ 
	To separate $\bX_t(\cdot)$ from $\be_t(\cdot)$, e.g., via the covariance function, even in the univariate case,
	% \st{To separate $\bSigma_0^X$ from $\bSigma_0^W$ even in the univariate case,}
	some special structures were imposed; see,
	e.g., %Imposing some parametrically specified structure in the univariate case,
	diagonal $\bSigma^e_0$ of \cite{yao2005} and banded
	%under the assumption that
	$\bSigma^e_0$ of \cite{descary2019}.
	In contrast, we do not impose any structures
	on $\bSigma^e_0$ in this paper, { and our estimation filters
		out the impact of $\be_t(\cdot)$ % via $\bSigma^e_h$
		automatically.} 
	% \st{ and our estimation filters out $\bSigma^e_0$ automatically.} {\color{blue} (XQ: our estimation filters out $\bSigma^e_h$ not $\bSigma^e_0$ automatically, so either make the above change or remove?)}
	
	The standard estimation procedures for univariate functional time series
	models usually consist of three steps \citep{aue2015}.
	% Due to the intrinsic infinite-dimensionality of functional data, 
	Dimension-reduction is performed  first %in the first step 
	via, e.g., functional
	principal components analysis
	(FPCA). Each observed curve is then approximated by a finite truncation. This 
	effectively transforms functional time series into a vector time series
	of FPC scores.
	In the second step the estimation of the function-valued parameters
	in the model is transformed  to that of some appropriate parameter {vectors/matrices} based on estimated FPC scores.
	Finally %in the third step
	the estimated principal component functions are utilized to obtain function-valued estimates based on the estimated parameter
	{vectors/matrices}. %{\color{blue} (XQ: I change matrices to vector/matrices, since vectors comes from functional covariates and scalar responses and matrices come from functional covariates and functional responses.)}
	% Estimation for high-dimensional functional time series models is often impossible without imposing some lower-dimensional structural assumption on the model parameters space. Under functional sparsity assumptions
	To overcome the difficulties caused by high-dimensionality (i.e. large $p$
	in relation to $n$),
	{ some functional sparsity assumptions are %commonly 
		imposed, %on the model parameter space, 
		which results in the estimation under block sparsity constraints in the second step}
	%we impose in the second step a functional sparsity assumption, which results in a block sparsity constraint 
	%{\color{blue} (XQ: I make this change since the functional sparsity assumptions are imposed on the functional model, not in the second step and different models may have different functional sparsity assumptions.)}
	in the sense that
	variables belonging to the same block (or group) are simultaneously included or excluded. In regression setups, the group-lasso penalized least squares estimation \citep{yuan2006} is often adopted in the second step to obtain block sparse estimates. %, from which the third step can recover functional sparse estimates. 
	Similar three-step procedures have been developed to estimate sparse
	high-dimensional functional models, see, e.g., vector functional
	autoregression (VFAR) \citep{guo2022}, scalar-on-function linear
	additive regression (SFLR) \citep{fan2015,kong2016,xue2020,fang2022} and
	function-on-function linear additive regression (FFLR)
	\citep{fan2014b,luo2017,fang2022}.
	% with serially dependent observations. 
	However,
	those estimation procedures are developed under an assumption that signal
	curves are observed directly. 
	
	In our setting the observed curves $\bW_t(\cdot)$ are subject to
	the error contamination as in model (\ref{error.model}). 
	Both FPCA and penalized least squares estimation based on the estimated
	covariance function %of $\bW_t(\cdot),$ i.e.
	$\widehat{\bf \Sigma}_0^W$ % the standard covariance-based procedure 
	are inappropriate since 
	$\bSigma^W_0 = \bSigma^X_0 + \bSigma^e_0$ and, hence,
	$\widehat{\bf \Sigma}^W_0$ is no longer a consistent estimator for $\bSigma^X_0.$
	Motivated from the fact that $\bSigma^W_h = \bSigma^X_h$ for any $h \neq
	0,$ which automatically removes the impact from the noise $\be_t(\cdot)$
	and ensures that the estimator for $\bSigma^W_h$ is also legitimate for
	$\bSigma^X_h,$ 
	we propose an autocovariance-based three-step learning procedure.
	Differing from FPCA based on the Karhunen--Lo\`eve expansion, % of $W_{tj}(\cdot),$ 
	our first dimension reduction step is formulated under an alternative %
	data-driven basis
	expansion of $X_{tj}(\cdot)$ based on the eigenanalysis of a
	positive-definite operator defined in terms of the
	autocovariance functions of $W_{tj}(\cdot).$
	Different from the penalized least squares estimation, our second step
	makes use of the autocovariance of basis coefficients to construct
	high-dimensional moment equations and then applies the proposed block
	regularized method to estimate the associated block sparse 
	parameter vectors/matrices.
	% based on estimated basis coefficients obtained in the first step.
	Our third step re-transforms block sparse estimates to
	functional sparse estimates via estimated basis functions obtained in the
	first step.
	
	Our theoretical development 
	% non-asymptotics There are several challenges in the theoretical analysis of our
	% autocovariance-based learning framework for high-dimensional functional
	% time series, which
	stands at the intersection between high-dimensional
	statistics and functional time series, facing
	several challenges due to
	non-asymptotics %\citep{Bwainwright2019} 
	and infinite-dimensionality with serial dependence. %\citep{jirak2016}. 
	Firstly, in the proposed second step we deal with the estimated basis
	coefficients % rather than true coefficients
	to produce block sparse
	estimates whereas the conventional sparse estimation is applied directly
	to observed data. Accounting for such approximation is a major
	undertaking. Secondly,
	under a high-dimensional and dependent setting, it
	is essential to develop non-asymptotic error bounds on the relevant estimated
	terms as a function of $n,$ $p$ and the truncated dimension, and to
	assess how the serial dependence affects non-asymptotic results.
	%non-asymptotic results.%under our autocovariance-based dimension reduction framework. 
	Thirdly, compared to non-functional data, the infinite-dimensional
	nature of functional data leads to the additional theoretical complexity
	that arises from specifying the block structure and controlling bias
	terms formed by truncation errors from the dimension reduction step.
	
	The main contribution of our paper is three-fold.
	\begin{enumerate}
		\item[1.] Our autocovariance-based learning framework can address the error contamination model (\ref{error.model}) in the presence of infinite-dimensional signal curve dynamics with the addition of `genuinely functional' noise. It makes the good use of the serial correlation information, which is the most relevant in the context of time series modelling.
		
		\item[2.] To provide theoretical guarantees for the first and the third steps and to verify imposed high-level regularity conditions in the second step, we %rely on functional stability measures \citep{guo2022,fang2022} to characterize the effect of serial dependence and 
		establish useful non-asymptotic error bounds on the relevant estimated terms under the autocovariance-based dimension reduction framework.
		
		\item[3.]  We utilize the autocovariance among basis coefficients to construct high-dimensional moment equations with partitioned group structure, based on which we formulate the second step in a novel block regularized minimum distance (RMD) estimation framework to produce block sparse estimates. %Within such framework, 
		The group information can be explicitly encoded in a convex optimization targeting at minimizing the block $\ell_1$ norm objective function subject to the block $\ell_{\infty}$ norm constraint. To theoretically support the second step, we investigate convergence properties of the block RMD estimator. Besides being useful in the second step, the block RMD estimation framework itself is of independent interest and can be applied more broadly.
		
		% \item[4.]  We illustrate the autocovariance-based three-step procedure using three examples of sparse high-dimensional functional time series models, i.e. SFLR, FFLR and VFAR. Theoretically, we study convergence rates of the associated estimators in these models. Empirically, we demonstrate the superiority of these autocovariance-based estimators relative to the covariance-based competitors.
	\end{enumerate}
	
	Our paper is set out as follows. 
	%In Section~\ref{sec.auto3}, we present our autocovariance-based three-step procedure and illustrate it using SFLR as an example. 
	In Section~\ref{sec.dr}, we present Step 1, i.e. 
	the autocovariance-based dimension
	reduction technique.
	We also establish some essential deviation bounds 
	on the relevant estimated terms. % used in subsequent analysis. 
	In Section~\ref{sec.rmd}, we use an example to illustrate the construction of high-dimensional moment equations. We then formulate a general block RMD estimation method (i.e. Step 2) and investigate its theoretical properties. In Section~\ref{sec_3app}, we illustrate the proposed three-step framework using three examples of sparse high-dimensional functional time series models, i.e. SFLR, FFLR and VFAR. Theoretically, we study convergence rates of the associated estimators in these models.
	%we illustrate the proposed autocovariance-based learning framework using three applications, and present the corresponding convergence analysis. 
	In Section~\ref{sec.emp}, we examine the finite-sample performance of the proposed estimators through
	both simulations and an analysis of a real financial dataset. All technical proofs are relegated to the Appendix.
	
	{\bf Notation}. For a positive integer $q,$ we denote $[q]=\{1, \dots,
	q\}.$ Let  $L_2(\cU)$ be a Hilbert space of square-integrable functions
	on a compact interval $\cU.$ The inner product of $f,g \in L_2(\cU)$ is defined as
	$\langle f,g \rangle=\int_{\cU} f(u)g(u)\,{\rm d}u$. For a Hilbert space $\mathbb{H} \subset L_2(\cU),$ we denote the $p$-fold Cartesian product by $\mathbb{H}^p = \mathbb{H} \times \cdots \times \mathbb{H}$ and the tensor product by $\mathbb{S} = \mathbb{H} \otimes \mathbb{H}.$
	For $\bbf=(f_1, \ldots, f_p)^{\T}$ and $\bg=(g_1, \dots,g_p)^{\T}$ in $\cH^p,$ we define $\langle\bbf,\bg\rangle=\sum_{i=1}^p\langle f_i, g_i\rangle.$ We use $\|\bbf\| = \langle \bbf, \bbf \rangle^{1/2}$ and $\|\bbf\|_0 = \sum_{i=1}^p I(\|f_i\| \neq 0)$ with $I(\cdot)$ being the indicator function to denote functional versions of induced norm and $\ell_0$-norm, respectively. For an integral operator $\bK: \mathbb{H}^p \rightarrow \mathbb{H}^q$ induced from the kernel function $\bK=(K_{ij})_{q \times p}$ with each $K_{ij} \in \mathbb{S},$ $\bK(\bbf)(u)=\{\sum_{j=1}^p\langle K_{1j}(u,\cdot),f_{j}(\cdot)\rangle, \dots, \sum_{j=1}^p\langle K_{qj}(u,\cdot),f_{j}(\cdot)\rangle\}^{\T} \in \mathbb{H}^q$ for any $\bbf=(f_1,\ldots,f_p)^{\T} \in \mathbb{H}^p.$ For notational economy, we will also use $\bK$ to denote both the kernel and the operator.
	We define functional versions of Frobenius and matrix $\ell_{\infty}$-norms by $\|\bK\|_{\tF} = (\sum_{i=1}^{q}\sum_{j=1}^p\|K_{ij}\|_{\cS}^2)^{1/2}$ and $\|\bK\|_{\infty} = \max_{i\in[q]} \sum_{j=1}^p \|K_{ij}\|_{\cS},$ respectively, where $\|K_{ij}\|_{\cS}=\{\int_{\cU}\int_{\cU}K_{ij}^2(u,v)\,{\rm d}u{\rm d}v\}^{1/2}$ denotes the Hilbert--Schmidt norm of $K_{ij}.$ For any real matrix $\bB=(b_{ij})_{q\times p}$, we write $\|\bB\|_{\max}=\max_{i\in[q],j\in[p]}|b_{ij}|$ and use $\|\bB\|_{\tF}=(\sum_{i=1}^q\sum_{j=1}^p|b_{ij}|^2)^{1/2}$ and $\|\bB\|_2=\lambda_{\max}^{1/2}(\bB^\T\bB)$ to denote its Frobenius norm and $\ell_2$-norm, respectively. For two sequences of positive numbers $\{a_n\}$ and $\{b_n\}$, we write $a_n\lesssim b_n$ or $b_n\gtrsim a_n$ if there exist a positive constant $c$ such that $a_n/b_n\leq c$. We write $a_n\asymp b_n$ if and only if $a_n\lesssim b_n$ and $b_n\lesssim a_n$ hold simultaneously.

	\section{Autocovariance-based dimension reduction}
	\label{sec.dr}
	
	\subsection{Methodology}
	\label{sec.dr.method}
	
	Our Step 1 is to approximate each curve $X_{tj}(\cdot)$ by a finite linear combination:
	we expand curve $X_{tj}(\cdot)$  using the data-driven
	orthonormal basis functions $\{\psi_{jl}(\cdot)\}_{l=1}^{\infty},$ and truncate
	the expansion to the first $d_j$ (to be specified in Section \ref{sec.emp}) terms:
	\begin{equation}
		\label{expansion}
		X_{tj}(\cdot) = \sum_{l=1}^{\infty}\eta_{tjl}\psi_{jl}(\cdot) \approx \bfeta_{tj}^{\T}\bpsi_{j}(\cdot)\,,~~ j \in [p]\,,
	\end{equation}
	where
	%$\{\psi_{jl}(\cdot)\}_{l=1}^{\infty}$ are orthonormal basis functions, $d_j$ is an appropriate integer,
	$\eta_{tjl} = \langle X_{tj}, \psi_{jl} \rangle,$
	$\bfeta_{tj}=(\eta_{tj1}, \dots, \eta_{tjd_j})^{\T}\in\mathbb{R}^{d_j}$
	and $\bpsi_j(\cdot)=\{\psi_{j1}(\cdot), \dots, \psi_{jd_j}(\cdot)\}^{\T}.$
	Different from the conventional Karhunen--Lo\`eve expansion, the eigenvalues $\lambda_{j1} \geq \lambda_{j2} \geq \dots >0$ and the corresponding eigenfunctions
	$\psi_{j1}(\cdot), \psi_{j2}(\cdot), \dots$ are taken from the
	spectral decomposition of an operator defined as
	\begin{equation}
		\label{df.K}
		K_{jj}(u,v) = \sum_{h=1}^L \int_{\cU}   \Sigma_{h,jj}^X(u,z) \Sigma_{h,jj}^X(v, z)\,\bd z\,,
		%=\sum_{h=1}^L \int_{\cU} \Sigma_{h,jj}^W(u,z) \Sigma_{h,jj}^W(v, z)\,\bd z\,,
	\end{equation}
	where $L>0$ is some prescribed fixed integer, and
	$\Sigma_{h,ij}^X(u,v)$ denotes the $(i,j)$-th element of
	$\bSigma_h^X(u,v)$ in (\ref{n12}). %See (\ref{n13}) below.
	Also denote by
	$\Sigma_{h,ij}^W$ and $\Sigma_{h,ij}^e$ the $(i,j)$-th
	element of, respectively, $\bSigma_h^W$ and $\bSigma_h^e$.
	The idea of using non-zero lagged autocovariances was initiated by  \cite{bathia2010}. A direct consequence is the identity
	\[
	K_{jj}(u,v)=\sum_{h=1}^L \int_{\cU} \Sigma_{h,jj}^W(u,z) \Sigma_{h,jj}^W(v, z)\,\bd z\,,
	\]
	since $\Sigma_{h,jj}^X(u,z) =\Sigma_{h,jj}^W(u,z) $ for all $(u,z) \in \cU^2$ and
	$h\ne 0.$
	This paves the way to estimate $K_{jj}$, and therefore also
	$\bpsi_{j}(\cdot)$, directly based on observations $W_{1j}(\cdot), \dots, W_{nj}(\cdot).$
	{\color{black}The impact of the noise terms $e_{tj}(\cdot)$ is filtered out automatically. 
		It is worth noting that we choose not to use autocovariance functions $\Sigma^W_{h,jj}$ directly in defining $K_{jj}$ as they are not nonnegative definite.
		The definition of $K_{jj}$ in (\ref{df.K})
		ensures that it is nonnegative definite, and there is
		no cancellation of the information
		accumulated from lags 1 to $L$. Hence the estimation is not 
		sensitive to the choice of $L$.}
	In practice, we choose small
	$L$ such as $1\le L \le 5$, as the most significant autocorrelations typically occur at small lags. 
	
	In the standard Karhunen--Lo\`eve expansion, $\{\psi_{jl}(\cdot)\}_{l=1}^\infty$ is deduced from the spectral
	decomposition of $\Sigma_{0,jj}^X$. Since
	\[
	\Sigma_{0,jj}^X(u,v) = \Sigma_{0,jj}^W(u,v) - \Sigma_{0,jj}^e(u,v)\,,
	\]
	some strong assumptions have to be imposed to eliminate the 
	impact of $\Sigma_{0,jj}^e(u,v)$ in order to obtain consistent
	estimates for $\psi_{jl}(\cdot)$. For example, \cite{hall2006} assumes
	that $W_{1j}(\cdot), \dots, W_{nj}(\cdot)$ are independent and the noise $e_{tj}(\cdot)$ goes to 0 as $n$ grows to $\infty.$ 
	{\color{black} Note that the dimension reduction via FPCA can also be performed based on the spectral decomposition of $\Sigma_{0,jj}^W$ instead of $\Sigma_{0,jj}^X,$ as any basis could be used for expanding the data. However, because of $\Sigma_{0,jj}^W=\Sigma_{0,jj}^X+\Sigma_{0,jj}^e,$ using $\Sigma_{0,jj}^W$ may require a larger truncated dimension to capture the sufficient signal information, leading to reduced statistical efficiency.} %Moreover, when the same truncated dimension is used, there is a risk that the leading FPCs selected by $\Sigma_{0,jj}^X$ are not selected by $\Sigma_{0,jj}^W,$ resulting in incorrect inference.}
%{\color{blue} To QY: How about adding the above comment to answer Referee~2's 2nd comment?}
%\QY{I would suggest to delete the last sentence, as the earlier argument would imply to use different number of terms with two covariance matrices -- QY}.
It is also worth mentioning that the penalized least squares approach adopted in the covariance-based second step is based on $\Sigma_{0,jk}^X(u,v)=\Sigma_{0,jk}^W(u,v)-\Sigma_{0,jk}^e(u,v)$ and hence is inappropriate under model~(\ref{error.model}).

With the available observations $\{\bW_t(\cdot)\}_{t \in [n]}$, 
%define the sample (auto)covariance function
% with mean zero and (auto)covariance functions
%$\bSigma_h^W(u,v)=\{\Sigma^W_{h,jk}(u,v)\}_{j,k \in [p]}$
% Define for integer $h \geq 0$ and $(u,v) \in \cU^2.$
% whose sample estimators are 
%With legitimate estimators $\widehat\Sigma_{h,jj}^W$ for positive integer $h$ in (\ref{est.Sigma}),
a natural estimator for $K_{jj}$ in (\ref{df.K}) is defined as
\begin{equation}
	\label{est.K}
	\begin{split}
		\widehat K_{jj}(u,v)=&~\sum_{h=1}^L \int_{\cU} \widehat \Sigma_{h,jj}^W(u,z) \widehat \Sigma_{h,jj}^{W}(v, z)\,\bd z \\
		% =&~\frac{1}{(n-L)^2} \sum_{h=1}^L \sum_{t,s=1}^{n-L} W_{tj}(u)W_{sj}(v) \langle W_{(t+h)j}, W_{(s+h)j}\rangle\,,
		=&~\frac{1}{(n-L)^2} \sum_{h=1}^L \sum_{t,s=h+1}^n W_{(t-h)j}(u)W_{(s-h)j}(v) \langle W_{tj}, W_{sj}\rangle\,,
	\end{split}
\end{equation}
where
\begin{equation}
	\label{est.Sigma}
	\widehat{\bf \Sigma}_h^W(u,v) = \frac{1}{n-h}\sum_{t=h+1}^n\bW_{t-h}(u)
	\bW_t(v)^{\T}=\{\widehat{\Sigma}_{h,jk}^W(u,v)\}_{j,k\in[p]}\,, \quad
	(u,v) \in \cU^2\,, \; h\geq 0\,.
\end{equation}
Performing the spectral decomposition
\begin{equation}
	\label{n13}
	\widehat K_{jj}(u,v) = \sum_{l=1}^\infty \hat{\lambda}_{jl} \hat \psi_{jl}(u) \hat \psi_{jl}(v)\,,
\end{equation}
where $\hat \lambda_{j1} \ge \hat \lambda_{j2} \ge \cdots > 0$ are the eigenvalues, and $\hat \psi_{j1}(\cdot), \hat \psi_{j2}(\cdot),
\cdots$ are the corresponding eigenfunctions.
Let $\eE\{\bfeta_{(t-h)j}\bfeta_{tk}^{\T}\}=\{\sigma^{(h)}_{jklm}\}_{l \in [d_j], m \in [d_k]}$ with its estimator $(n-h)^{-1}\sum_{t=h+1}^n \widehat\bfeta_{(t-h)j} \widehat\bfeta_{tk}^{\T}=\{\hat\sigma^{(h)}_{jklm}\}_{l \in [d_j], m \in [d_k]}$ for $j,k \in [p]$ and $h \geq 0,$ where $\widehat \bfeta_{tj}=(\hat \eta_{tj1}, \dots,
\hat \eta_{tjd_j})^{\T}.$
Our proposed autocovariance-based Step~2 and Step~3 explicitly rely on the sample autocovariance among estimated basis coefficients, $\{\hat \sigma_{jklm}^{(h)}: j,k \in [p], l \in [d_j], m \in [d_k], h \in [L]\}$, and the estimated basis functions $\{\hat \psi_{jl}(\cdot): j \in [p], l \in [d_j]\},$ respectively. See details in Sections~\ref{sec.ill.ex} and ~\ref{sec_3app}. Their convergence properties in elementwise $\ell_{\infty}$-norm under high-dimensional scaling are investigated in Section~\ref{sec.dr.theory} below. %It is noteworthy that the penalized least squares approach commonly adopted in the covariance-based Step~2 \citep{kong2016} is inappropriate under model~(\ref{error.model}), because it is based on  %$\{\hat\sigma_{jklm}^{(0)}: j,k \in [p], l \in [d_j], m \in [d_k]\},$ the inconsistent estimator $\widehat \bSigma_{0}^W$ for $\bSigma_{0}^X.$

\subsection{Rates in elementwise $\ell_{\infty}$-norm}
\label{sec.dr.theory}
To characterize the effect of serial dependence on the relevant estimated terms, we will use the functional stability measure of $\{\bW_t(\cdot)\}_{t \in \mathbb{Z}}$ \citep{guo2022}.

\begin{condition}
	\label{cond_fsm}
	For $\{\bW_t(\cdot)\}_{t \in \mathbb{Z}}, $ the spectral density operator $\bbf_\theta^W=(2\pi)^{-1}\sum_{h\in \mathbb{Z}}\bSigma_h^We^{-ih\theta}$ for $\theta \in [-\pi,\pi]$ exists and
	the functional stability measure defined in {\rm(\ref{def.fsm})} is finite, i.e.
	\begin{equation}
		\label{def.fsm}
		\cM^W = 2\pi \cdot \esssup\limits_{\theta \in [-\pi,\pi],\bPhi\in \mathbb{H}_0^p}\frac{\langle\bPhi,\bbf_\theta^W(\bPhi)\rangle}{\langle\bPhi,\bSigma_0^W(\bPhi)\rangle} <\infty\,,
	\end{equation}
	where $\mathbb{H}_0^p = \{\bPhi \in \mathbb{H}^p:\langle\bPhi,\bSigma_0^W(\bPhi)\rangle \in (0, \infty)\}.$
\end{condition}

The quantity $\cM^W$ in (\ref{def.fsm}) is expressed proportional to functional Rayleigh quotients of $\bbf_{\theta}^W$ relative to $\bSigma_0^W.$ 
%Hence it can more precisely capture the effect of eigenvalues of $\bbf_{\theta}^W$ relative to small decaying eigenvalues of $\bSigma_0^W,$
Hence it can more precisely capture the effect of small decaying eigenvalues of $\bSigma_0^W$ on the numerator in (\ref{def.fsm}),
which is essential to handle truly infinite-dimensional functional objects $\{W_{tj}(\cdot)\}.$
%{\color{blue} To QY: 
	%The direct extension of Basu's stability measure is
	%$$
	%2\pi \cdot \esssup\limits_{\theta \in [-\pi,\pi],\bPhi\in \mathbb{H}_0^p}\frac{\langle\bPhi,\bbf_\theta^W(\bPhi)\rangle}{\langle\bPhi,\bPhi\rangle} <\infty\,,
	%$$
	%which assumes that the biggest eigenvalues of $\bbf_\theta^W$ to be bounded. However, such measure cannot handle infinite-dimensional functional data well as it does not measure the effect of small decaying eigenvalues. Instead, we make use of the generalized Rayleigh quotients in the sense that we control the biggest eigenvalues of $(\bSigma_0^W)^{-1}\bbf_\theta^W(\bPhi)$ to be bounded. Choosing $\Phi$ as eigenfunctions of $\bSigma_0^W,$ the denominator becomes eigenvalues, which decay to zero for infinite-dimensional functional data. The ratio then becomes infinite. However, it is intuitive that small eigenvalues of $\bSigma_0^W$ typically correspond with small eigenvalues of $\bbf_\theta^W,$ thus making ratio to be bounded very plausible. Such formulation can facilitate our technical derivation in particular for infinite-dimensional functional data. For i.i.d. data, the ratio degenerates to 1.}
%with $\sum_{l=1}^{\infty}\omega_{jl}<\infty$ for each $j.$
We next define the functional stability measure of all $k$-dimensional subsets of $\{\bW_t(\cdot)\}_{t \in \eZ},$ i.e. $\{(W_{tj}(\cdot):j \in J)^\T\}_{t \in \eZ}$ for $J\subset [p]$ with cardinality $|J| \leq k,$ by
\begin{equation} \label{def_sub_fsm}
	\cM_k^W = 2\pi \cdot \underset{\theta\in [-\pi, \pi],\|\bPhi\|_0 \le k,\bPhi \in \cH_0^p}{\text{ess}\sup} \frac{\langle \bPhi, \bbf_{\theta}^W(\bPhi)\rangle}{\langle \bPhi, \bSigma_0^W(\bPhi)\rangle}\,, ~~k\in[p]\,.
\end{equation} Under Condition~\ref{cond_fsm}, it is easy to verify that $\cM_k^W \leq \cM^W <\infty$. % which will be used in subsequent  analysis. %See also \cite{guo2022} for further discussions.

%Provided that 
Our non-asymptotic results are developed using the infinite-dimensional analog of Hanson--Wright inequality \citep{rudelson2013} in a general Hilbert space $\mathbb{H}$, for which we need to impose the sub-Gaussian condition.

\begin{definition}
	\label{def_subG}
	Let $Z_t(\cdot)$ be a mean zero random variable in $\mathbb{H}$ for any fixed $t$ and $\Sigma_0: \mathbb{H} \to \mathbb{H}$ be a covariance operator. Then $Z_t(\cdot)$ is a sub-Gaussian process if there exists a constant $c>0$ such that $\mathbb{E}(e^{\langle x, Z \rangle}) \leq e^{c^2\langle x,\Sigma_0(x)\rangle/2}$ for all $ x \in \mathbb{H}$.
\end{definition}

\begin{condition}
	\label{cond_flp}
	{\rm(i)}
	$\{\bW_t(\cdot)\}_{t \in \eZ}$ is a sequence of multivariate functional linear processes with sub-Gaussian errors, namely sub-Gaussian functional linear processes,
	$\bW_t(\cdot) = \sum_{l=0}^{\infty}\bB_l(\bvarepsilon_{t-l})$ for any $t \in \mathbb{Z}$,
	where $\bB_l=(B_{l,jk})_{p\times p}$ with each $B_{l,jk} \in \mathbb{S},$
	$\bvarepsilon_t(\cdot)=\{\varepsilon_{t1}(\cdot),\dots,\varepsilon_{tp}(\cdot)\}^{\T} \in \cH^p$ and the components in $\{\bvarepsilon_t(\cdot)\}_{t \in \eZ}$ are independent sub-Gaussian processes satisfying Definition~{\rm\ref{def_subG}}; {\rm(ii)}
	The coefficient functions satisfy $\sum_{l=0}^{\infty}\|\bB_l\|_{\infty} =O(1);$
	{\rm(iii)}
	$\omega_0^{\varepsilon} = \max_{j\in[p]}\int_\cU \Sigma_{0,jj}^{\varepsilon}(u,u)\,\bd u = O(1)$, where $\Sigma_{0,jj}^{\varepsilon}(u,u)={\rm Cov}\{\varepsilon_{tj}(u),\varepsilon_{tj}(u)\}$.
\end{condition}

The multivariate functional linear process can be seen as the generalization of functional linear process \citep{Bbosq1} to the multivariate setting and also the extension of multivariate linear process \citep{hamilton1994} to the functional domain.  Condition~\ref{cond_flp}(ii) ensures the stationarity of $\{\bW_t(\cdot)\}_{t \in \eZ}$ and, together with Condition~\ref{cond_flp}(iii), implies that  $\omega_0^{W} = \max_{j\in[p]}\int_\cU \Sigma_{0,jj}^{W}(u,u)\,\bd u =O(1)$ (see Lemma~\ref{lemma_eigflp} in Appendix~\ref{ap_proof}), which is essential in deriving non-asymptotic results.  
%In general, we can assume that %$\sum_{l=0}^{\infty}\|\bB_l\|_{\infty}$ and %$\omega_0^{\varepsilon}$ grow at very slow rates as $p$ %increases, then our established deviation bounds will depend %on these terms. 
\textcolor{black}{The sub-Gaussian condition is imposed on the functional process to facilitate the use of Hanson--Wright-type inequality in our non-asymptotic analysis. We believe that a Nagaev-type concentration bound can be established to accommodate functional linear process with functional errors under a weaker finite polynomial moments condition. It is also interesting to develop non-asymptotic results for more general non-Gaussian functional time series under other commonly adopted dependence framework.}

%\begin{condition} \label{con_sub_coefficient}
%The coefficient functions satisfy
% for any $j=1,\dots,p,$
%$\sum_{l=0}^{\infty}\|\bB_l\|_{\infty} < \infty.$
%= \max_{1\leq j \leq p}\sum_{l=0}^\infty \|\bA_{l}\|_{\infty} < \infty.$
%\end{condition}
%\textcolor{black}{Rosenblatt is on the stationarity of time series, what about the result under the functional domain? Can this condition be weaker?(definition of the operator norm? too much) implied from the proof of Lemma~\ref{lm_sigma_sub}}

%\begin{condition} \label{con_e}
%(i) The marginal-covariance functions of $\{\bvarepsilon_t(\cdot)\},$ $\Sigma_{0,jj}^{\varepsilon}(u,v)$'s, are continuous on $\cU^2$ and uniformly bounded over $j \in \{1,\dots,p\};$
%(ii) $\omega_0^{\varepsilon} = \max_{1 \leq j \leq p}\int_\cU \Sigma_{0,jj}^{\varepsilon}(u,u)du < \infty .$
%\end{condition}

\begin{condition}
	\label{cond_eigen}
	{\rm(i)}
	For each $j \in [p]$,  $\lambda_{j1}> \lambda_{j2} > \cdots >0,$ and there exist some positive constants $c_0$ and $\alpha>1$ such that $ \lambda_{jl}-\lambda_{j(l+1)} \geq c_0 l^{-\alpha-1}$ for $l \geq 1;$ {\rm(ii)}
	For each $j \in [p],$ the linear space spanned by  $\{\nu_{jl}(\cdot)\}_{l=1}^\infty$ (i.e. eigenfunctions of $\Sigma_{0,jj}^X$) is the same as that spanned by $\{\psi_{jl}(\cdot)\}_{l=1}^\infty.$
\end{condition}

Condition~\ref{cond_eigen}(i) controls the lower bound of eigengaps with larger values of $\alpha$ yielding tighter gaps between adjacent eigenvalues. %and also implies that $\lambda_{jl} \geq c_0 \alpha^{-1} l^{-\alpha}.$ %for some positive constant $c$. 
See similar conditions in \cite{hall2007} and \cite{kong2016}. To simplify notation, we assume the same $\alpha$ across $j,$ but this condition can be relaxed by allowing $\alpha$ to depend on $j$ and our theoretical results can be generalized accordingly.

We next establish the deviation bounds on estimated eigenpairs, $\{\hat{\lambda}_{jl}, \hat{\psi}_{jl}(\cdot)\},$ and the sample autocovariance among estimated basis coefficients, $\{\hat{\sigma}_{jklm}^{(h)}\}$, in elementwise $\ell_{\infty}$-norm.
%, which play a crucial role in further convergence analysis under high-dimensional scaling.
\begin{theorem}
	%\footnote{\color{black}Need to specify the relationship between $d$ and $n$. Same as Theorem 2.}
	\label{thm_eigen}
	Let Conditions~{\rm\ref{cond_fsm}--\ref{cond_eigen}} hold, and $d$ be a positive integer possibly depending on $(n,p)$. %If $\log(pd)d^{4\alpha+2}(\cM_{1}^Z)^2/n \rightarrow 0$ as $n, p \rightarrow \infty,$
	For $n \gtrsim \log p$, there exist some positive constants $c_1$ and $c_2$ independent of $(n,p,d)$ such that
	\begin{equation}
		\label{bd_max_eigen}
		\max_{j \in [p],l \in [d]} \bigg\{|\hat \lambda_{jl} - \lambda_{jl}|
		+ \bigg\|\frac{ \hat \psi_{jl} - \psi_{jl}}{l^{\alpha + 1}}\bigg\| \bigg\}
		\lesssim \cM_1^W \sqrt{\frac{\log p}{n}}
	\end{equation}
	holds with probability greater than
	$ 1- c_1 p^{-c_2}$, where $\cM_1^W$ is defined in \eqref{def_sub_fsm}.
\end{theorem}

\begin{theorem}
	\label{thm_score}
	Let conditions in Theorem~{\rm\ref{thm_eigen}} hold and $h \geq 1$ be fixed.
	For $n \gtrsim d^{2\alpha+2} (\cM_1^W)^2\log p,$  there exist some positive constants $c_3$ and $c_4$  independent of $(n,p,d)$ such that
	\begin{equation}
		\label{bd_score_max}
		\underset{j,k \in [p], l,m \in [d]}{\max} \frac{|\hat \sigma_{jklm}^{(h)} - \sigma_{jklm}^{(h)}|}{(l \vee m)^{\alpha+1}}\lesssim \cM_1^W \sqrt{\frac{\log p}{n}}
	\end{equation}
	holds with probability greater than $ 1- c_3p^{-c_4}$, where $\cM_1^W$ is defined in \eqref{def_sub_fsm}.
\end{theorem}

\begin{remark}
	\label{rmk_thm1}
	(i) The parameter $d$ in Theorems~{\rm\ref{thm_eigen}} and {\rm\ref{thm_score}} can be understood as the truncated dimension of infinite-dimensional functional objects under the expansion in {\rm(\ref{expansion})}. In general, $d$ can depend on $j,$ say $d_j,$ then the maximums in {\rm(\ref{bd_max_eigen})} and {\rm(\ref{bd_score_max})} are taken over $j,k \in [p],$ $l \in [d_j], m \in [d_k]$ and the corresponding right-sides remain the same.\\
	\textcolor{black}{(ii) Compared with the normalized deviation bounds %on estimated eigenpairs 
		under FPCA framework %, $\{(\hat \omega_{jl},\hat \nu_{jl})\},$ 
		%and sample autocovariance among estimated FPC scores 
		established in \cite{guo2022}, we obtain slower rates in (\ref{bd_max_eigen}) and (\ref{bd_score_max}) for decaying eigenvalues. 
		%Intuitively, as opposed to the expansion of $X_{tj}$ through $\psi_{j1}, \psi_{j2}, \dots$ with correlated coefficients
		Note that $\{\nu_{jl}(\cdot)\}_{l \geq 1}$ provides the unique basis with respect to which $X_{tj}(\cdot)$ can be expressed as Karhunen--Lo\`eve expansion with uncorrelated coefficients. It gives the most rapidly convergent representation of $X_{tj}(\cdot)$ in the $L_2$ sense. By comparison, the expansion of $X_{tj}(\cdot)$ through $\{\psi_{jl}(\cdot)\}_{l \geq 1}$ in (\ref{expansion}) results in a suboptimal convergent representation with correlated coefficients.
		From a theoretical viewpoint, whether the rates in (\ref{bd_max_eigen}) and (\ref{bd_score_max}) are minimax optimal is of interest and requires further investigation.
	}
\end{remark}

\section{Block RMD estimation framework}
\label{sec.rmd}
Resulting from Step 1, the estimation of sparse function-valued parameters is transformed to the block sparse estimation of parameter vectors/matrices in Step~2. To identify these parameters, we choose $\{\widehat\bfeta_{(t-h)k}: h \in [L], k \in [p]\}$ as vector-valued instrumental variables and construct autocovariance-based moment equations, which is illustrated using an example of SFLR in Section~\ref{sec.ill.ex}. We then formulate a general block RMD estimation method in Section~\ref{sec.rmd.method} and study its theoretical properties in Section~\ref{sec.thm_BRMD}.

\subsection{An illustrative example}
\label{sec.ill.ex}
We illustrate via the high-dimensional SFLR:
\begin{equation}
	\label{model_sflr}
	Y_t = \sum_{j=1}^p \int_{\cU} X_{tj}(u)\beta_{0j}(u)\,{\rm d}u + \varepsilon_t\,, ~~t \in [n]\,,
\end{equation}
{where % $p$ functional covariates
	$\{X_{tj}(\cdot)\}_{t \in [n], j \in [p]}$ satisfy model~(\ref{error.model}),
	%and are independent, 
	$\{\varepsilon_t\}_{t \in [n]}$ are i.i.d. and mean-zero random errors, and $\{X_{tj}(\cdot)\}$ and $\{\varepsilon_t\}$ are independent.
	% $\{\varepsilon_t\}_{t \in [n]}.$ %, and $\{\beta_{0j}(\cdot)\}_{j \in [p]}$ are unknown functional coefficients. %{\color {red} check this sentence??? -- QY.}
	Given observations $\{(\bW_t(\cdot), Y_t)\}_{t \in [n]},$ our goal is to estimate $p$  functional coefficients $\bbeta_0(\cdot)=\{\beta_{01}(\cdot), \dots, \beta_{0p}(\cdot)\}^{\T}.$} To guarantee a feasible solution under high-dimensional scaling, we assume that $\bbeta_0(\cdot)$ is functional $s$-sparse, i.e. $s$ components in $\bbeta_{0}(\cdot)$ are nonzero with $s \ll p.$
%with support $S=\{j \in [p]: \|\beta_{0j}\| \neq 0\}$ and cardinality $s=|S|.$

Resulting from the truncated expansion of $X_{tj}(\cdot)$ via (\ref{expansion}) in Step~1, (\ref{model_sflr}) can be rewritten as
\begin{equation*}
	\label{rep_sflr}
	Y_t = \sum_{j=1}^p\bfeta_{tj}^{\T}\bbb_{0j} + r_t + \varepsilon_t\,,
\end{equation*}
where %$\bfeta_{tj}=(\eta_{tj1}, \dots, \eta_{tjd_j})^{\T},$
%$\bbb_{0j}=(b_{0j1}, \dots, b_{0jd_j})^{\T}$ with each $b_{0jl}=\langle \psi_{jl}, \beta_{0j} \rangle$ and
%$\bpsi_j=(\psi_{j1}, \dots, \psi_{jd_j})^{\T},$
$\bbb_{0j} = \int_{\cU}\bpsi_j(u)\beta_{0j}(u)\,{\rm d}u\in\mathbb{R}^{d_j}$ and
$r_t=\sum_{j=1}^p\sum_{l=d_j+1}^{\infty}\eta_{tjl}\langle \psi_{jl}, \beta_{0j} \rangle$ is the truncation error. %Let $\bbb_0 = (\bbb_{01}^{\T}, \dots, \bbb_{0p}^{\T})^{\T} \in {\eR}^{\sum_{j=1}^p d_j}.$ and
Given some prescribed positive integer $L$, in Step~2, we choose $\{\bfeta_{(t-h)k}: h \in [L], k \in [p]\}$ as vector-valued instrumental variables. %, which are assumed to be uncorrelated with the random error $\varepsilon_t$ in (\ref{rep_sflr}). %Within the framework of moment equations in (\ref{b_moment}), $q=pL$ and
Then $\bbb_0 = (\bbb_{01}^{\T}, \dots, \bbb_{0p}^{\T})^{\T} \in {\eR}^{\sum_{j=1}^p d_j}$ can be identified by the following moment equations:
\begin{equation}
	\label{moment_sflr}
	\eE\{\bfeta_{(t-h)k}\varepsilon_t\}=\bg_{hk}(\bbb_0) + \bR_{hk}= {\bf 0}\,, ~~k \in [p]\,, \,h \in [L]\,,
	%=E\Big\{\bfeta_{(t+h)k}Y_t-\sum_{j=1}^p\bfeta_{(t+h)k}\bfeta_{tj}^{\T}\bbb_{0j}-\bfeta_{(t+h)k}r_t\Big\}={\bf 0} \in {\eR}^{q_k}.
\end{equation}
where
$\bg_{hk}(\bbb_0) =  \eE\{\bfeta_{(t-h)k}Y_t\}-\sum_{j=1}^p \eE\{\bfeta_{(t-h)k}\bfeta_{tj}^{\T}\bbb_{0j}\}$ and the bias term $\bR_{hk}=-\eE\{\bfeta_{(t-h)k}r_t\}.$ %formed by the truncation error.

With $\{\widehat\bfeta_{tj}\}_{t \in [n], j \in [p]}$ and $\{\widehat\bpsi_j(\cdot)\}_{j \in [p]}$ obtained in Step~1, for any $\bbb=(\bbb_1^\T,\ldots,\bbb_p^\T)^\T\in\mathbb{R}^{\sum_{j=1}^pd_j}$, we define
\begin{equation}
	\label{emp_sflr}
	\widehat \bg_{hk}(\bbb)=  \frac{1}{n-h}\sum_{t=h+1}^n \widehat\bfeta_{(t-h)k}Y_t - \frac{1}{n-h}\sum_{t=h+1}^n \sum_{j=1}^p\widehat\bfeta_{(t-h)k}\widehat\bfeta_{tj}^{\T}\bbb_{j}\,, ~~k \in [p]\,, \,h \in [L]\,,
\end{equation}
which provides the empirical version of $\bg_{hk}(\bbb)=  \eE\{\bfeta_{(t-h)k}Y_t\}-\sum_{j=1}^p \eE\{\bfeta_{(t-h)k}\bfeta_{tj}^{\T}\bbb_{j}\}$. %It follows from \eqref{moment_sflr} that
%\begin{align}\label{eq:1}
%\widehat{\bg}_{hk}(\bbb_0)\approx{\bf 0}\,,~~k \in [p]\,, \,h \in [L]\,.
%\end{align}
%Based on \eqref{eq:1}, 
Applying the block RMD estimation introduced in Section \ref{sec.rmd.method} below results in a block sparse estimator $\widehat \bbb = (\widehat \bbb_1^{\T}, \dots, \widehat \bbb_p^{\T})^{\T}$.

\subsection{A general estimation procedure}
\label{sec.rmd.method}
%Performing the autocovariance-based dimension reduction in the first step transforms the problem of modelling of $\bX_t(\cdot)$ into that of modelling of $(\sum_{j=1}^p d_j)$-dimensional time series, $\widehat\bfeta_t=(\widehat\eta_{t11}, \dots, \widehat\eta_{t1d_1}, \dots, \widehat\eta_{tp1}, \dots, \widehat\eta_{tpd_p})^{\T}.$
%The second step involves the estimation of block sparse vector- or matrix-valued parameters under high-dimensional scaling. The classical covariance-based least squares approach relies on $\{\widehat\sigma_{jklm}^{(0)}\},$ whose components are not consistent estimators for those in $\{\sigma_{jklm}^{(0)}\},$ and hence is inappropriate. Inspired from the fact that, for each $(j,k,l,m),$ $\widehat\sigma_{jklm}^{(h)}$ is a consistent estimator for $\sigma_{jklm}^{(h)}$ when $h \neq 0,$ we formulate the minimum distance estimation problem within an high-dimensional instrumental variables framework based on autocovariance among scores $\{\sigma_{jklm}^{(h)}\}$ in (\ref{autoscore}). Finally, in the third step, the functional sparse estimates are recovered from those block sparse estimates obtained in the second step.

In this section, we present the proposed Step~2 in a general block RMD estimation framework. %Motivated from the autocovariance-based moment condition, we propose a block regularized minimum distance (RMD) estimator for the block sparse coefficient.
%It follows from the dimension reduction step that the estimation of function-valued parameters involved in sparse high-dimensional functional models  
%As mentioned in Section~\ref{sec.auto3}, 
Note that Step~2 considers
the block sparse estimation of some %vector- or 
matrix-valued parameters, $\btheta_0 = (\btheta_{01}^{\T}, \dots, \btheta_{0p}^{\T})^{\T} \in \eR^{\sum_{j=1}^p d_j \times \tilde d}$ with each $\btheta_{0j}\in\mathbb{R}^{d_j\times \tilde{d}}.$
For SFLR with a scalar response, $\tilde d=1.$ 
%For FFLR and VFAR, $\tilde d \geq 1$ is the truncated dimension of the functional response.
Given some prescribed positive integer $L$ and $q=pL$ target moment functions $\btheta \mapsto \bg_i(\btheta)$ mapping $\btheta \in \eR^{\sum_{j=1}^p d_j \times \tilde d}$ to $\bg_i(\btheta)\in \eR^{d_k \times \tilde d}$ %for $i=L(j-1)+1, \dots, Lj$ with $j \in [p],$ %with their empirical versions given by $\widehat\bg_i : \eR^{\sum_{j=1}^p d_j \times \widetilde d} \mapsto \eR^{d_i \times \widetilde d},$
with $i=(h-1)p+k$ and $k \in [p]$ for $h \in [L],$
where both $p$ and $q$ are large, we assume that $\btheta_0$ can be identified by the following moment equations:
\begin{equation}
	\label{b_moment}
	\bg_i(\btheta_0) + \bR_i = {\bf 0}\,, ~~ i \in [q]\,,
\end{equation}
where $\bR_i$'s are formed by autocovariance-based truncation errors due to finite approximations in Step~1. %For SFLR, $\bR_{i}=-E(\bfeta_{(t+h)k}r_t)$ according to (\ref{moment_sflr1}).
We are interested in estimating the block sparse $\btheta_0$ based on empirical mappings $\btheta\mapsto\widehat\bg_i(\btheta)$ of $\btheta\mapsto\bg_i(\btheta)$ for $i \in [q].$ See Sections \ref{sec.ill.ex} and \ref{sec_3app} for detailed expressions of $\bg_i(\cdot)$ and $\widehat\bg_i(\cdot)$ in some exemplified models.

It follows from (\ref{b_moment}) that 
\begin{equation}
	\label{b_moment.app}
	\widehat \bg_i(\btheta_0) \approx {\bf 0}\,,~~i \in [q]\,.
\end{equation}
Based on (\ref{b_moment.app}),
we define the block RMD estimator $\widehat \btheta =(\widehat \btheta_1^{\T}, \dots, \widehat \btheta_p^{\T})^{\T} \in \eR^{{\sum_{j=1}^p d_j \times \tilde d}}$ as a solution to the following convex optimization problem:
\begin{equation}
	\label{opt_BRMD}
	\widehat \btheta = \underset{\btheta}{\arg\min} \sum_{j=1}^p \|\btheta_j\|_{\tF} ~~ \text{subject to}~~\underset{i \in [q]}{\max} \|\widehat \bg_i(\btheta)\|_{\tF} \leq \gamma_n\,,
\end{equation}
where $\gamma_n \geq 0$ is a regularization parameter. % For SFLR or FFLR with $\tilde d=1,$ the matrix Frobenius norm in (\ref{opt_BRMD}) degenerates to the vector $\ell_2$-norm. For FFLR and VFAR with $\tilde d>1,$ the corresponding optimization tasks are formulated under the matrix Frobenius norm. 
The group information is encoded in the objective function, which forces the elements of $\widehat \btheta_j$ to either all be zero or nonzero, thus producing the block sparsity in $\widehat\btheta.$
It is worth noting that, without the bias terms $\bR_i$'s in (\ref{b_moment}), our proposed block RMD estimation framework can be seen as a blockwise generalization of the RMD estimation \citep{belloni2018} by replacing $|\cdot|$ by  $\|\cdot\|_{\tF}$.
To solve the large-scale convex optimization problem in (\ref{opt_BRMD}), we use the R package \verb"CVXR" \citep{fu2020}, which is easy to implement and converges fast.  In Sections~\ref{sec_sflr}, \ref{sec_fflr} and \ref{sec_vfar}, we will illustrate our proposed autocovariance-based block RMD estimation framework using examples of SFLR, FFLR and VFAR, respectively.
%, in the context of high-dimensional functional time series.

\subsection{Theoretical properties}
\label{sec.thm_BRMD}
%We begin with some notation that will be used in this section. 
For a block matrix $\bB = (\bB_{ij})_{i \in [N_1], j \in [N_2]} \in \eR^{N_1 m_1 \times N_2 m_2}$ with the $(i,j)$-th block $\bB_{ij} \in {\eR}^{m_1 \times m_2},$  let $\|\bB\|_{\max}^{(m_1, m_2)} = \max_{i\in[N_1],j\in[N_2]} \|\bB_{ij}\|_{\tF},$ 
and  $\|\bB\|_{1}^{(m_1,m_2)} = \sum_{i=1}^{N_1} \|\bB_{i}\|_{\tF}$ when
$N_2=1$. To simplify notation in this section and theoretical analysis in Section~\ref{sec_3app}, we assume the same truncated dimension $d_{j}=d$ across $j \in [p],$ but our theoretical results can be extended naturally to the more general setting where $d_{j}$'s are different.

%Note $\bg_i: \eR^{pd \times \widetilde d} \mapsto \eR^{d \widetilde d},$ we
Let $\bg(\btheta) = \{\bg_1(\btheta)^{\T}, \dots, \bg_q(\btheta)^{\T}\}^{\T}$ and $\bR =(\bR_1^{\T}, \dots, \bR_q^{\T})^{\T} \in \eR^{qd \times \tilde d}.$ We focus on the case of which the moment function $\btheta \mapsto \bg(\btheta)$ mapping from $\eR^{pd \times \tilde d}$ to $\eR^{qd \times \tilde d}$ is linear with respect to $\btheta$ in the form of $\bg(\btheta)=\bG\btheta + \bg({\bf 0})$ for some $\bG \in \eR^{qd \times pd}.$ This together with
(\ref{b_moment}) implies that
\begin{equation}
	\label{linear.g}
	\bG\btheta_0 + \bg({\bf 0}) +\bR = {\bf 0}\,,
\end{equation}
the form of which can be easily verified for, e.g., SFLR, FFLR and VFAR models considered in Section \ref{sec_3app}.
%\begin{equation}
%\label{linear.g}
%{\bf 0}=\bg(\btheta) +\bR = \bG\btheta + \bg({\bf 0}) +\bR\,.
%\end{equation}
%Here $\bg(\btheta)=\bG\btheta + \bg({\bf 0})$ is linear with respect to $\btheta$ and $\bR$ is the bias term.
%It is easy to verify that SFLR, FFLR and VFAR can be represented in the form of (\ref{linear.g}).
Now we reformulate the optimization task in (\ref{opt_BRMD}) as
\begin{equation}
	\label{opt_BRMD1}
	\widehat \btheta = \underset{\btheta}{\arg\min} \|\btheta\|_1^{(d,\tilde d)} ~~ \text{subject to}~~ \|\widehat \bg(\btheta)\|_{\max}^{(d,\tilde d)} \leq \gamma_n\,,
\end{equation}
where
$\widehat\bg(\btheta)= \widehat \bG \btheta + \widehat \bg({\bf 0})$ is the empirical version of $\bg(\btheta).$
It is worth noting that $\btheta_0$ is block $s$-sparse with support $S=\{j \in [p]: \|\btheta_{0j}\|_{\tF} \neq 0\}$ and its cardinality $s=|S|.$

Before presenting properties of the block RMD estimator $\widehat \btheta$, we impose some high-level regularity conditions.

\begin{condition}
	\label{cond_emc}
	(i) There exist $\epsilon_{n1}, \delta_{n1}>0$ such that
	$
	\|\widehat \bG - \bG\|_{\max}^{(d,d)} \vee \|\widehat \bg({\bf 0}) - \bg({\bf 0})\|_{\max}^{(d, \tilde d)} \leq \epsilon_{n1}$ with probability at least $1-\delta_{n1};$
	%\end{condition}
	%\begin{condition}
	%    \label{cond_bd.R} 
	(ii) There exists $\epsilon_{2} > 0$ such that
	$\|\bR \|_{\max}^{(d, \tilde d)} \leq \epsilon_{2};$
	%\end{condition}
	%
	%\begin{condition}
	%\label{cond_reg_par}
	(iii) There exists $\delta_{n2}>0$ such that
	$\|\widehat \bg (\btheta_0)\|_{\max}^{(d,\tilde d)} \leq \gamma_n$ with probability at least $1-\delta_{n2}.$
\end{condition}

Conditions~\ref{cond_emc}(i) and (ii) together ensure that the empirical moment functions are nicely concentrated around the target moment functions. Using our derived non-asymptotic results in Section~\ref{sec.dr.theory}, we can easily specify the concentration bounds in Condition~\ref{cond_emc}(i) for SFLR, FFLR and VFAR. With further imposed smoothness conditions on coefficient functions, Condition~\ref{cond_emc}(ii) can also be verified. Condition~\ref{cond_emc}(iii) indicates that $\btheta_0$ is feasible in the optimization problem (\ref{opt_BRMD1}) with high probability, in which case a solution $\widehat\btheta$ of (\ref{opt_BRMD1}) exists and satisfies $\|\widehat\btheta\|_1^{(d,\tilde d)} \leq \|\btheta_0\|_1^{(d,\tilde d)}.$ The non-block version of such property typically plays a crucial role to tackle high-dimensional models in the literature.

Let $\bdelta = \btheta - \btheta_0.$ We define a block $\ell_1$-sensitivity coefficient
\begin{equation}
	\label{def.kappa}
	\kappa(\btheta_0) = \inf_{T:\, |T| \leq s} \inf_{\bdelta \in C_T:\, \| \bdelta \|^{(d, \tilde d)}_1 > 0} \frac{ \| \bG \bdelta \|^{(d, \tilde d)}_{\max} }{\| \bdelta \|^{(d, \tilde d)}_{1}}  \,,
\end{equation}
where $ C_T=\{ \bdelta \in \eR^{pd \times \tilde d}: \|\bdelta_{T^{\rm c}}\|^{(d, \tilde d)}_1 \leq \|\bdelta_{T}\|^{(d, \tilde d)}_1\}$ for $T \subset [p].$ Provided that $\widehat \bdelta=\widehat\btheta-\btheta_0 \in C_S$ under Condition~\ref{cond_emc}(iii) as justified in Lemma~\ref{lemma_delta_S_bd} in Appendix~\ref{ap_proof}, the lower bound of $\kappa(\btheta_0)$ is useful to establish the error bound for $\|\widehat\bdelta\|_1^{(d, \tilde d)}.$ See also %\cite{belloni2018} and 
\cite{gautier2019} for non-block $\ell_q$-sensitivity quantities to handle high-dimensional instruments. We then need Condition~\ref{cond_svd_bd} below to determine such lower bound. 
%Before presenting this condition, we introduce some notation. Let 
%For $J \subset [q]$ and $M \subset [p],$ let $\bG_{J,M}=(\bG_{jk})_{j \in J, k \in M}$ with each $\bG_{jk} \in \eR^{d\times d}$ being the block submatrix of $\bG$ consisting of all block rows $j \in J$ and all block columns $k \in M$ of $\bG.$ 
Note that $\bG$ can be divided into $q\times p$ blocks of the size $d\times d$. Let  $\bG_{J,M}$ be the submatrix of $\bG$ consisting of all the $(j, k)$-blocks with
$j\in J\subset [q]$ and $k \in M \subset[p].$
For an integer $m \geq s,$ let \begin{equation*}
	\sigma_{\min}(m, \bG)=\min_{|M| \le m} \max_{|J| \le m} \sigma_{\min}(\bG_{J,M})~~\text{and}~~
	\sigma_{\max}(m, \bG)=\max_{|M| \le m} \max_{|J| \le m} \sigma_{\max}(\bG_{J,M})\,,
\end{equation*}
where $\sigma_{\min}(\bG_{J,M})$ and $\sigma_{\max}(\bG_{J,M})$ are the smallest and largest singular values of $\bG_{J,M}.$

\begin{condition}
	\label{cond_svd_bd} There exist universal constants $c_5>0$ and $\mu >0$ such that  $\sigma_{\max}(m, \bG)\geq c_5$ and $\sigma_{\min}(m, \bG)/\sigma_{\max}(m, \bG) \geq \mu$ for $m=16s/\mu^2.$
\end{condition}

In Condition~\ref{cond_svd_bd}, the quantity $\mu$ serves as a key factor to determine the lower bound of $\kappa(\btheta_0),$ which is justified in Lemma~\ref{lemma_kappa_bd2} in Appendix~\ref{ap_proof}. When $\mu$ is bounded away from zero, we have a strongly-identified model. When $\mu \rightarrow 0,$ it corresponds to the scenario with weak instruments. See also \cite{belloni2018} for similar conditions. %We are now ready to present the theorem on the convergence rate of $\widehat\btheta.$

\begin{theorem}
	\label{thm_BRMD}
	Let Conditions~{\rm\ref{cond_emc}--\ref{cond_svd_bd}} hold. If $\|\btheta_0\|_{1}^{(d,\tilde d)} \leq K$ for some $K >0$ and the regularization parameter $\gamma_n \lesssim (K+1) \epsilon_{n1} + \epsilon_2,$ then with probability at least $1-(\delta_{n1}+\delta_{n2}),$ the block RMD estimator $\widehat\btheta$ satisfies
	\begin{equation}
		\label{rate_BRMD}
		\|\widehat\btheta -\btheta_0\|_{1}^{(d,\tilde d)}  \lesssim s\mu^{-2}\{(K +1)\epsilon_{n1}  + \epsilon_2\}\,.
	\end{equation}
\end{theorem}

\begin{remark}
	\label{rmk_thm2}
	(i) The error bound in (\ref{rate_BRMD}) has the familiar variance-bias tradeoff as commonly considered in nonparametrics statistics, suggesting us to carefully select the truncated dimension $d$ so as to balance the variance and bias terms for the optimal estimation.
	
	(ii) With commonly imposed smoothness conditions on functional coefficients, it is easy to verify that $K \vee \epsilon_2=o(s)$ for SFLR, FFLR and VFAR in Section \ref{sec_3app}.
	
	(iii) For three examples we consider, $\bG$ is formed by $\{\sigma_{jklm}^{(h)}: j,k \in [p], l,m \in [d], h \in [L]\}$ with the
	components $\sigma_{jklm}^{(h)}$ satisfying $|\sigma_{jklm}^{(h)}|\leq \{\eE(\eta_{(t-h)jl}^2)\}^{1/2}[\eE\{\eta_{tkm}^2\}]^{1/2}=\lambda_{jl}^{1/2}\lambda_{km}^{1/2} \rightarrow 0$ as $l,m\rightarrow\infty.$ Consider a general cross-covariance matrix $\bG= \mathbb{E}(\bx\by^{\T}) \in \eR^{qd \times pd}$ with entries decaying to zero as $d\rightarrow\infty$, where $\bx=(x_1, \dots, x_{qd})^{\T}$ with $\eE(\bx)={\bf 0}$ and $\by =(y_1, \dots, y_{pd})^{\T}$ with $\eE(\by)={\bf 0},$ it is more sensible to impose Condition~\ref{cond_svd_bd} on its normalized version
	$\widetilde \bG  = \bD_x \bG \bD_y$ instead of $\bG$ itself,
	where $\bD_{x}={\rm{diag}}\{{\rm Var}(x_1)^{-1/2}, \dots, {\rm Var}(x_{qd})^{-1/2}\}$ and $\bD_{y}={\rm{diag}}\{{\rm{Var}}(y_1)^{-1/2}, \dots, {\rm Var}(y_{pd})^{-1/2}\}.$ For three exemplified models, $\bD_x$ and $\bD_y$ are formed by $\{\lambda_{jl}^{-1/2}: j \in [p], l \in [d]\}.$
\end{remark}

Remark~\ref{rmk_thm2}(iii) motivates us to present the following proposition that will be used in the theoretical analysis of associated estimators for SFLR, FFLR and VFAR in Section~\ref{sec_3app}.

\begin{proposition}
	\label{prop_BRMD}
	Suppose that all conditions in Theorem~{\rm\ref{thm_BRMD}} hold except that Condition~{\rm\ref{cond_svd_bd}} holds for $\widetilde \bG,$ then with probability at least $1-(\delta_{n1}+\delta_{n2}),$ the block RMD estimator $\widehat\btheta$ satisfies
	\begin{equation}
		\label{rate_BRMD1}
		\|\widehat\btheta -\btheta_0\|_{1}^{(d,\tilde d)}  \lesssim  s\mu^{-2}\|\bD_x\|_{\max}\|\bD_y\|_{\max}\{(K +1)\epsilon_{n1}  + \epsilon_2\}\,.
	\end{equation}
\end{proposition}

\section{Applications}
\label{sec_3app}
%In Section~\ref{sec_3app}, we provide theoretical analysis of the general block RMD estimator, based on which we then investigate convergence properties of the autocovariance-based estimates for SFLR, FFLR and VFAR in Section~\ref{sec.thm.autoest}.
In this section, we illustrate the proposed 
%autocovariance-based 
estimation procedures with the three concrete models, namely
%the corresponding convergence analysis using applications of 
SFLR, FFLR and VFAR. %models %under high-dimensional scaling
%in Sections~\ref{sec_sflr}, \ref{sec_fflr} and \ref{sec_vfar}, respectively.
%and derived theoretical results to study the convergence properties for high-dimensional SFLR, FFLR and VFAR models in Sections~\ref{sec_sflr}, \ref{sec_fflr} and \ref{sec_vfar}, respectively.
%to estimate high-dimensional SFLR, FFLR and VFAR models and present the corresponding convergence properties in Sections~\ref{sec_sflr}, \ref{sec_fflr} and \ref{sec_vfar}, respectively.

\subsection{High-dimensional SFLR}
\label{sec_sflr}
Consider the high-dimensional SFLR in (\ref{model_sflr}), we first perform autocovariance-based dimension reduction on $\{W_{tj}(\cdot)\}_{t\in [n]}$ for each $j \in [p].$ 
According to Section~\ref{sec.ill.ex} and following the optimization framework in (\ref{opt_BRMD}), we then develop the block RMD estimator $\widehat \bbb$ as a solution to the constrained optimization problem:
\begin{equation*}
	\label{opt_sflr}
	\widehat \bbb = \underset{\bbb}{\arg\min} \sum_{j=1}^p \|\bbb_j\|_2 ~~ \text{subject to}~~\underset{k \in [p], h \in [L]}{\max} \|\widehat \bg_{hk}(\bbb)\|_2 \leq \gamma_n\,,
\end{equation*}
where $\gamma_n \geq 0$ is a regularization parameter and $\widehat\bg_{hk}(\bbb)$ is defined in (\ref{emp_sflr}).
Given that the recovery of functional sparsity in $\bbeta_{0}(\cdot)$ is equivalent to estimating the block sparsity in $\bbb_{0},$ in Step~3, we estimate functional sparse coefficients by
\begin{equation}
	\label{beta_est_sflr}
	\hat\beta_j(\cdot)=\widehat\bpsi_j(\cdot)^{\T}\widehat\bbb_j\,,~~ j\in [p]\,.
\end{equation}
%Finally, we obtain estimated functional coefficients $\{\hat\beta_j(\cdot)\}_{j \in [p]}$ in (\ref{beta_est_sflr}).
%Finally, our estimated functional coefficients are obtained by
%\begin{equation}
%\label{beta_est_sflr}
%\widehat\beta_j(\cdot)=\widehat\bpsi_j(\cdot)^{\T}\widehat\bbb_j ~~\text{for } j \in [p].
%\end{equation}

%Before stating these proofs, we give some notation. For a block matrix $\bB = (\bB_{jk}) \in \eR^{p_1d \times p_2 d}$, $p_1, p_2 \in \eN^+,$ with its $(j,k)$-th block $\bB_{jk} \in {\eR}^{d \times d},$  we define its $d$-block versions of Frobenius, elementwise $\ell_{\infty}$ and matrix $\ell_1$ norms by $\|\bB\|_{\tF} = \big(\sum_{j,k}\|\bB_{jk}\|_{\tF}^2\big)^{1/2},$  $\|\bB\|_{\max}^{(d)} = \max_{j,k} \|\bB_{jk}\|_{\tF}$ and $ \|\bB\|_{1}^{(d)} = {\max}_{k}\sum_{j} \|\bB_{jk}\|_{\tF},$ respectively.
%\subsubsection{SFLR}
%\label{sec.thm_sflr}
We next present the convergence analysis of $\{\hat\beta_j(\cdot)\}_{j \in [p]}.$ To simplify the notation, we assume the same truncated dimension $d_{j}=d$ across $j \in [p].$ We rewrite (\ref{moment_sflr}) in the form of (\ref{linear.g}), where
$\bg=(\bg_{11}^{\T}, \dots, \bg_{1p}^{\T}, \dots, \bg_{L1}^{\T}, \dots, \bg_{Lp}^{\T})^{\T},$
$\bR=(\bR_{11}^{\T}, \dots, \bR_{1p}^{\T}, \dots, \bR_{L1}^{\T}, \dots, \bR_{Lp}^{\T})^{\T}$ and
$\bG = (\bG_{ij}) \in \eR^{pLd \times pd}$ whose $(i,j)$-th block is $\bG_{ij}= \mathbb{E}\{\bfeta_{(t-h)k} \bfeta_{tj}^{\T}\} \in \eR^{d \times d}$ with $i=(h-1)p+k$ and $k \in [p]$ for $h \in [L].$ Applying Theorem~\ref{thm_score} and Proposition~\ref{prop_score_cross_mix} in Appendix~\ref{asec.fcsm} on $\widehat\bG$ and $\widehat\bg({\bf 0}),$ respectively, we can verify Condition~\ref{cond_emc}(i) with the choice of $\epsilon_{n1} \asymp \cM_{W,Y} d^{\alpha+2}(n^{-1}\log p)^{1/2},$ where $\cM_{W,Y}$ is specified in Proposition \ref{prop_score_cross_mix}. Before presenting the main theorem, we list the regularity conditions below.

\begin{condition}
	\label{cond_sflr}
	(i) For each $j \in S=\{j \in [p]: \|\beta_{0j}\| \neq 0\},$ $\beta_{0j}(\cdot) = \sum_{l=1}^\infty a_{jl}\psi_{jl}(\cdot)$ and there exists some positive constant $\tau > \alpha + 1/2$ such that $|a_{jl}| \lesssim l^{-\tau}$ for $l \geq 1;$
	%\end{condition}
	%
	%\begin{condition}
	%\label{cond_G_sflr}
	(ii) Let $\widetilde \bG=(\widetilde \bG_{ij})$ be the normalized version of $\bG=(\bG_{ij})$ by replacing each $\bG_{ij}$ by
	$\widetilde \bG_{ij}= \mathbb{E}\{\bD_{k}\bfeta_{(t-h)k} \bfeta_{tj}^{\T}\bD_{j}\},$ $i=(h-1)p+k,$ $k \in [p]$ for $h \in [L]$ and $j \in [p],$ where
	$\bD_j = {\rm{diag}}(\lambda_{j1}^{-1/2}, \dots, \lambda_{jd}^{-1/2}).$ There exist universal constants $c_{6}>0$ and $\mu >0$ such that  $\sigma_{\max}(m, \widetilde\bG)\geq c_{6}$ and $\sigma_{\min}(m, \widetilde \bG)/\sigma_{\max}(m, \widetilde \bG) \geq \mu$ for $m=16s/\mu^2.$
\end{condition}

Condition~\ref{cond_sflr}(i) restricts each component in $\{\beta_{0j}(\cdot): j \in S\}$ based on its expansion through basis $\{\psi_{jl}(\cdot)\}_{l \geq 1}$. The parameter $\tau$ determines the decay rate of basis coefficients and hence controls the level of smoothness with large values yielding smoother functions in $\{\beta_{0j}(\cdot): j \in S\}.$ %Condition~\ref{cond_coef_sflr} also ensures $\|\bbb_0\|_1^{(d,1)} \vee \|\bR\|_{\max}^{(d,1)}=O(s).$
See similar conditions in \cite{hall2007} and \cite{kong2016}.
Noting that components of $\bG$ decay to zero as $d$ grows to infinity, we impose Condition~\ref{cond_sflr}(ii) on $\widetilde \bG,$ which can be viewed as the normalized counterpart of Condition~\ref{cond_svd_bd} for SFLR.

Applying Proposition~\ref{prop_BRMD} and Theorem~\ref{thm_eigen} yields the convergence rate of the SFLR estimate $\widehat\bbeta(\cdot) = \{\hat\beta_{1}(\cdot), \dots, \hat\beta_p(\cdot)\}^{\T}$ under functional $\ell_1$ norm in the following theorem.

\begin{theorem}
	\label{thm_sflr}
	Suppose that Conditions~{\rm\ref{cond_fsm}--\ref{cond_eigen}}, {\rm\ref{cond_sflr}} and {\rm\ref{cond_cfsm}(ii)} in Appendix~\ref{asec.fcsm} hold, and $\{Y_t\}_{t \in [n]}$ is sub-Gaussian linear process. If the regularization parameter $\gamma_n \asymp  s \{ d^{\alpha+2}\cM_{W,Y}(n^{-1}\log p)^{1/2}$ $+ d^{-\tau+1/2}\},$ then the estimate $\widehat \bbeta(\cdot)$ satisfies
	\begin{equation}
		\label{rate_sflr}
		\sum_{j=1}^p\|\hat\beta_j -\beta_{0j}\| = O_\p\bigg\{ \mu^{-2} s^2 \bigg( d^{2\alpha+2}\cM_{W,Y} \sqrt{\frac{\log p}{n}} + d^{\alpha - \tau + 1/2}\bigg) \bigg\}\,.
	\end{equation}
\end{theorem}

\begin{remark}
	\label{rmk_thm4}
	%\begin{enumerate} [(a)]
	%\item
	(i) The rate of convergence in (\ref{rate_sflr}) is governed by both dimensionality parameters $(n, p, s)$ and internal parameters $(\cM_{W,Y}, d, \alpha, \tau, \mu).$ Typically, the rate is better when $\tau, \mu$ are large and $\cM_{W,Y}, \alpha$ are small.  To balance the variance and bias terms in (\ref{rate_sflr}) for the optimal estimation, we can choose the optimal truncated dimension $d \asymp (\cM^2_{W,Y}n^{-1}\log p)^{-1/(2\tau+2\alpha+3)}.$\\
	\textcolor{black}{
		(ii) Note that our convergence analysis relies on (\ref{bd_score_max}) rather than the normalized deviation bounds in \cite{guo2022}, the rate in (\ref{rate_sflr}) is slightly slower than that in \cite{fang2022} by a multiplicative factor $d^{\alpha/2}.$ For univariate functional linear regression, we similarly observe a slower rate for the autocovariance-based generalized methods-of-moments estimator \citep{chen2022} compared to the covariance-based least squares estimator \citep{hall2007}. %We leave the investigation on the optimality of (\ref{rate_sflr}) as a topic for future research.
	}
	%\end{enumerate}
\end{remark}

\subsection{High-dimensional FFLR}
\label{sec_fflr}

%\subsubsection{Estimation}
%\label{sec_est_fflr}
%In Sections~\ref{sec_est_fflr} and \ref{sec_est_vfar}, we illustrate the proposed autocovariance-based group RMD estimation framework using two examples of FFLR and VFAR.
Consider high-dimensional FFLR in the form of
\begin{equation}
	\label{model_fflr}
	Y_t(v) = \sum_{j=1}^p \int_{\cU} X_{tj}(u)\beta_{0j}(u,v)\,{\rm d}u + \varepsilon_t(v)\,, ~~t \in [n]\,, \, v \in \cal V\,,
\end{equation}
where $\{\bX_t(\cdot)\}_{t \in [n]}$ satisfy model~(\ref{error.model}) and are independent of i.i.d. mean-zero functional errors $\{\varepsilon_t(\cdot)\}_{t \in [n]},$ and $\{\beta_{0j}(\cdot,\cdot)\}_{j \in [p]}$ are functional coefficients to be estimated. With observed data $\{(\bW_t(u), Y_t(v)): (u,v) \in \cU \times {\cal V}, t \in [n]\},$ we target to estimate $\bbeta_0=\{\beta_{01}(\cdot,\cdot), \dots, \beta_{0p}(\cdot,\cdot)\}^{\T}$ under a functional sparsity constraint when $p$ is large. Specifically, we assume $\bbeta_0$ is functional $s$-sparse with support $S=\{j \in [p]: \|\beta_{0j}\|_{\cS} \neq 0\}$ and cardinality $s=|S| \ll p.$ 

Provided that each observed $Y_t(\cdot)$ is decomposed into the sum of dynamic and white noise components in (\ref{model_fflr}), we approximate $Y_t(\cdot)$ under the Karhunen--Lo\`eve expansion truncated at $\tilde d,$ i.e.
$Y_{t}(\cdot) \approx \bzeta_t^{\T}\bphi(\cdot),$
where  $\bzeta_t=(\zeta_{t1}, \dots, \zeta_{t\tilde d})^{\T}$ and $\bphi(\cdot) = \{\phi_1(\cdot), \dots, \phi_{\tilde d}(\cdot)\}^{\T}.$
Note that we can relax the independence assumption for $\{\varepsilon_t(\cdot)\}_{t \in [n]}$ and model observed responses via $\widetilde Y_t(\cdot) = Y_t(\cdot) + e_t^Y(\cdot),$ where $Y_t(\cdot)$ and $e_t^Y(\cdot)$ correspond to the dynamic signal and white noise elements, respectively. Then $Y_t(\cdot)$ can be approximated under the autocovariance-based
expansion in the sense of (\ref{expansion}) and our subsequent analysis still follow.

For each $j \in [p],$ we expand $X_{tj}(\cdot)$ according to (\ref{expansion}) truncated at $d_j.$ Some specific calculations lead to the representation of  (\ref{model_fflr}) as
\begin{equation}
	\label{rep_fflr}
	\bzeta_t^{\T} = \sum_{j=1}^p\bfeta_{tj}^{\T}\bB_{0j} + \br_t^{\T} + \bvarepsilon_t^{\T}\,,
\end{equation}
where $\bB_{0j}=\int_{\cU \times {\cal V}}\bpsi_{j}(u)\beta_{0j}(u,v)\bphi(v)^{\T}\,{\rm d}u{\rm d}v \in {\eR}^{d_j \times \tilde d}$ and $\br_{t}=(r_{t1}, \dots, r_{t\tilde d})^{\T}$ is the truncation error with $r_{tm}=\sum_{j=1}^p\sum_{l=d_j+1}^{\infty}\eta_{tjl}\langle\langle \psi_{jl}, \beta_{0j} \rangle, \phi_m \rangle$ for $m \in [\tilde d].$ Let $\bB_0 = (\bB_{01}^{\T}, \dots, \bB_{0p}^{\T})^{\T} \in {\eR}^{\sum_{j=1}^p d_j \times \tilde d}.$ We choose $\{\bfeta_{(t-h)k}: h \in [L], k \in [p]\}$ as vector-valued instrumental variables, which are assumed to be uncorrelated with the random error $\bvarepsilon_t$ in (\ref{rep_fflr}). Within the framework of (\ref{b_moment}), we %choose $q=pL$ and
assume that $\bB_0$ is the unique solution to the following moment equations:
\begin{equation}
	\label{moment_fflr}
	{\bf 0} = \eE\{\bfeta_{(t-h)k}\bvarepsilon_t^{\T}\}=\bg_{hk}(\bB_0) + \bR_{hk}\,, %\in {\eR}^{q_k \times \widetilde d}
	~~h \in [L]\,,~ k \in [p]\,,
	%E\Big\{\bfeta_{(t+h)k} \bzeta_t^{\T} - \sum_{j=1}^p \bfeta_{(t+h)k}\bfeta_{tj}^{\T}\bB_{0j}-\bfeta_{(t+h)k} \br_t^{\T}\Big\}
\end{equation}
where $\bg_{hk}(\bB_0) = \eE\{\bfeta_{(t-h)k}\bzeta_t^{\T}\}-\sum_{j=1}^p \eE\{\bfeta_{(t-h)k}\bfeta_{tj}^{\T}\bB_{0j}\}$ and $\bR_{hk}= -\eE\{\bfeta_{(t-h)k}\br_t^{\T}\}.$
%In a similar fashion to moment equations in (\ref{b_moment}), assume $\bB_0$ is the unique solution to
%\begin{equation}
%    \label{moment_fflr1}
%    g_{hk}(\bB_0) + \bR_{hk}= {\bf 0} ~~\text{for}~~k=1, \dots, p, h=1, \dots, L.
%\end{equation}
%where $L$ is some prescribed positive integer. See \cite{chen2022} for the selection of $L$ in practice.

Given the recovery equivalence between functional sparsity in $\bbeta_{0}$ and the block sparsity in $\bB_{0},$ we aim to estimate the block sparse matrix $\bB_0$ using the empirical versions $\bB \mapsto \widehat \bg_{hk}(\bB)$ for $h \in [L]$ and $k \in [p],$
\begin{equation*}
	\label{emp_fflr}
	\widehat \bg_{hk}(\bB)=  \frac{1}{n-h}\sum_{t=h+1}^n \widehat\bfeta_{(t-h)k}\widehat\bzeta_t^{\T} - \frac{1}{n-h}\sum_{t=h+1}^n \sum_{j=1}^p\widehat\bfeta_{(t-h)k}\widehat\bfeta_{tj}^{\T}\bB_{j}\,,
\end{equation*}
where $\widehat\bzeta_t=(\hat\zeta_{t1}, \dots, \hat\zeta_{t \tilde d})^{\T}$ with $\hat\zeta_{tm}=\langle Y_t, \hat \phi_m\rangle$ for $m \in [\tilde d]$ and $\{\widehat \bfeta_{tj}\}_{t\in [n], j \in [p]}$ are obtained in Step~1. In Step~2, according to (\ref{opt_BRMD}),  we formulate the block RMD estimator $\widehat \bB$ by solving the convex optimization problem below:
\begin{equation*}
	\label{opt_fflr}
	\widehat \bB = \underset{\bB}{\arg\min} \sum_{j=1}^p \|\bB_j\|_{\tF} ~~ \text{subject to}~~\underset{k\in [p], h \in [L]}{\max} \|\widehat \bg_{hk}(\bB)\|_{\tF} \leq \gamma_n\,,
\end{equation*}
where $\gamma_n \geq 0$ is a regularization parameter. In Step~3, we estimate the coefficient functions by
\begin{equation}
	\label{beta_est_fflr}
	\hat\beta_j(u,v)=\widehat\bpsi_j(u)^{\T}\widehat\bB_j\widehat\bphi(v)\,,~~ (u,v) \in \cU \times {\cal V}\,, j \in [p]\,,
\end{equation}
where $\{\widehat\bpsi_j(u)\}_{j \in [p]}$ and $\widehat\bphi(v) = \{\hat\phi_1(v), \dots, \hat\phi_{\tilde d}(v)\}^{\T}$ are obtained in Step~1.

%\subsubsection{FFLR}
%\label{sec.thm_fflr}
In the following, we investigate the convergence property of $\{\hat\beta_j(\cdot,\cdot)\}_{j \in [p]}$ in (\ref{beta_est_fflr}). To simplify the notation, we assume the same truncated dimension $d_{j}=d$ across $j \in [p]$.
We first rewrite (\ref{moment_fflr}) in the form of (\ref{linear.g})
and apply Theorem~\ref{thm_score} and Proposition~\ref{prop_score_cross} in Appendix~\ref{asec.fcsm} on $\widehat\bG$ and $\widehat\bg({\bf 0})$ to verify Condition~\ref{cond_emc}(i) with the choice of $\epsilon_{n1} \asymp \cM_{W,Y} d^{\alpha \vee \tilde \alpha + 2}(n^{-1}\log p)^{1/2}$, where $\cM_{W,Y}$ is specified in Proposition \ref{prop_score_cross}. In a similar fashion to $\alpha,$ the parameter $\tilde \alpha$ as specified in Condition~\ref{cond_Y_eigen} in Appendix~\ref{asec.fcsm} determines the tightness of eigengaps of the covariance function of $\{Y_t(\cdot)\}.$ We then impose the following smoothness condition on nonzero coefficient functions.

\begin{condition}
	\label{cond_coef_fflr}
	For each $j \in S,$ $\beta_{0j}(u,v) = \sum_{l,m=1}^\infty a_{jlm}\psi_{jl}(u)\phi_m(v)$ and there exists some positive constant $\tau > \alpha \vee \tilde \alpha + 1/2$ such that $|a_{jlm}| \lesssim (l+m)^{-\tau-1/2}$ for $l,m \geq 1.$
\end{condition}

We are now ready to present the convergence rate of the FFLR estimate $\widehat\bbeta(\cdot,\cdot)=\{\hat\beta_{1}(\cdot,\cdot), \dots, \hat\beta_p(\cdot,\cdot)\}^{\T}$ under functional $\ell_1$ norm in Theorem~\ref{thm_fflr}.

\begin{theorem}
	\label{thm_fflr}
	Suppose that Conditions~{\rm\ref{cond_fsm}--\ref{cond_eigen}}, {\rm \ref{cond_sflr}(ii), \ref{cond_coef_fflr}} and {\rm\ref{cond_cfsm}(i), \ref{cond_Y_eigen}} in Appendix~{\rm\ref{asec.fcsm}} hold, and $\{Y_t(\cdot)\}_{t\in[n]}$ is sub-Gaussian functional linear process. %$\{\bW_t(\cdot)\}$ and $\{Y_t(\cdot)\},$ and also Conditions~{\rm \ref{cond_sflr}(ii), \ref{cond_coef_fflr}} hold. 
	Let $d\asymp \tilde d.$ If the regularization parameter $\gamma_n \asymp  s \{d^{\alpha \vee \tilde\alpha + 2}\cM_{W,Y}(n^{-1}\log p)^{1/2}$ $+d^{-\tau+1/2}\},$ then the estimate $\widehat \bbeta(\cdot,\cdot)$ satisfies
	\begin{equation}
		\label{rate_fflr}
		\sum_{j=1}^p\|\hat\beta_j -\beta_{0j}\|_{\cS} = O_\p\bigg\{ \mu^{-2} s^2 \bigg( d^{\alpha +\alpha\vee\tilde\alpha +2}\cM_{W,Y} \sqrt{\frac{\log p}{n}} + d^{\alpha - \tau + 1/2}\bigg) \bigg\}\,.
	\end{equation}
\end{theorem}
\begin{remark}
	\label{rmk.thm5}
	(i) With the same expression of $\bG$ for both SFLR and FFLR, Condition~\ref{cond_sflr}(ii) is required in both Theorems~\ref{thm_sflr} and \ref{thm_fflr}. Note we can further remove the assumption $d\asymp \tilde d,$ and establish the general convergence rate as a function of $d, \tilde d $ and other parameters.\\
	\textcolor{black}{(ii) The rate for the autocovariance-based estimator in (\ref{rate_fflr}) is slightly slower than that for the covariance-based estimator in \cite{fang2022} by a multiplicative factor $d^{\alpha/2}.$}
	%where the convergence results can be established under a general setting. %present the convergence result in a general setting.
\end{remark}

\subsection{High-dimensional VFAR}
\label{sec_vfar}
%\subsubsection{Estimation}
%\label{sec_est_vfar}
The high-dimensional VFAR of a fixed lag order $H,$ namely VFAR($H$), takes the form of
\begin{equation}
	\label{model_vfar}
	\bX_{t}(v) = \sum_{h' = 1}^H \int_{\cU}\bA_{0}^{(h')}(u,v) \bX_{t-h'}(u)\,{\rm d}u  + \bvarepsilon_{t}(v)\,, ~~t =H+1, \ldots,n\,,
\end{equation}
where $\{\bX_t(\cdot)\}$ satisfy model~(\ref{error.model}), the errors $\bvar_t(\cdot)=\{\varepsilon_{t1}(\cdot), \dots, \varepsilon_{tp}(\cdot)\}^{\T}$ are i.i.d. sampled from a $p$-vector of mean-zero random functions, independent of $\bX_{t-1}(\cdot), \bX_{t-2}(\cdot), \dots,$ and
$\bA_0^{(h')}=\{A_{0,jj'}^{(h')}(\cdot,\cdot)\}_{j,j' \in [p]}$ is the unknown functional transition matrix at lag $h'.$ \textcolor{black}{
	In the special case $H=1$ with $\bA_0=\bA_0^{(1)},$ Theorem~3.1 of \cite{Bbosq1} ensures the stationarity of $\{\bX_t(\cdot)\}$ 
	if there exists an integer $l_0$ such that $\sup_{\|\bbf\|\leq1} \| \bA_0^{l_0}(\bbf) \| <1$ for $\bbf \in \mathbb{H}^p.$ According to \cite{guo2022}, all VFAR($H$) models can be reformulated as a VFAR(1) model and hence it is not hard to adjust the stationarity condition for the general case $H>1.$} To make a feasible fit to (\ref{model_vfar}) under a high-dimensional regime based on observed curves $\{\bW_t(\cdot)\}_{t \in [n]},$ we assume $\{\bA_0^{(h')}\}_{h'\in [H]}$ is rowwise functional $s$-sparse with $s=\max_{j\in[p]} s_j \ll p.$ To be specific, for the $j$-th row of components in $\{\bA_0^{(h')}\},$ we denote the set of nonzero functions by $S_j=\{(j',h') \in [p] \times [H]: \|A_{0,jj'}^{(h')}\|_{\cS} \neq 0 \}$ and its cardinality by $s_j=|S_j|$ for $j \in [p].$

For each $j \in [p],$ we approximate $X_{tj}(\cdot)$ based on the expansion in (\ref{expansion}) truncated at $d_j.$ %According to \cite{guo2022}, 
With some specific calculations, model~(\ref{model_vfar}) can be rowwisely rewritten as
\begin{equation}
	\label{rep_vfar}
	\bfeta_{tj}^\T= \sum_{h' = 1}^H \sum_{j'=1}^p
	\bfeta_{(t-h')j'}^{\T}\bOmega_{0, jj'}^{(h')} + \br_{tj}^{\T} + \bvarepsilon_{tj}^{\T}\,, ~~j \in [p]\,,
\end{equation}
where $\bOmega_{0,jj'}^{(h')}=\int_{\cU^2}\bpsi_{j'}(u) A_{0,jj'}^{(h')}(u,v)\bpsi_j(v)^{\T}\,{\rm d}u{\rm d}v \in {\eR}^{d_{j'} \times d_j}$ and
$\br_{tj}=(r_{tj1}, \dots, r_{tjd_j})^{\T}$ is the truncation error with each $r_{tjm}=\sum_{h'=1}^H\sum_{j'=1}^p\sum_{l=d_{j'}+1}^{\infty}\eta_{(t-h')j'l}\langle\langle \psi_{j'l}, A_{0,jj'}^{(h')} \rangle, \psi_{jm} \rangle$ for $m \in [d_j].$ Let $\bOmega_{0j} =\{(\bOmega_{0,j1}^{(1)})^{\T}, \dots, (\bOmega_{0,jp}^{(1)})^{\T}, \dots, (\bOmega_{0,j1}^{(H)})^{\T}, \dots, (\bOmega_{0,jp}^{(H)})^{\T})\}^{\T} \in {\eR}^{H\sum_{j'=1}^p d_{j'} \times d_j}.$ We choose $\{\bfeta_{(t-H-h)k}: h \in [L], k \in [p]\}$ as vector-valued instrumental variables, which are assumed to be uncorrelated with the random error $\bvarepsilon_{tj}$ in (\ref{rep_vfar}). Within the framework of (\ref{b_moment}), we %choose $q=pL$  and
assume that $\bOmega_{0j}$ is the unique solution to the following moment equations:
\begin{equation}
	\label{moment_vfar}
	{\bf 0} = \eE\{\bfeta_{(t-H-h)k}\bvarepsilon_{tj}^{\T}\} = \bg_{j,hk}(\bOmega_{0j}) + \bR_{j,hk}\,,~~h \in [L]\,, k \in [p]\,,
	%E\Big\{\bfeta_{(t+h)k} \bzeta_t^{\T} - \sum_{j=1}^p \bfeta_{(t+h)k}\bfeta_{tj}^{\T}\bB_{0j}-\bfeta_{(t+h)k} \br_t^{\T}\Big\}
\end{equation}
where $\bg_{j,hk}(\bOmega_{0j}) = \eE\{\bfeta_{(t-H-h)k}\bfeta_{tj}^{\T}\} - \sum_{h'=1}^H \sum_{j'=1}^p \eE\{\bfeta_{(t-H-h)k}\bfeta_{(t-h')j'}^{\T}\bOmega_{0,jj'}^{(h')}\}$ and $\bR_{j,hk} = -\eE\{\bfeta_{(t-H-h)k}\br_{tj}^{\T}\}.$

%The first step implements the autocovariance-based dimension reduction approach to obtain $\{\widehat\bfeta_{tj}\}_{t \in [n], j \in [p]}$ and $\{\widehat \bpsi_{j}(\cdot)\}_{j \in [p]}.$
Given that estimating the functional sparsity in the $j$-th row of $\{\bA_0^{(h')}\}_{h' \in [H]}$ is equivalent to estimating the block sparsity in $\bOmega_{0j}$ for each $j,$ our goal is to estimate the block sparse matrix $\bOmega_{0j}$ using the empirical versions $\bOmega_j \mapsto \widehat \bg_{j,hk}(\bOmega_j)$ for $h \in [L]$ and $k \in [p],$ where
\begin{equation*}
	\label{emp_vfar}
	\widehat \bg_{j,hk}(\bOmega_j)=  \frac{1}{n-H-h}\sum_{t=H+h+1}^n \widehat\bfeta_{(t-H-h)k}\widehat\bfeta_{tj}^{\T} - \frac{1}{n-H-h}\sum_{t=H+h+1}^n \sum_{h'=1}^H\sum_{j'=1}^p\widehat\bfeta_{(t-H-h)k}\widehat\bfeta_{(t-h')j'}^{\T}\bOmega_{jj'}^{(h')}
\end{equation*}
and $\{\widehat\bfeta_{tj}\}_{t \in [n], j \in [p]}$ are obtained in Step~1.
Step~2 follows (\ref{opt_BRMD}) to formulate the block RMD estimator $\widehat \bOmega_j$ by solving the following optimization task:
\begin{equation*}
	\label{opt_vfar}
	\widehat \bOmega_j = \underset{\bOmega_j}{\arg\min} \sum_{h'=1}^H\sum_{j'=1}^p \|\bOmega_{jj'}^{(h')}\|_{\tF} ~~ \text{subject to}~~\underset{k\in [p], h \in [L]}{\max} \|\widehat \bg_{j,hk}(\bOmega_j)\|_{\tF} \leq \gamma_{nj}\,,
\end{equation*}
where $\gamma_{nj} \geq 0$ is a regularization parameter. Step~3 estimates functional transition matrices by
\begin{equation*}
	\label{beta_est_vfar}
	\hat A_{jj'}^{(h')}(u,v)=\widehat\bpsi_{j'}(u)^{\T}\widehat\bOmega_{jj'}^{(h')}\widehat\bpsi_j(v)\,,~~ (u,v) \in \cU^2\,,~ j,j' \in [p]\,,~ h' \in [H]\,,
\end{equation*}
where $\{\widehat \bpsi_{j}(\cdot)\}_{j \in [p]}$ are obtained in Step~1.
%\subsubsection{VFAR}
%\label{sec.thm_vfar}

We next present convergence analysis of $\{\hat A_{jj'}^{(h')}(\cdot, \cdot): j,j' \in [p], h' \in [H]\}.$ To simplify the notation, we assume the same truncated dimension $d_{j}=d$ across $j \in [p]$. 
For each $j,$ we first express (\ref{moment_vfar}) as below: 
$$  \bg_j(\bOmega_{0j})+\bR_j = \bG_j\bOmega_{0j} + \bg_j({\bf 0}) +\bR_j ={\bf 0}\,,$$
where
$\bg_j = (\bg_{j,11}^{\T}, \dots, \bg_{j,1p}^{\T}, \dots, \bg_{j,L1}^{\T}, \dots, \bg_{j,Lp}^{\T})^{\T},$ $\bR_j = (\bR_{j,11}^{\T}, \dots, \bR_{j,1p}^{\T}, \dots, \bR_{j,L1}^{\T}, \dots, \bR_{j,Lp}^{\T})^{\T}$ and
$\bG_j = (\bG_{j,ii'}) \in \eR^{pLd \times pHd}$ whose $(i,i')$-th block is $\bG_{j,ii'}= \mathbb{E}\{\bfeta_{(t-H-h)k} \bfeta_{(t-h')j'}^{\T}\} \in \eR^{d \times d}$ with $i=(h-1)p+k,$ $k \in [p]$ for $h \in [L]$ and $i'=(h'-1)p+j',$ $j' \in [p]$ for $h' \in [H].$
Applying Theorem~\ref{thm_score} on $\widehat\bG_j$ and $\widehat\bg_j({\bf 0}),$ we can verify Condition~\ref{cond_emc}(i) with the choice of $\epsilon_{n1} \asymp \cM_1^{W} d^{\alpha+2}(n^{-1}\log p)^{1/2}.$ Similar to Condition~\ref{cond_sflr} for SFLR, we then give the following regularity conditions.

\begin{condition}
	\label{cond_vfar}
	(i) For each $j \in [p]$ and $(j',h') \in S_j,$ $A_{0,jj'}^{(h')}(u,v) = \sum_{l,m=1}^{\infty} a^{(h')}_{jj'lm}\psi_{j'm}(u)\psi_{jl}(v)$ and there exists some constant $\tau >\alpha+1/2$ such that $|a^{(h')}_{jj'lm}|\lesssim (l+m)^{-\tau-1/2}$ for $l,m \geq 1;$
	%\end{condition}
	%\begin{condition}
	%\label{cond_G_vfar}
	(ii) For each $j \in [p],$ let $\widetilde \bG_j = (\widetilde \bG_{j,ii'})$ be the normalized version of $\bG_j=(\bG_{j,ii'})$ by replacing each $\bG_{j,ii'}$ by
	$\widetilde \bG_{j,ii'}= \mathbb{E}\{\bD_{k}\bfeta_{(t-H-h)k} \bfeta_{(t-h')j'}^{\T}\bD_{j'}\}$ for $i=(h-1)p+k$ and $i'=(h'-1)p+j'$
	with $k,j' \in [p],$ $h \in [L]$ and $h' \in [H],$ where
	$\bD_j = \text{diag}(\lambda_{j1}^{-1/2}, \dots, \lambda_{jd}^{-1/2}).$ There exist universal constants $\tilde c_j>0$ and $\mu_j >0$ such that  $\sigma_{\max}(m, \widetilde\bG_j)\geq \tilde c_j$ and $\sigma_{\min}(m, \widetilde \bG_j)/\sigma_{\max}(m, \widetilde \bG_j) \geq \mu_j$ for $m=16s_j/\mu_j^2.$
\end{condition}

We finally establish the convergence rate of the VFAR estimate $\{\hat A_{jj'}^{(h')}(\cdot,\cdot)\}_{j,j' \in [p], h' \in [H]}$ in the sense of functional matrix $\ell_{\infty}$ norm as follows.

\begin{theorem}
	\label{thm_vfar}
	Suppose that Conditions~{\rm\ref{cond_fsm}--\ref{cond_eigen}} %hold for sub-Gaussian functional linear process $\{\bW_t(\cdot)\},$ 
	and~{\rm\ref{cond_vfar}} hold. If the regularization parameters satisfy $\gamma_{nj} \asymp  s_j \{ d^{\alpha+2}\cM_1^{W}(n^{-1}\log p)^{1/2} + d^{-\tau+1/2}\}$ for $j \in [p]$ and $\mu=\min_{j\in[p]}\mu_j,$ the estimate $\{\hat A_{jj'}^{(h')}(\cdot,\cdot)\}$ satisfies
	\begin{equation}
		\label{rate_vfar}
		\max_{j \in [p]}\sum_{j'=1}^p\sum_{h'=1}^H\|\hat A_{jj'}^{(h')} - A_{0,jj'}^{(h')}\|_{\cS} = O_\p\bigg\{ \mu^{-2} s^2 \bigg( d^{2\alpha+2}\cM_1^{W} \sqrt{\frac{\log p}{n}} + d^{\alpha - \tau + 1/2}\bigg) \bigg\}\,.
	\end{equation}
\end{theorem}

\textcolor{black}{
	\begin{remark}
		Similar to Remarks~\ref{rmk_thm4}(ii) and \ref{rmk.thm5} (ii) for SFLR and FFLR respectively, the rate for $\{\hat A_{jj'}^{(h')}(\cdot,\cdot)\}$ in (\ref{rate_vfar}) is slightly slower than that for the covariance-based estimator in \cite{guo2022} by the same factor $d^{\alpha/2}.$
	\end{remark}
}

\section{Empirical studies}
\label{sec.emp}

\subsection{Simulation study}
\label{sec.sim}
In this section, we conduct a number of simulations to evaluate the finite-sample performance of the proposed autocovariance-based estimators for SFLR, FFLR and VFAR models.

In each simulated scenario, to mimic the infinite-dimensional feature of signal curves, we generate $X_{tj}(u) = \sum_{l=1}^{25} \eta_{tjl}\psi_{l}(u) = \bfeta_{tj}^{\T}\bpsi(u)$ with $\bfeta_{tj}=(\eta_{tj1},\dots,\eta_{tj25})^{\T}$ and $\bpsi(\cdot)=\{\psi_{1}(\cdot), \dots, \psi_{25}(\cdot)\}^{\T}$ for $t \in [n], j \in [p]$ and $u \in \cU=[0,1],$ where $\{\psi_l(u)\}_{1\leq l \leq 25}$ is formed by $25$-dimensional Fourier basis functions, $1, \sqrt{2} \cos(2\pi l u), \sqrt{2} \sin(2\pi l u)$ %\footnote{50 basis functions?}
for $l = 1,\dots, 12$ and each $\bfeta_{t}=(\bfeta_{t1}^{\T}, \dots, \bfeta_{tp}^{\T})^{\T} \in \eR^{25p}$ is generated from a stationary vector autoregressive (VAR) model, $\bfeta_t = \bOmega \bfeta_{t-1} + \beps_t,$ with block transition matrix $\bOmega = (\bOmega_{jk})_{j,k \in [p]} \in \eR^{25p \times 25p}$ and $\beps_t=(\epsilon_{t1}, \dots, \epsilon_{t 25})^{\T},$ whose components are sampled independently according to $\epsilon_{tj} \sim \mathcal{N}(0, 0.7 - 0.1j)$ for $j=1, \dots, 5$ and $\mathcal{N}(0, j^{-2})$ for $j=6, \dots, 25.$
%components in $\beps_t$ being independent $N(0, 0.5)$ variables.
Therefore, $\bX_t(\cdot)$ follows a VFAR(1) model $\bX_t(v)=\int_{\cU}\bA(u,v)\bX_{t-1}(u)\,{\rm d}u + \bvarepsilon_t(v),$ where
$\bvarepsilon_{tj}(v) = \bpsi(v)^{\T} \beps_{tj}$ and autocoefficient functions satisfy $A_{jk}(u,v) = \bpsi(v)^{\T} \bOmega_{jk} \bpsi(u)$
%\footnote{Should be $\bpsi(v)^{\T} \bOmega_{jk} \bpsi(u)$?}
for $j,k \in [p]$ and $u,v \in \cU.$ In our simulations, we generate $n=100, 200, 400$ serially dependent observations of $p=40, 80$ functional variables. %and set $d=3.$ %It is worth noting thateven the ``low-dimensional" setting with $p=40$ and $n=400,$ the estimation results in high-dimensional estimation problems. For example SFLR requires
%estimating VFAR(1) results in a very high-dimensional task, since, e.g. even under $p=40,n=400,$ one needs to estimate $40^2 \times 3^2=14,400$ parameters based on only 400 observations.
The observed curves are generated from $W_{tj}(u) = X_{tj}(u) + e_{tj}(u)$, where white noise curves $e_{tj}(u) = \sum_{l=1}^5 z_{tjl}\psi_{l}(u),$ $\bz_{tj}=(z_{tj1}, \dots, z_{tj5})^{\T}$ and $\{\bz_{tj}\}_{t\in [n]}$ are sampled independently from  multivariate normal distribution with mean zero and covariance matrix $\text{diag}(1, 0.8, 0.3, 1.5, 1.6).$ For each of the three models, the data is generated as follows.

{\bf VFAR}: We generate block sparse $\bOmega$ with $5\%$ or $10\%$ nonzero blocks for $p=80$ or $p=40,$ respectively. %selected at random.
Specifically, for the $j$-th block row, we set the diagonal block $\bOmega_{jj} = \text{diag}(0.60, 0.59, 0.58, 0.3, 0.2, 6^{-2}, \dots, 25^{-2})$
and randomly choose one off-diagonal block being $0.4 \bOmega_{jj}$ and two off-diagonal blocks being $0.1 \bOmega_{jj}.$
%the off-diagonal blocks $\bOmega_{jk} = 0.5 \bOmega_{jj}$ for $j \neq k$.
Such block sparse design on $\bOmega$ can guarantee the stationarity of the generated VFAR(1) process.
%To guarantee the stationarity of the VFAR(1) process, we rescale $\bOmega$ by $0.9\bOmega/\rho(\bOmega),$ where $\rho(\bOmega)$ is the spectral radius of $\bOmega.$
It is worth noting that estimating VFAR(1) results in a very high-dimensional task, since, e.g. even under the most `low-dimensional' setting with $p=40, n=400$ and truncated dimension $d=3$, one needs to estimate $(40 \times 3)^2=14,400$ parameters based on only 400 observations.
The $p$-vector of functional covariates $\{\bX_{t}(\cdot)\}_{t \in [n]}$ for SFLR and FFLR below are generated in the same way as those for VFAR.

{\bf SFLR}: We generate the scalar responses $\{Y_t\}_{t \in [n]}$ from model~(\ref{model_sflr}), where $\varepsilon_t$'s are independent $\mathcal{N}(0,1)$ variables. For each $j \in S=\{1, \dots, 5\},$ we generate $\beta_{j}(u) =\sum_{l=1}^{25} b_{jl} \psi_l(u)$ for $u \in \cU,$ where $b_{j1}, b_{j2}, b_{j3}$ are sampled from the uniform distribution with support $[-1, -0.5] \cup [0.5, 1]$ and $b_{jl}=(-1)^l l^{-2}$ for $l=4, \dots, 25.$  For $j \in [p] \setminus S,$ we let $\beta_{j}(u)=0.$

{\bf FFLR}: We generate the functional responses $\{Y_t(v): v \in {\cal V}\}_{t \in [n]}$ with ${\cal V}=[0,1]$ from model~(\ref{model_fflr}), where $\varepsilon_t(v)=\sum_{m=1}^{5} g_{tm}\psi_m(v)$ with $g_{tm}$'s being independent $\mathcal{N}(0,1)$ variables. For $j \in S,$ we generate $\beta_{j}(u,v)=\sum_{l,m=1}^{25} b_{jml} \psi_l(u) \psi_m(v)$ for $(u,v) \in \cU \times {\cal V},$ where components in $\{b_{jlm}\}_{1 \leq l,m \leq 3}$ are sampled from the uniform distribution with support $[-1, -0.5] \cup [0.5, 1]$
and $b_{jlm}=(-1)^{l+m} (l+m)^{-2}$ for $l$ or $m = 4, \dots, 25.$ For $j \in [p] \setminus S,$ we let $\beta_{j}(u,v)=0$.

Implementing our proposed autocovariance-based learning framework (AUTO) requires choosing $L$ and $d_{j}$'s. As our simulated results suggest that the estimators are not sensitive to the choice of $L$, we set $L=3$ in simulations. To select $d_j,$ we take the standard approach by selecting the largest $d_j$ eigenvalues of $\widehat K_{jj}$ in (\ref{est.K}) such that %the selected leading eigenvalues can explain more than 90\% of the variation in the trajectory
the cumulative percentage of selected eigenvalues exceeds 90\%. %\footnote{In covariance method, we choose 85\%. The referee may criticize here.}
%When $n$ is sufficiently large, we can also implement the bootstrap method \citep{bathia2010} to select $d_j.$
To choose the regularization parameter(s) for each model and comparison method, there are several possible methods one could adopt such as AIC, BIC and cross-validation. The AIC and BIC methods require the calculation of the effective degrees of freedom, which leads to a very challenging task given the high-dimensional, functional and dependent nature of the model structure and hence is left for future research. In our simulations, we generate a training sample of size $n$ and a separate validation sample of the same size. Using the training data, we compute a series of estimators with 30 different values of the regularization parameters, i.e.
$\{\widehat\bbb_j^{(\gamma_n)}\}_{j \in [p]}$ (or
$\{\widehat\bB_j^{(\gamma_n)}\}_{j \in [p]}$) as a function of $\gamma_n$ for SFLR (or FFLR) and $\{\widehat \bOmega_{jk}^{(\gamma_{nj})}\}_{k \in [p]}$ as a function of $\gamma_{nj}$ for VFAR,
%$\{\widehat\beta_{j}^{(\gamma_n)}\}_{j \in [p]}$ as a function of $\gamma_n$ for SFLR or FFLR and $\{\widehat A_{jk}^{(\gamma_{nj})}\}_{j,k \in [p]}$ as a function of $\gamma_{nj}$ for VFAR,
calculate the squared error between observed and fitted values on the validation set, i.e.
$\sum_{t=1}^n[Y_t -\sum_{j=1}^p \{\widehat \bbb_j^{(\gamma_n)}\}^{\T} \widehat\bfeta_{tj}]^2$ for SFLR,
$\sum_{t=1}^n\|\widehat\bzeta_t -\sum_{j=1}^p \{\widehat \bB_j^{(\gamma_n)}\}^{\T} \widehat\bfeta_{tj}\|^2$ for FFLR and
$\sum_{t=1}^n\|\widehat\bfeta_{tj} -\sum_{k=1}^p (\widehat \bOmega_{jk}^{(\gamma_{nj})})^{\T} \widehat\bfeta_{(t-1)k}\|^2$ for VFAR, and choose the one with the smallest error.
%i.e. $\sum_t\big(Y_t -\sum_{j=1}^p \int_{\cU} W_{tj}(u)\widehat \beta_{j}^{(\gamma_n)}(u)du\big)^2$ for SFLR, $\sum_t\|(Y_t(\cdot) -\sum_{j=1}^p \int_{\cU} W_{tj}(u)\widehat \beta_{j}^{(\gamma_n)}(u,\cdot)du \|^2$ for FFLR and $\sum_t\|(W_{tj}(\cdot) -\sum_{k=1}^p \int_{\cU} W_{(t-1)k}(u)\widehat A_{jk}^{(\gamma_{nj})}(u,\cdot)du \|^2$ for VFAR, and choose the one with the smallest error.

We compare AUTO with the standard covariance-based estimation framework (COV), which proceeds in the following three steps.
%For the covariance-based method,
The first step performs FPCA on $\{W_{tj}(\cdot)\}_{t \in [n]}$ for each $j \in [p],$
where the truncated dimension was selected in the same way as $d_j.$
%chosen such that the selected principal components can explain more than 85\% of the variation in the trajectory.
%The first step performs FPCA so that each $W_{tj}(\cdot)$ is approximated by the finite dimensional truncation, where the truncated dimension for each $j$ was chosen such that the selected principal components can explain more than 85\% of the variation in the trajectory.
Therefore, estimating SFLR and FFLR models are transformed into fitting multiple linear regressions with the univariate response \citep{kong2016} and the multivariate response \citep{fang2022}, respectively and the VFAR estimation is converted to the VAR estimation \citep{guo2022}. The second step considers minimizing the covariance-based criterion, essentially the least squares with the addition of a group lasso type penalty. Such criterion can be optimized using an efficient block fast iterative shrinkage-thresholding algorithm developed in \cite{guo2022}, which converges faster than the commonly adopted block coordinate descent algorithm \citep{fan2015}. The third step recovers functional sparse estimates using estimated eigenfunctions.

\begin{figure}
	\centering
	\includegraphics[width=6.2in]{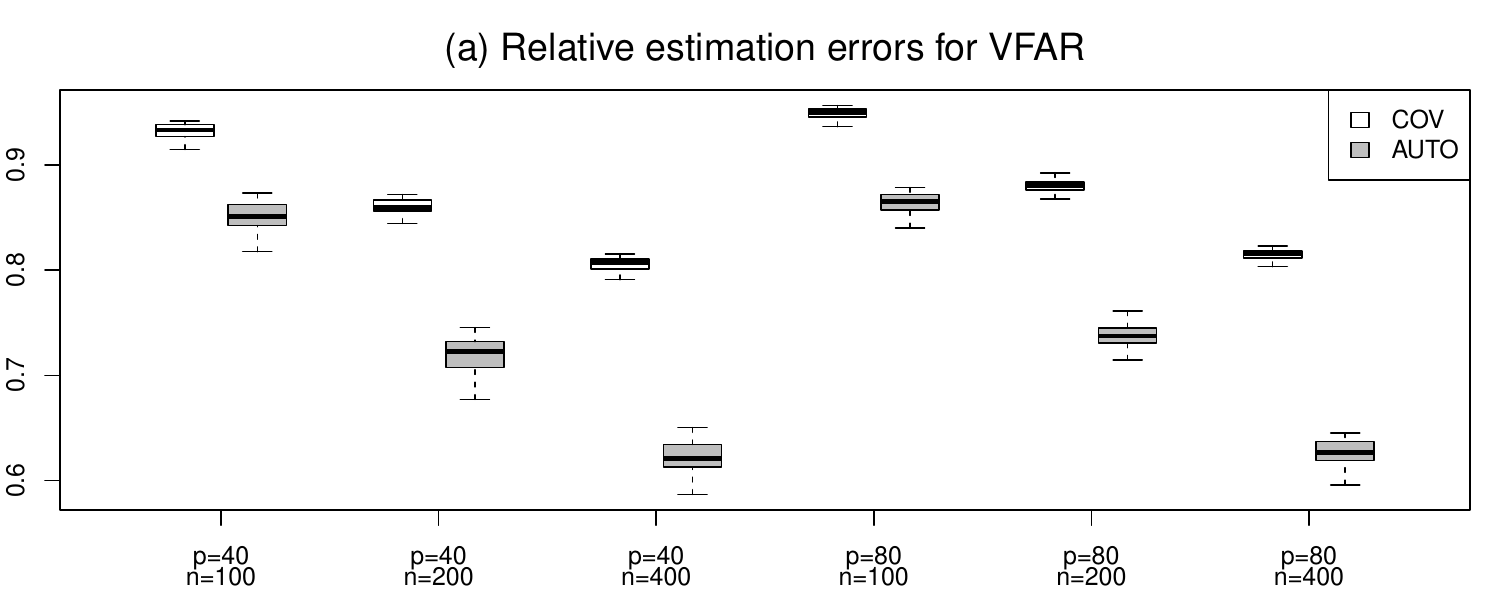}
	\vspace{0.2in}
	
	\includegraphics[width=6.2in]{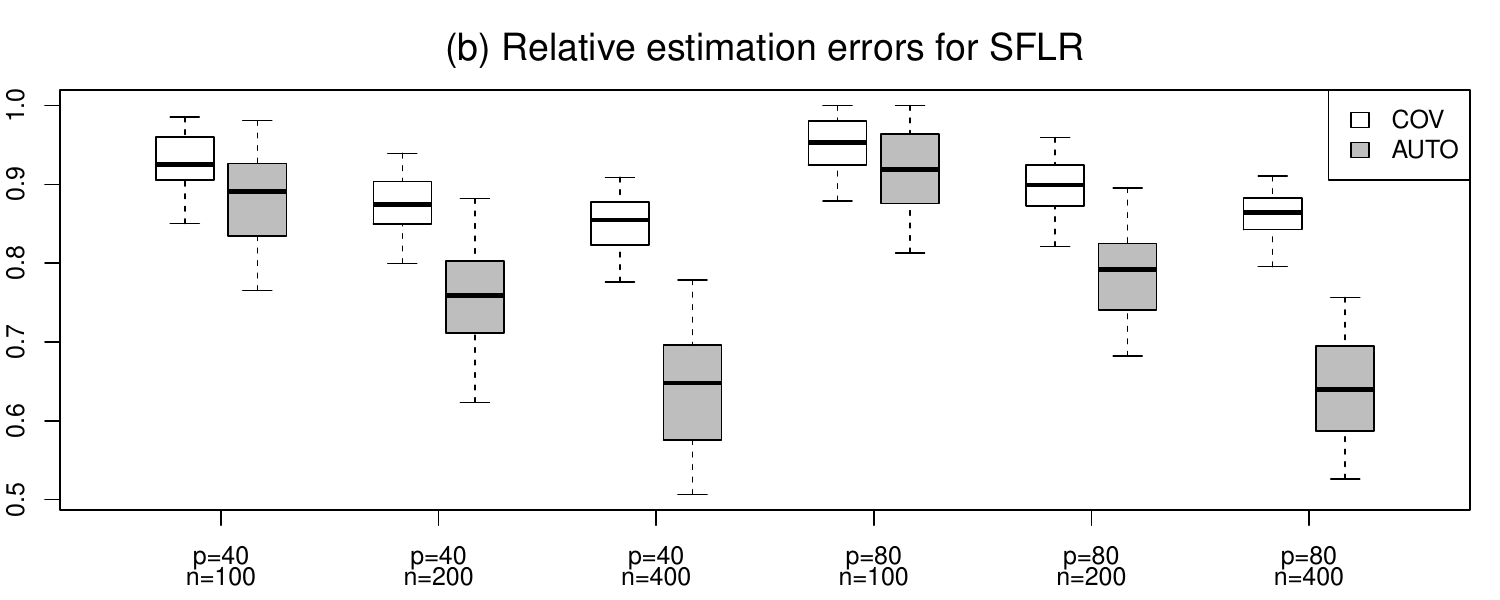}
	\vspace{0.2in}
	
	\includegraphics[width=6.2in]{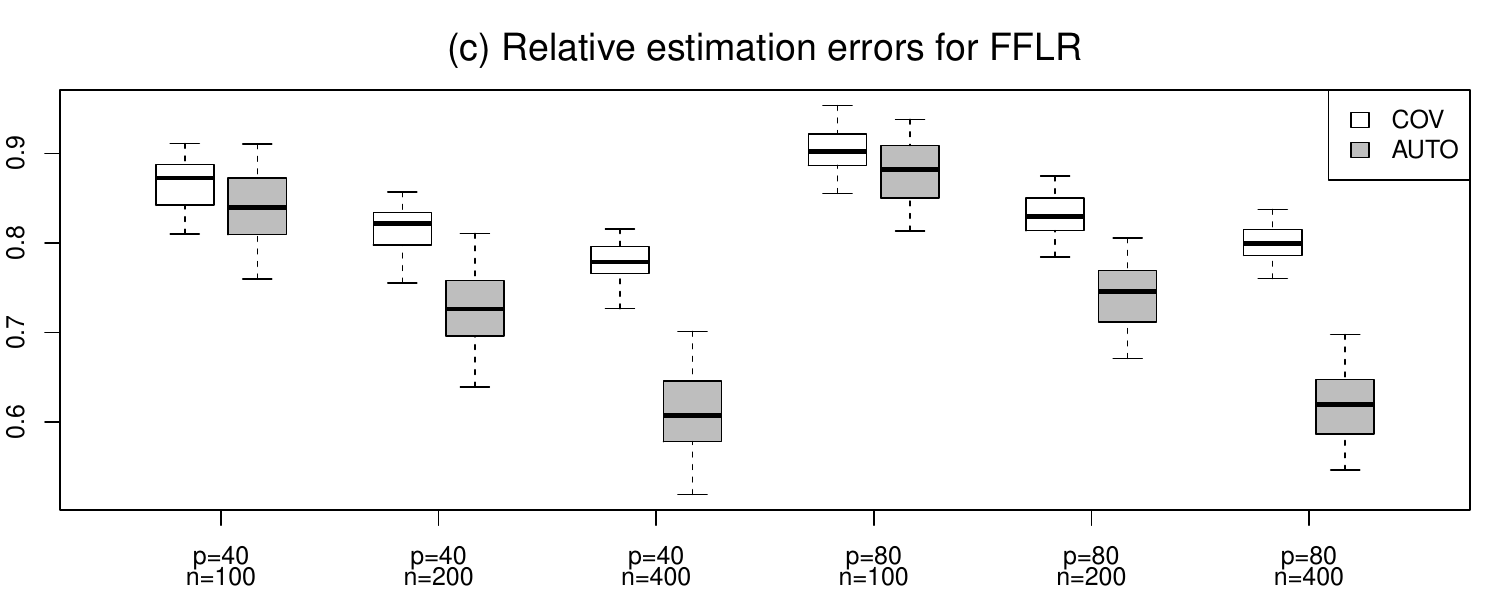}
	
	\caption{\label{box_sim}  The boxplots of relative estimation errors for (a) VFAR, (b) SFLR and (c) FFLR.}
\end{figure}

We examine the performance of COV and AUTO for three models in terms of relative estimation errors, i.e. $\|\widehat\bA -\bA\|_{\tF}/\|\bA\|_{\tF}$ for VFAR, $(\sum_{j=1}^p\|\hat\beta_j -\beta_{0j}\|^2)^{1/2}/(\sum_{j=1}^p\|\beta_{0j}\|^2)^{1/2}$ for SFLR and $(\sum_{j=1}^p\|\hat\beta_j -\beta_{0j}\|_{\cS}^2)^{1/2}/(\sum_{j=1}^p\|\beta_{0j}\|_{\cS}^2)^{1/2}$ for FFLR. We ran each simulation 100 times. Figure~\ref{box_sim} displays boxplots of relative estimation errors for three models. %, while Table~\ref{table_sim} in Appendix~\ref{ap_sim} gives numerical summaries. 
Several conclusions can be drawn from Figure~\ref{box_sim}.
First, AUTO significantly outperforms COV for three models under all scenarios we consider. Second, as discussed in Section~\ref{sec.dr.method}, AUTO provides consistent estimates, while the consistency of COV estimates is jeopardized by the white noise contamination. This can be demonstrated by our empirical results that AUTO provides more substantially improved estimates over COV as $n$ increases from 100 to 400. %especially for SFLR and FFLR. 
Third, the performance of AUTO slightly deteriorates as $p$ increases from 40 to 80, providing empirical evidence to support that the rates in~(\ref{rate_sflr}), (\ref{rate_fflr}) and (\ref{rate_vfar}) for SFLR, FFLR and VFAR models, respectively, all depend on the $(\log p)^{1/2}$ term.

\subsection{Real data analysis}
\label{sec.real}
In this section, we illustrate our developed methodology using a public financial dataset, which was obtained from the WRDS database and consists of high-frequency observations of prices for S\&P~100 index and component stocks (list available in Table~\ref{name.tb} in Appendix~\ref{ap_sim}, we removed several stocks for which the data were not available so that $p=98$ in our analysis) in year~2017 comprising 251 trading days. We obtain one-minute resolution prices \textcolor{black}{by using the last transaction price in each one-minute interval after removing the outliers,}
and hence convert the trading period (9:30--16:00) to minutes $[0, 390].$ We construct cumulative intraday return (CIDR) trajectories \citep{horvath2014}, in percentage, by $W_{tj}(u_k)=100[\log\{P_{tj}(u_k)\} - \log\{P_{tj}(u_1)\}],$ where $P_{tj}(u_k)$ $(t\in[n], j\in[p], k\in[N])$ denotes
the price of the $j$-th stock at the $k$-th minute after the opening time on the $t$-th trading day. \textcolor{black}{We work with mildly smoothed CIDRs obtained by expanding the data with respect to a $45$-dimensional B-spline basis}.
The CIDR curves always start from zero and have nearly the same shape as the original price curves, but make the stationarity assumption more plausible. {\color{black} We performed the functional KPSS test \citep{horvath2014} on CIDR curves for each stock using the R package \verb"fsta" \citep{shang2013}. The p-values are all larger than 1\%, which indicates that there is no overwhelming evidence against the stationarity.}

Our target is to predict the intraday return of the S\&P~100 index based on observed CIDR trajectories of component stocks, $W_{tj}(u), u \in \cU= [0,N]$ up to time $N,$ where, e.g., $N=$ 360 corresponds to 30 minutes prior to the closing time of the trading day. With this in mind, we construct a sparse SFLR model with erroneous functional covariates as follows
\begin{equation}
	\label{model_real}
	Y_t = \sum_{j=1}^p \int_{\cU} X_{tj}(u)\beta_{0j}(u)\,{\rm d}u + \varepsilon_t\,, ~~ W_{tj}(u) = X_{tj}(u) + e_{tj}(u)\,,~~t\in [n]\,, \,j \in [p]\,,
\end{equation}
where $Y_t$ is the intraday return of the S\&P~100 index on the $t$-th trading day, $X_{tj}(\cdot)$ and $e_{tj}(\cdot)$ represent the signal and noise components in $W_{tj}(\cdot),$ respectively. We split the whole dataset into three subsets: training, validation and test sets consisting of the first 171, the subsequent 40 and the last 40 observations, respectively. We apply the validation set approach to select the regularization parameters for AUTO and COV, based on which we estimate sparse functional coefficients in (\ref{model_real}) and calculate the mean squared prediction errors (MSPEs) on the test set. For comparison, we also implement autocovariance-based generalized method-of-moments (AGMM) \citep{chen2022} and covariance-based least squares method (CLS) \citep{hall2007} to fit the unvariate version of (\ref{model_real}) for each component stock, among which we choose the best models leading to the lowest test MSPEs. Finally, we include the null model using the mean of training responses to predict test responses.

\begin{table}
	\caption{\label{table.real} MSPEs up to different current times, $N=$ 300, 315, 330, 345, 360, 370 and 380 minutes, for AUTO and four competing methods. All entries have been multiplied by 100 for formatting reasons. The lowest MSPE for each value of $N$ is in bold font.}
	\begin{center}
		%\vspace{-0.8cm}
		\resizebox{6.0in}{!}{
			\begin{tabular}{c|ccccccc}
				\hline
				Method & $u \leq 300$ & $u \leq 315$ & $u \leq 330$ & $u \leq 345$ & $u \leq 360$ & $u \leq 370$ & $u \leq 380$ \\ \hline
				AUTO & {\bf 5.068} & {\bf 4.936} & {\bf 4.814} & {\bf 4.161} & {\bf 3.892} & {\bf 3.798} & {\bf 3.726} \\
				COV  & 5.487 & 5.360 & 5.222 & 5.090 & 4.976 & 4.927 & 4.882 \\
				AGMM & 6.506 & 6.470 & 6.454 & 6.441 & 6.408 & 6.385 & 6.364 \\
				CLS  & 6.859 & 6.798 & 6.730 & 6.655 & 6.583 & 6.546 & 6.507 \\
				Mean  & 8.832 & 8.832 & 8.832 & 8.832 & 8.832 & 8.832 & 8.832 \\
				\hline
			\end{tabular}
		}
	\end{center}
\end{table}

The resulting test MSPEs for different values of $N$ and all comparison approaches are presented in Table~\ref{table.real}. We observe a few apparent patterns. First, in all scenarios we consider, AUTO provides the best predictive performance, while the autocovariance-based methods are superior to the covariance-based counterparts. Second, the predictive accuracy for functional regression type of methods improves as $N$ approaches to 390 providing more recent information into the covariates. Third, AUTO and COV significantly outperform AGMM and CLS, while Mean gives the worst results. This indicates that using multiple selected functional covariates from the trading histories indeed improves the prediction results.

%\section*{Acknowledgements}

\begin{appendix}

	\section*{Appendix}
	
	This appendix contains further non-asymptotic results in Section~\ref{asec.fcsm}, all technical proofs in Section~\ref{ap_proof} %additional technical lemmas and their proofs in Appendix~\ref{ap_lemma} 
	and list of S\&P 100 stocks in Section~\ref{ap_sim}. 
	
	\appendix
	
	\setcounter{equation}{0}
	\renewcommand{\theequation}{A.\arabic{equation}}
	\section{Further non-asymptotic results}
	\label{asec.fcsm}
	
	To provide theoretical guarantees for the proposed estimators in Sections~\ref{sec_sflr} and \ref{sec_fflr}, we present essential non-asymptotic error bounds on the relevant estimated cross-(auto)covariance terms  based on the functional cross-spectral stability measure \citep{fang2022} between $\{\bW_t(\cdot)\}_{t\in \eZ}$ and  $\tilde p$-vector of mean-zero functional time series (or scalar time series) $\{\bY_t(\cdot)\}_{t\in \eZ}$ (or $\{\bZ_t\}_{t \in \eZ}$). Define $\bSigma_h^{W,Y}(u,v)=\cov\{\bW_{t-h}(u), \bY_t(v)\}$ and
	$\bSigma_h^{W,Z}(u)=\cov\{\bW_{t-h}(u), \bZ_t\}$ for $h \in \eZ$ and $(u,v) \in \cU \times \cal V.$
	
	\begin{condition}
		\label{cond_cfsm}
		{\rm(i)} For $\{\bW_t(\cdot)\}_{t\in\mathbb{Z}}$ and $\{\bY_t(\cdot)\}_{t\in\mathbb{Z}},$ the cross-spectral density function $
		\bbf_\theta^{W,Y} = (2\pi)^{-1}\sum_{h\in \mathbb{Z}}\bSigma_h^{W,Y}e^{-ih\theta}$ for $\theta\in[-\pi,\pi]
		$
		exists and the functional cross-spectral stability measure defined in {\rm(\ref{df_cfsm})} is finite, i.e.
		\begin{equation} \label{df_cfsm}
			\cM^{W,Y} = 2\pi \cdot \esssup\limits_{\theta \in [-\pi,\pi],\bPhi_1\in \mathbb{H}_0^p,\bPhi_2\in \mathbb{H}_0^{\widetilde p}}\frac{|\langle\bPhi_1,\bbf_\theta^{W,Y}(\bPhi_2)\rangle|}{\sqrt{\langle\bPhi_1,\bSigma_0^{W}(\bPhi_1)\rangle}\sqrt{\langle\bPhi_2,\bSigma_0^{Y}(\bPhi_2)\rangle}}<\infty\,,
		\end{equation}
		where $\mathbb{H}_0^p = \{\bPhi \in \mathbb{H}^p:\langle\bPhi,\bSigma_0^W(\bPhi)\rangle \in (0, \infty)\}$ and $\mathbb{H}_0^{\tilde p} = \{\bPhi \in \mathbb{H}^{\tilde p}:\langle\bPhi,\bSigma_0^Y(\bPhi)\rangle \in (0, \infty)\}.$\\
		{\rm(ii)} For $\{\bW_t(\cdot)\}_{t\in\mathbb{Z}}$ and $\{\bZ_t\}_{t\in\mathbb{Z}},$ the cross-spectral density function $
		\bbf_\theta^{W,Z} =(2\pi)^{-1}\sum_{h\in \mathbb{Z}}\bSigma_h^{W,Z}e^{-ih\theta}$ for $\theta\in[-\pi,\pi]
		$
		exists and the functional cross-spectral stability measure defined in {\rm(\ref{df_m_cfsm})} is finite, i.e.
		\begin{equation} \label{df_m_cfsm}
			\mathcal{M}^{W,Z} = 2\pi \cdot \esssup\limits_{\theta \in [-\pi,\pi],\bPhi\in \mathbb{H}_0^p,\bv\in \mathbb{R}_0^{\tilde p}}\frac{|\langle \bPhi, \bbf_\theta^{W,Z}\bv\rangle|}{\sqrt{\langle\bPhi,\bSigma_0^{X}(\bPhi)\rangle}\sqrt{\bv^{\T}\bSigma_0^{Z}\bv}} <\infty\,,
		\end{equation}
		where $\mathbb{R}_0^{\tilde p} = \{\bnu \in \mathbb{R}^{\tilde p}: \bv^{\T}\bSigma_0^Z\bv \in (0, \infty)\}.$
	\end{condition}
	
	In analogy to (\ref{def_sub_fsm}), we can define the functional cross-spectral stability measure of all $k_1$-dimensional subsets of $\{\bW_t(\cdot)\}$ and $k_2$-dimensional subsets of $\{\bY_t(\cdot)\}$ (or $\{\bZ_t\}$) as  $\cM_{k_1,k_2}^{W,Y}$ (or $\cM_{k_1,k_2}^{W,Z}$). It is easy to verify that
	$\cM_{k_1,k_2}^{W,Y} \leq \cM^{W,Y}<\infty$ (or $\cM_{k_1,k_2}^{W,Z} \leq \cM^{W,Z}<\infty$) for $k_1 \in [p]$ and $k_2 \in [\tilde p].$ %See \cite{fang2022} for further discussions. 
	For scalar time series $\{\bZ_t\},$ the non-functional stability measure degenerates to
	%\begin{equation*} \label{def_v_fsm}
	$$\cM^Z = 2 \pi \cdot \underset{\theta\in [-\pi, \pi],\bv \in \mathbb{R}_0^{\tilde p}}{\text{ess}\sup} \frac{\bv^{\T} \bbf_{\theta}^Z\bv}{\bv^{\T}\bSigma_0^Z \bv}\,,$$
	%\end{equation*} 
	which is equivalent to that proposed in \cite{basu2015a}.
	The stability measure of all $k$-dimensional subsets of $\{\bZ_t\},$ i.e. $\cM^Z_{k}$ for $k \in [\tilde p],$ can be defined similarly according to (\ref{def_sub_fsm}). 
	
	For each $k \in [\tilde p],$ we represent $Y_{tk}(\cdot)=\sum_{m=1}^{\infty} \zeta_{tkm} \phi_{km}(\cdot)$ under the Karhunen--Lo\`eve expansion, where $\zeta_{tkm} = \langle Y_{tk} , \phi_{km}\rangle$ and $\{(\theta_{km}, \phi_{km})\}_{m \geq 1}$ are the pairs of eigenvalues and eigenfunctions of $\Sigma_{0,kk}^Y.$ Let $\{(\hat\theta_{km}, \hat\phi_{km})\}_{m \geq 1}$ be the estimated eigenpairs of $\widehat\Sigma_{0,kk}^Y$ and $\hat\zeta_{tkm} = \langle Y_{tk} , \hat\phi_{km}\rangle.$ We next impose a condition on the eigenvalues $\{\theta_{km}\}_{m \geq 1}$ and then develop the deviation bound in elementwise $\ell_{\infty}$-norm on how $\hat{\sigma}_{h,jklm}^{W,Y} = (n-h)^{-1}\sum_{t=h+1}^n\hat{\eta}_{(t-h)jl}\hat{\zeta}_{tkm}$ concentrates around $\sigma_{h,jklm}^{W,Y} = \cov\{\eta_{(t-h)jl},\zeta_{tkm}\},$ which plays a crucial role in the convergence analysis of the FFLR estimate in Section~\ref{sec_fflr}.
	
	\begin{condition}
		\label{cond_Y_eigen}
		(i) For each $k \in [\tilde p],$ $\theta_{k1}> \theta_{k2} > \cdots >0,$ and there exist some positive constants $\tilde c$ and $\tilde \alpha>1$ such that $\theta_{km}-\theta_{k(m+1)} \geq \tilde c m^{-\tilde \alpha-1}$ for $m \geq 1;$ (ii) $\max_{k \in [\tilde p]} \sum_{m=1}^{\infty}\theta_{km} = O(1).$
	\end{condition}
	
	\begin{proposition}
		\label{prop_score_cross}
		Suppose that Conditions~{\rm\ref{cond_fsm}--\ref{cond_eigen}}, {\rm\ref{cond_cfsm}(i)} and {\rm\ref{cond_Y_eigen}} hold,
		%for sub-Gaussian functional linear processes, $\{\bW_t(\cdot)\},$ $\{\bY_t(\cdot)\},$ 
		$\{\bY_t(\cdot)\}_{t\in[n]}$ is sub-Gaussian functional linear process and $h$ is fixed. Let $d$ and $\tilde d$ be positive integers possibly depending on $(n,p,\tilde p)$ and $\cM_{W,Y}= \cM_1^W +\cM_1^Y + \cM_{1,1}^{W,Y}.$
		%If $\log(pdM_1M_2)(M_1^{4\alpha_1+2} \vee M_2^{4\alpha_2+2})\cM_{X,Y}^2/n \rightarrow 0$ as $n,p,d \rightarrow \infty,$
		For $n \gtrsim (d^{2\alpha+2} \vee {\tilde d}^{2\tilde \alpha + 2}) (\cM_{W,Y})^2\log(p\tilde p),$ there exist some positive constants $c_7$ and $c_{8}$ independent of $(n,p, \tilde p, d, \tilde d)$ such that %, with probability greater than $ 1- c_5(p \tilde p d \tilde d)^{-c_6},$ the estimates $\{\widehat \sigma_{h,jklm}^{W,Y}\}$ satisfy
		\begin{equation}
			\label{eq_FPCscores_XY}
			\underset{j\in [p], k \in [\tilde p], l \in [d], m \in [\tilde d]}{\max} \frac{|\hat \sigma_{h,jklm}^{W,Y} - \sigma_{h,jklm}^{W,Y}|}{l^{\alpha+1} \vee m^{\tilde \alpha+1}}\lesssim \cM_{W,Y}\sqrt{\frac{\log(p \tilde p)}{n}}
		\end{equation}
		holds with probability greater than
		$ 1- c_7(p \tilde p)^{-c_{8}}.$
	\end{proposition}
	
	We next consider a mixed process scenario consisting of $\{\bW_t(\cdot)\}$ and $\{\bZ_t\}$ and establish the deviation bound %in elementwise $\ell_{\infty}$-norm 
	on how $\hat\varrho_{h,jkl}^{X,Z} =(n-h)^{-1}\sum_{t=h+1}^n\hat{\eta}_{(t-h)jl}Z_{tk}$ concentrates around $\varrho_{h,jkl}^{X,Z} =\cov\{{\eta}_{(t-h)jl},Z_{tk}\},$ which is essential in deriving the convergence rate of the SFLR estimate in Section~\ref{sec_sflr}.
	
	\begin{proposition}
		\label{prop_score_cross_mix}
		Suppose that Conditions~{\rm\ref{cond_fsm}--\ref{cond_eigen}} and {\rm\ref{cond_cfsm}(ii)} hold, %for sub-Gaussian functional linear process $\{\bW_t(\cdot)\},$ 
		$\{\bZ_t\}_{t \in [n]}$ is sub-Gaussian linear process and $h$ is fixed.
		Let $d$ be a positive integer possibly depending on $(n,p, \tilde p)$ and $\cM_{W,Z}= \cM_1^W +\cM_1^Z + \cM_{1,1}^{W,Z}.$
		%with $\min_{1\leq j\leq p,1\leq k\leq q}\omega_{jM}^\bX\omega_{kM}^\bY >0.$
		%If $\log(pdM_1)M_1^{3\alpha_1+2}\cM_{X,Z}^2/n \rightarrow 0$ as $n,p,d \rightarrow \infty,$
		For $n \gtrsim (\cM_{W,Z})^2 \log(p \tilde p)$, there exist some positive constants $c_9$ and $c_{10}$ independent of $(n, p, \tilde p, d)$ such that %with probability greater than  $ 1- c_7(p \tilde p d)^{-c_8},$ the estimates $\{\widehat \varrho_{h,jkl}^{W,Z}\}$ satisfy
		\begin{equation}
			\label{eq_FPCscores_XY_partial}
			\underset{j\in [p], k \in [\tilde p], l \in [d]}{\max} \frac{|\hat \varrho_{h,jkl}^{W,Z} - \varrho_{h,jkl}^{W,Z}|}{l^{\alpha+1}}\lesssim \cM_{W,Z}\sqrt{\frac{\log(p \tilde p)}{n}}\,,
		\end{equation}
		holds with probability greater than  $ 1- c_9(p \tilde p)^{-c_{10}}.$
	\end{proposition}
\end{appendix}

\setcounter{equation}{0}
\renewcommand{\theequation}{B.\arabic{equation}}
\section{Technical proofs}
\label{ap_proof}
Throughout, we use $c, \bar{c}$, $\tilde{c}$, $\check{c}$ and $\dot{c}$ to denote generic positive finite constants that may be different in different uses.

\subsection{Auxiliary lemmas}
\begin{lemma}
	\label{lemma_delta_S_bd}
	Suppose that Condition~{\rm\ref{cond_emc}(iii)} holds. Then $\|\widehat \bdelta_{S^{\rm c}}\|^{(d, \tilde d)}_1 \leq \| \widehat \bdelta_{S}\|^{(d, \tilde d)}_1$
	with probability at least $1-\delta_{n2}$.
\end{lemma}
\noindent {\it Proof}. It follows from Condition~\ref{cond_emc}(iii) and $\btheta_{0,S^{\rm c}}={\bf 0}$ by definition that with probability at least $1-\delta_{n2},$ $\|\widehat{\btheta}\|^{(d, \tilde d)}_1 \leq \|\btheta_0\|^{(d, \tilde d)}_1=\|\btheta_{0,S}\|^{(d, \tilde d)}_1,$ which implies that
$
\|\btheta_{0,S}\|^{(d, \tilde d)}_1\geq\|\widehat{\btheta}_{S}\|^{(d, \tilde d)}_1 + \|\widehat{\btheta}_{S^{\rm c}}\|^{(d, \tilde d)}_1     \geq\|\btheta_{0,S}\|^{(d, \tilde d)}_1 - \| \widehat{\btheta}_{S} - \btheta_{0,S}\|^{(d, \tilde d)}_1 + \|\widehat{\btheta}_{S^{\rm c}}\|^{(d, \tilde d)}_1$. 
By cancelling $\|\btheta_{0,S}\|^{(d, \tilde d)}_1$ on both sides above, we obtain $\| \widehat{\btheta}_{S^{\rm c}} - \btheta_{0,S^{\rm c}}\|^{(d, \tilde d)}_1 \leq \| \widehat{\btheta}_{S} - \btheta_{0,S}\|^{(d, \tilde d)}_1$. $\hfill\Box$

\begin{lemma}\label{lemma_svd_ieq}
	For $\bA \in \eR^{q \times p}$ with ${\rm rank}(\bA) \leq \min(p,q)$ and $\bx \in \eR^{p \times d},$ let $\bA=\bU \bLambda \bV^{\T}$ be the singular value decomposition of $\bA$ with $\bLambda = {\rm diag}\{\sigma_1, \dots, \sigma_r\}$ and $\sigma_1 \geq \dots \geq \sigma_r > 0.$ Then we have 
	$
	\sigma_r \|\bx\|_{\tF} \leq \|\bA \bx\|_{\tF} \leq \sigma_1 \|\bx\|_{\tF}$. 
\end{lemma}
\noindent{\it Proof.} Let $\bv_j$ denotes the $j$-th row of $\bV^{\T}\bx$ for $j \in [r].$ Write
$
\sigma_r^2 \|\bx\|_{\tF}^2 \leq \|\bA \bx\|_{\tF}^2
= \tr( \bx^{\T} \bA^{\T} \bA \bx)
= \tr( \bx^{\T} \bV \bLambda^2 \bV^{\T} \bx )=( \sum_{j=1}^r \sigma^2_j \bv_j^{\T} \bv_j)^{1/2} \leq  \sigma_1^2 \|\bx\|_{\tF}^2$, 
where, in the inequalities above, we have used $\|\bV^{\T} \bx\|_{\tF} = \|\bx\|_{\tF}$ due to the orthonormality of $\bV.$ Taking the squared root on both sides completes the proof of this lemma. $\hfill\Box$
\medskip

To simplify the notation, we will use $\sigma_{\min}(m)$ and $\sigma_{\max}(m)$ to represent $\sigma_{\min}(m, \bG)$ and $\sigma_{\max}(m, \bG),$ respectively.
\begin{lemma}
	\label{lemma_kappa_bd1}
	It holds that
	\begin{equation*}
		\label{eq_kappa_bd1}
		\kappa(\btheta_0) \geq \max_{m \geq s} \bigg\{ \frac{\sigma_{\min}(m)}{\sqrt{m}} - \frac{2\sigma_{\max}(m)}{\sqrt{m}}  \sqrt{\frac{s}{m}}\bigg\} \frac{s^{-1/2}}{2(1 + 2\sqrt{s/m})}\,.
	\end{equation*}
\end{lemma}
\noindent{\it Proof.} Let $T \subset [p]$ and $\|\bdelta_{T^{\rm c}}\|_1^{(d,\tilde d)} \leq \|\bdelta_T\|_1^{(d,\tilde d)}  $ by (\ref{def.kappa}).
% Without lose of generality, after some rearrangement, we assume $\|\bdelta_{1}\| \ge \|\bdelta_{2}\| \ge \dots.$
Let $T_1$ denote the largest $m$ components of $\{ \| \bdelta_i \|_{\tF} \}_{i \in [p]}$, and $T_2$ be the subsequent $m$-largest, etc. Let $\bV_{\bmu} = \text{diag}(\bmu \otimes \mathbf{1}_d)$ where $\bmu \in \mathbb R^q$ with $\sum_{i=1}^q I( |\mu_i| \neq 0) \le m$ and $\mathbf{1}_d=(1, \dots, 1)^{\T} \in {\eR}^d.$ We let $\|\bmu\| = (\sum_{i=1}^q \mu_i^2)^{1/2}$ and $\|\bmu\|_{\infty}=\max_{i \in [q]}|\mu_i|.$ Then, we have
\begin{align}
	\| \bG \bdelta \|^{(d, \tilde d)}_{\max}
	=&~ \max_{\bmu} \bigg\| \frac{1}{\| \bmu \|} \bV_{\bmu} \bG \bdelta \bigg\|_{\tF}
	\geq \bigg\| \frac{1}{\sqrt{m}\| \bmu \|_\infty} \bV_{\bmu} \bG \bdelta \bigg\|_{\tF} \nonumber\\
	\geq&~ \bigg\| \frac{1}{\sqrt{m}\| \bmu \|_\infty} \bV_{\bmu} \bG_{\cdot, T_1} \bdelta_{T_1} \bigg\|_{\tF} -  \sum_{j \geq 2}\bigg\|\frac{1}{\sqrt{m}\| \bmu \|_\infty} \bV_{\bmu} \bG_{\cdot, T_j} \bdelta_{T_j} \bigg\|_{\tF}\,, \label{Gd_lbd}
\end{align}
where $\bG_{\cdot, T_j}$ is the block submatrix of $\bG$ consisting of all rows and all block columns in $T_j$ of $\bG$ for $j \geq 1.$

Define $\tilde J_1 = \arg \max_{|J| \le m} \sigma_{\min} (\bG_{J,T_1})$. We can let $\bmu = (\mu_i)$ with $\mu_i = 1$ if $i \in \tilde J_1$ and $0$ otherwise, so that $\|\bmu\|_\infty = 1.$ Then the first term in (\ref{Gd_lbd}) becomes
\begin{equation}
	\label{first.min}
	\begin{split}
		\bigg\| \frac{1}{ \sqrt{m}\| \bmu \|_\infty} \bV_{\bmu} \bG_{\cdot T_1} \bdelta_{T_1} \bigg\|_{\tF}
		=&~\bigg\| \frac{1}{\sqrt{m}}\bG_{\tilde J_1,T_1} \bdelta_{T_1}  \bigg\|_{\tF} \\
		\geq&~\frac{\sigma_{\min}(\bG_{\tilde J_1, T_1})}{\sqrt{m}} \|\bdelta_{T_1}\|_{\tF} = \frac{1}{\sqrt{m}} \max_{|J| \le m} \sigma_{\min} (\bG_{J,T_1})\| \bdelta_{T_1} \|_{\tF}\\
		\geq&~\frac{1}{\sqrt{m}} \min_{|M| \le m} \max_{|J| \le m} \sigma_{\min}\left(\bG_{J,M} \right)\| \bdelta_{T_1} \|_{\tF} %\\
		= \frac{\sigma_{\min}(m)}{\sqrt{m}} \| \bdelta_{T_1} \|_{\tF}\,,
	\end{split}
\end{equation}
where the first inequality comes from Lemma~\ref{lemma_svd_ieq}. Define $\tilde J_j = \arg \max_{|J| \le m} \sigma_{\max} (\bG_{J,T_j})$ for each $j \geq 2.$ By the similar arguments as above, the second term in (\ref{Gd_lbd}) becomes
\begin{equation}
	\label{second.max}
	\begin{split}
		\sum_{j \geq 2} \bigg\| \frac{1}{\sqrt{m} \| \bmu \|_\infty} \bV_{\bmu} \bG_{\cdot, T_j} \bdelta_{T_j} \bigg\|_{\tF}
		=&~\frac{1}{\sqrt{m}} \sum_{j \geq 2} \| \bG_{\tilde J_j, T_j} \bdelta_{T_j} \|_{\tF} \\
		\leq&~\frac{1}{\sqrt{m}} \sum_{j \geq 2} \sigma_{\max}(\bG_{\tilde J_j,T_j})\|\bdelta_{T_j}\|_{\tF} = \frac{1}{\sqrt{m}} \sum_{j \geq 2}\max_{|J| \le m} \sigma_{\max} (\bG_{J,T_j})\| \bdelta_{T_j} \|_{\tF} \\
		\leq&~\frac{1}{\sqrt{m}}\max_{|M| \le m} \max_{|J| \le m} \sigma_{\max}\left(\bG_{J,M} \right) \sum_{j \geq 2}\|\bdelta_{T_j} \|_{\tF} =\frac{\sigma_{\max}(m)}{\sqrt{m}} \sum_{j \geq 2} \| \bdelta_{T_j} \|_{\tF}\,.
	\end{split}
\end{equation}
By the construction of sets $\{{T_j}\}_{j\geq1},$ we have $\|\bdelta_{T_j}\|^{(d, \tilde d)}_1 = \sum_{l \in T_j} \|\bdelta_{l}\|_{\tF} \geq m \|\bdelta_{T_{j+1}}\|^{(d, \tilde d)}_{\max} \geq \sqrt{m} \|\bdelta_{T_{j+1}}\|_{\tF},$ which implies that
\begin{equation}
	\label{second.bd}
	\sum_{j \geq 2} \|\bdelta_{T_j}\|_{\tF} \le \frac{1}{\sqrt{m}}\sum_{j \geq 1} \|\bdelta_{T_j}\|^{(d, \tilde d)}_1 \leq \frac{\|\bdelta\|^{(d, \tilde d)}_1}{\sqrt{m}}\,.
\end{equation}
Combining (\ref{first.min}), (\ref{second.max}) and (\ref{second.bd}) yields
\begin{equation}
	\begin{split}
		\| \bG \bdelta \|^{(d, \tilde d)}_{\max}
		&\geq \frac{\sigma_{\min}(m)}{\sqrt{m}} \| \bdelta_{T_1} \|_{\tF} - \frac{\sigma_{\max}(m)}{\sqrt{m}} \|\bdelta\|^{(d, \tilde d)}_1/\sqrt{m} \\
		&\geq \frac{\sigma_{\min}(m)}{\sqrt{m}} \| \bdelta_{T_1} \|_{\tF} - \frac{\sigma_{\max}(m)}{\sqrt{m}} 2 \sqrt{\frac{s}{m}} \|\bdelta_T\|_{\tF} \\
		&= \bigg\{ \frac{\sigma_{\min}(m)}{\sqrt{m}}  - 2\frac{\sigma_{\max}(m)}{\sqrt{m}} \sqrt{\frac{s}{m}} \frac{\|\bdelta_T\|_{\tF}}{\| \bdelta_{T_1} \|_{\tF}}\bigg\} \| \bdelta_{T_1} \|_{\tF}
	\end{split}
	\label{Gd_bd}
\end{equation}
where the second inequality comes from $\|\bdelta\|^{(d, \tilde d)}_1 \leq 2 \|\bdelta_T\|^{(d, \tilde d)}_1 \leq 2\sqrt{s}\|\bdelta_T\|_{\tF}$ with $|T|\leq s.$
This fact together with (\ref{second.bd}) implies that
\begin{equation}
	\label{d_bd}
	\|\bdelta\|_{\tF} \leq \|\bdelta_{T_1}\|_{\tF} + \sum_{j \geq 2} \|\bdelta_{T_j}\|_{\tF} \leq \|\bdelta_{T_1}\|_{\tF} + 2\sqrt{s/m}\|\bdelta_T\|_{\tF} \leq (1 + 2\sqrt{s/m})\|\bdelta_{T_1}\|_{\tF}\,.
\end{equation}
Combining (\ref{Gd_bd}) and (\ref{d_bd}) yields that
\begin{equation}
	\label{kappa_bd1}
	\begin{split}
		\| \bG \bdelta \|^{(d, \tilde d)}_{\max}
		&\geq \bigg\{ \frac{\sigma_{\min}(m)}{\sqrt{m}}  - 2\frac{\sigma_{\max}(m)}{\sqrt{m}} \sqrt{\frac{s}{m}} \bigg\} \frac{\| \bdelta \|_{\tF}}{(1+2\sqrt{s/m})} \\
		&\geq \bigg\{ \frac{\sigma_{\min}(m)}{\sqrt{m}}  - 2\frac{\sigma_{\max}(m)}{\sqrt{m}} \sqrt{\frac{s}{m}} \bigg\} \frac{\| \bdelta \|^{(d, \tilde d)}_1 / \sqrt{s}}{2(1+2\sqrt{s/m})}\,,
	\end{split}
\end{equation}
where the second inequality comes from $\|\bdelta\|_{\tF} \geq \|\bdelta_T\|_{\tF} \geq \|\bdelta_T\|^{(d, \tilde d)}_1/\sqrt{s} \geq \|\bdelta\|^{(d, \tilde d)}_1/\sqrt{4s}.$ We complete our proof by (\ref{def.kappa}) and dividing $\|\bdelta\|^{(d, \tilde d)}_1$ on both sides of (\ref{kappa_bd1}). $\hfill\Box$

\begin{lemma}
	\label{lemma_kappa_bd2} Suppose that Condition~{\rm\ref{cond_svd_bd}} holds. Then there exists some positive constant $c$ such that
	$
	\kappa(\btheta_0) \geq c\mu^2/(24s)$. 
\end{lemma}
\noindent{\it Proof.} Applying Lemma~\ref{lemma_kappa_bd1} and choosing $m=16s/\mu^2$ yields that
\begin{align*}
	% \frac{\| \widetilde{\bG} \bdelta \|^{(d)}_{\max}}{\|\bdelta\|^{(d)}_1}
	\kappa(\btheta_0)
	\geq&~\max_{m \geq s} \frac{\sigma_{\max}(m, \bG)}{\sqrt{m}}\bigg\{\frac{\sigma_{\min}(m, \bG)}{\sigma_{\max}(m, \bG)} - \frac{\mu}{2}\bigg\}\frac{s^{-1/2}}{2(1+\mu/2)} \\
	\geq&~\frac{c \mu}{4\sqrt{s}}\bigg(\mu - \frac{\mu}{2}\bigg)\bigg\{2\bigg(1+\frac{\mu}{2}\bigg)\bigg\}^{-1}s^{-1/2} \geq \frac{c\mu^2}{24s}\,,
\end{align*}
which completes our proof. $\hfill\Box$

For each $j \in [p],$ let $\omega_{j1} \geq \omega_{j2} \geq \cdots > 0$ be the eigenvalues of $\Sigma_{0,jj}^X$ with the corresponding eigenfunctions $\nu_{j1}(\cdot),\nu_{j2}(\cdot), \dots.$ Similarly, let $\{(\omega_{jl}^W, \nu_{jl}^W(\cdot))\}_{l=1}^{\infty}$ be the eigenpairs of $\Sigma_{0,jj}^W.$
%we assume that $X_{tj}(\cdot)$ admits the Karhunen--Lo\`eve expansion$X_{tj}(\cdot) = \sum_{l=1}^{\infty} \xi_{tjl} \nu_{jl}(\cdot),$ where $\xi_{tjl}=\langle X_{tj} , \nu_{jl} \rangle$ corresponds to a sequence of random variables with $\mathbb{E}(\xi_{tjl})=0$ and $\cov(\xi_{tjl},\xi_{tjl'})=\omega_{jl} I(l=l').$ Here $\omega_{j1} \geq \omega_{j2} \geq \cdots \geq 0$ are eigenvalues of $\Sigma_{0,jj}^X$ and $\nu_{j1}(\cdot),\nu_{j2}(\cdot), \dots$ are the corresponding orthonormal eigenfunctions satisfying $\int_{\cU}\Si

\begin{lemma}\label{lemma_eigflp}
	Suppose that Condition~{\rm\ref{cond_flp}} holds. Then we have
	$\omega_0^W = \max_{j}\sum_{l=1}^{\infty} \omega_{jl}^W =O(1).$
\end{lemma}
\noindent{\it Proof.} This lemma follows directly from Lemma~2 of \cite{fang2022} and hence the proof is omitted here. $\hfill\Box$

%\subsection{Lemma~\ref{lemma_Sigmah} and its proof}
\begin{lemma}\label{lemma_Sigmah} For $p \times p$ lag-$h$ autocovariance function of $\{\bW_t(\cdot)\},$ $\{\Sigma_{h,jk}^W(\cdot,\cdot)\}_{j,k \in [p]},$ we have
	$\|\Sigma_{h,jk}^W\|_{\cS} \leq \omega_0^W$ and $\|\Sigma_{h,jk}^W(\psi_{km})\|_{\cS} \leq \omega_{km}^{1/2}(\omega_0^W)^{1/2}$ for $m \geq 1$. 
\end{lemma}
\noindent{\it Proof.} 
This lemma follows directly from Lemma~8 of \cite{guo2022} and hence the proof is omitted here. $\hfill\Box$

%\subsection{Lemma~\ref{lemma_svd_ieq} and its proof}

\subsection{Proof of Theorem~\ref{thm_eigen}}
Along the line of the proofs of Theorem~1 in \cite{fang2022} and Proposition~1 in \cite{guo2022}, we can obtain that for $h \geq 1$
\begin{equation}
	\label{thm_sig_k2}
	\mathbb{P}\bigg\{\bigg|\frac{\langle \bPhi_1, (\widehat{\bf \Sigma}_h^W - \bSigma_h^W)(\bPhi_2)\rangle}
	{\langle \bPhi_1, \bSigma_0^W(\bPhi_1)\rangle +
		\langle \bPhi_2, \bSigma_0^W(\bPhi_2)\rangle}\bigg| >    2\cM_k^W \delta\bigg\} \le 8\exp\{- c n \min(\delta^2, \delta)\}\,.
\end{equation}
For each $j \in [p],$ consider the spectral decomposition $\Sigma_{0,jj}^W(u,v) = \sum_{l=1}^{\infty} \omega_{jl}^W \nu_{jl}^W(u)\nu_{jl}^W(v)$ and $\omega_0 = \max_{j}\sum_{l=1}^{\infty} \omega_{jl}^W =O(1),$ implied from Lemma~\ref{lemma_eigflp}. For each $(j,k,l,m),$ choosing $\bPhi_{1} = (0, \ldots,0,(\omega_{jl}^W)^{-1/2}\nu_{jl}^W, 0,\ldots,0)^{\T}$ and $\bPhi_{2} = (0, \ldots,0,(\omega_{km}^W)^{-1/2}\nu_{km}^W, 0,\ldots,0)^{\T}$ on (\ref{thm_sig_k2}) and following the same procedure to prove Theorem~2 of \cite{guo2022} with the choice of suitable constant $\bar{c}$, we can obtain that
\begin{equation}
	\label{thm_Sigma_compt}
	\mathbb{P}\big\{\|\widehat{\Sigma}_{h,jk}^W - \Sigma_{h,jk}^W\|_{\cS} >   \cM_1^W \delta\big\} \le 8 \exp\{- \bar{c} n \min(\delta^2, \delta)\}\,.
\end{equation}
By (\ref{df.K}), (\ref{est.K}) and Cauchy--Schwartz inequality, we have 
$
\|\widehat K_{j} -K_{j}\|_{\cS}^2 \leq  2 L \sum_{h=1}^L \|\widehat \Sigma_{h,jj}^W -\Sigma_{h,jj}^W\|_{\cS}^2 \|\Sigma_{h,jj}^W\|_{\cS}^2 + L \sum_{h=1}^L \|\widehat \Sigma_{h,jj}^W -\Sigma_{h,jj}^W\|_{\cS}^4$. Let $\Omega_{\omega,jk}^{(h)}=\{\|\widehat \Sigma_{h,jk}^W-\Sigma_{h,jk}^W\|_{\cS}\leq \omega_0\}$ and $\Omega_{jk}^{(h)}=\{\|\widehat \Sigma_{h,jk}^W-\Sigma_{h,jk}^W\|_{\cS}\leq \cM_1^W \delta\}.$ On the event $\Lambda_j=\Omega_{\omega,jj}^{(1)} \cap \dots \cap\Omega_{\omega,jj}^{(L)}\cap\Omega_{jj}^{(1)} \cap \dots \cap \Omega_{jj}^{(L)},$ it follows from the above results and %Lemma~8 of \cite{guo2022} with $\|\Sigma_{h,jj}^W\|_{\cS} \leq \omega_0$
Lemma~\ref{lemma_Sigmah} %for each $j$ 
that
\begin{equation}
	\label{K.bd}
	\|\widehat K_{j} -K_{j}\|_{\cS} \leq \sqrt{3} L \omega_0 \cM_1^W\delta\,.
\end{equation}
Applying (\ref{thm_Sigma_compt}) and choosing $\delta = (\cM_1^W)^{-1}\omega_0$ for $\Omega_{\omega,jj}^{(1)}, \dots, \Omega_{\omega,jj}^{(L)}$ yields %that
$\mathbb{P}(\Lambda_j^{\rm c}) \leq 8L\exp\{- cn \min(\delta^2,\delta)\} + 8L\exp[-c n \min\{(\cM_1^W)^{-2} \omega_0^2, (\cM_1^W)^{-1}\omega_0\}]$. 
Combining the above results, we obtain
\begin{equation}
	\label{coneq_K}
	\mathbb{P}\big(\|\widehat K_{j} -K_{j}\|_{\cS} >  \cM_1^W \delta \big) \leq \bar{c}\exp\{- c n \min(\delta^2,\delta)\} + \bar{c}\exp(- cn)\,.
\end{equation}	
For each $j \in [p],$ it follows from Lemma~4.3 of \cite{Bbosq1} and Condition~\ref{cond_eigen} with $\min_{k \in [l]}\{\lambda_{jk}-\lambda_{j(k+1)}\} \geq c_0 l^{-\alpha-1}$ that
\begin{equation}
	\label{eigenbd}
	\max_{l \in [d]}|\hat\lambda_{jl} - \lambda_{jl}| \leq \|\widehat K_{j}-K_{j}\|_{\cS} ~~ \text{and} ~~ \max_{l \in [d]}\big(\|\hat\psi_{jl}-\psi_{jl}\|/l^{\alpha+1}\big)\leq 2 \sqrt{2} c_0^{-1}\|\widehat K_{j}-K_{j}\|_{\cS}\,.
\end{equation}
Combining (\ref{coneq_K}), (\ref{eigenbd}) and the union bound of probability yields that
\begin{align*}
	&\mathbb{P}\bigg(\max_{j \in [p], l \in [d]}|\hat \lambda_{jl}-\lambda_{jl}| >  \cM_1^W \delta \bigg) \vee \mathbb{P}\bigg\{\max_{j \in [p], l \in [d]} (\|\hat \psi_{jl}-\psi_{jl}\|/l^{\alpha+1}) >  2\sqrt{2}c_0^{-1}\cM_1^W \delta \bigg\} \\
	&~~~~~~~~~~~~~~\leq \bar{c}p \exp\{- cn \min(\delta^2,\delta)\} + \bar{c} p\exp(- cn)\,.
\end{align*}
Let $\delta = \rho \sqrt{n^{-1}\log p} \leq 1.$  Choosing suitable positive constants $\tilde{c}$ and $\check{c}=c\rho^2-1,$ we obtain that (\ref{bd_max_eigen}) holds with probability greater than $1-\tilde{c}p^{-\check{c}},$ which completes the proof of Theorem~\ref{thm_eigen}. $\hfill\Box$
%which can be achieved for sufficiently large $n$ Choosing $\delta=\rho\sqrt{\log(pd)/n}\leq 1$ and suitable constants $c_5, c_6>0,$ we obtain that (\ref{bd_max_eigen}) holds with probability greater than $1-c_5(pd)^{-c_6},$ which completes the proof of Theorem~\ref{thm_eigen}. $\square$

\subsection{Proof of Theorem~\ref{thm_score}}
For each $(j,k,l,m)$ and $h \geq 1,$ we write
\begin{align*}
	\hat \sigma_{jklm}^{(h)}- \sigma_{jklm}^{(h)} =&~\langle \hat \psi_{jl}, \widehat\Sigma_{h,jk}^W(\widehat\psi_{km})\rangle - \langle \psi_{jl}, \Sigma_{h,jk}^W(\psi_{km})\rangle\\
	=&~\langle(\hat\psi_{jl}-\psi_{jl}),\widehat\Sigma_{h,jk}^W(\hat\psi_{km}-\psi_{km})\rangle + \langle \psi_{jl}, (\widehat\Sigma_{h,jk}^W-\Sigma_{h,jk}^W)(\psi_{km})\rangle \\
	&+ \{\langle(\hat\psi_{jl}-\psi_{jl}), (\widehat\Sigma_{h,jk}^W-\Sigma_{h,jk}^W)(\psi_{km})\rangle  + \langle \psi_{jl},(\widehat\Sigma_{h,jk}^W - \Sigma_{h,jk}^W)(\hat \psi_{km} - \psi_{km})\rangle \} \\
	&+ \{\langle(\hat\psi_{jl}-\psi_{jl}), \Sigma_{h,jk}^W(\psi_{km})\rangle  + \langle \psi_{jl}, \Sigma_{h,jk}^W(\hat \psi_{km} - \psi_{km})\rangle \} \\
	=&:J_1 + J_2 + J_3 +J_4\,.
\end{align*}
On the event $\widetilde \Omega_{jk} = \Omega_{\omega,jk}^{(h)} \cap \Omega_{jk}^{(h)} \cap \Lambda_{j} \cap \Lambda_{k},$ it follows from  Lemma~\ref{lemma_Sigmah}, (\ref{K.bd}), (\ref{eigenbd}), the orthonormality of $\{\psi_{jl}\}, \{\psi_{km}\}$
that
$\max_{l,m \in [d]} \{|J_1|/ (l \vee m)^{2(\alpha+1)} \}\lesssim (\cM_1^W)^2 \delta^2$,
$\max_{l,m \in [d]}|J_2|\leq \cM_1^W \delta$,
$\max_{l,m \in [d]} \{|J_3|/ (l \vee m)^{\alpha+1} \}\lesssim \cM_1^W \delta$ and
$\max_{l,m \in [d]} \{|J_4|/(l \vee m)^{\alpha+1} \}\lesssim\cM_1^W \delta$.
Then 
$
\max_{l,m \in [d]}\{\sum_{i=1}^4 |J_i|/(l \vee m)^{\alpha+1}\} \leq c\cM_1^W \delta + \tilde{c} d^{\alpha+1}(\cM_1^W)^2 \delta^2$. 
Applying (\ref{thm_Sigma_compt}) and choosing $\delta = (\cM_1^W)^{-1}\omega_0$ for $\Omega_{\omega,jk}^{(h)},$ $\Omega_{\omega,jj}^{(1)}, \dots, \Omega_{\omega,jj}^{(L)}, \Omega_{\omega,kk}^{(1)}, \dots, \Omega_{\omega,kk}^{(L)}$ yields that 
$
\mathbb{P}(\widetilde \Omega_{jk}^{\rm c}) \leq (16L+8)\exp\{- cn \min(\delta^2,\delta)\}+ (16L+8)\exp[-cn \min\{(\cM_1^W)^{-2} \omega_0^2, (\cM_1^W)^{-1}\omega_0\}]$. 
Combining the above results, choosing suitable positive constants $\bar{c}, \tilde{c}, \check{c},$ and applying the union bound of probability yields that
\begin{equation}
	\label{coneq_cov}
	\begin{split}
		&\mathbb{P}\bigg\{ \underset{j,k\in [p], l,m \in [d]}{\max}\bigg| \frac{\hat \sigma_{jklm}^{(h)} - \sigma_{jklm}^{(h)}}{(l \vee m)^{\alpha+1}}\bigg| >
		\cM_1^W \delta + \bar{c} d^{\alpha+1} (\cM_1^W)^2 \delta^2 \bigg\}\le \tilde{c} p^2 [\exp\{- \check{c} n \min(\delta^2,\delta) \} +  \exp(- \check{c} n)]\,.
	\end{split}
\end{equation}	
Choosing $\delta=\rho_1\sqrt{n^{-1}\log p} \leq 1$ and $1+\bar{c}d^{\alpha+1}\cM_1^W\delta \leq \rho_2$ for some positive constants $\rho_1, \rho_2,$ which can be achieved for sufficiently large $n \gtrsim d^{2\alpha+2} (\cM_1^W)^2\log p,$ it follows from (\ref{coneq_cov}) that there exist positive constants $c,\dot{c}$ such that, with probability greater than $1-c p^{-\dot{c}},$
$$\underset{j,k\in [p], l,m \in [d]}{\max}\bigg| \frac{\hat \sigma_{jklm}^{(h)} - \sigma_{jklm}^{(h)}}{(l \vee m)^{\alpha+1}}\bigg| \leq \rho_1\rho_2\cM_1^W \sqrt{\frac{\log p}{n}}\,,$$
which completes the proof of Theorem~\ref{thm_score}. $\hfill\Box$

\subsection{Proof of Proposition~\ref{prop_score_cross}}
For each $(h, j, k, l, m),$ we write
\begin{align*}
	\hat \sigma_{h,jklm}^{W,Y} - \sigma_{h,jklm}^{W,Y}=& \,\langle(\hat\psi_{jl}-\psi_{jl}),\widehat\Sigma_{h,jk}^{W,Y}(\hat\phi_{km}-\phi_{km})\rangle + \langle \psi_{jl}, (\widehat\Sigma_{h,jk}^{W,Y}-\Sigma_{h,jk}^{W,Y})(\phi_{km})\rangle \\
	&+ \{\langle(\hat\psi_{jl}-\psi_{jl}), (\widehat\Sigma_{h,jk}^{W,Y}-\Sigma_{h,jk}^{W,Y})(\phi_{km})\rangle  + \langle \psi_{jl},(\widehat\Sigma_{h,jk}^{W,Y} - \Sigma_{h,jk}^{W,Y})(\hat \phi_{km} - \phi_{km})\rangle \} \\
	&+ \{ \langle(\hat\psi_{jl}-\psi_{jl}), \Sigma_{h,jk}^{W,Y}(\phi_{km})\rangle  + \langle \psi_{jl}, \Sigma_{h,jk}^{W,Y}(\hat \phi_{km} - \phi_{km})\rangle\} \\
	=&: I_1 + I_2 + I_3 +I_4\,.
\end{align*}
Let $\Omega_{0kk}^Y=\{\|\widehat\Sigma_{0,kk}^Y-\Sigma_{0,kk}^Y\|_{\cS} \leq \cM_1^Y \delta\}$ and $\Omega_{hjk}^{W,Y}=\{\|\widehat\Sigma_{h,jk}^{W,Y}-\Sigma_{h,jk}^{W,Y}\|_{\cS} \leq \cM_{W,Y} \delta\}.$ On the event $\Lambda_j \cap \Omega_{0,kk}^Y \cap \Omega_{h,jk}^{W,Y},$ it follows from 
$\|\langle\Sigma_{h,jk}^{W,Y},\phi_{km}\rangle\|\leq \omega_0^{1/2}\theta_{km}^{1/2}$ and $\|\langle \psi_{jl},\Sigma_{h,jk}^{W,Y}\rangle\|\leq \omega_0^{1/2}\theta_{0}^{1/2},$ derived by the similar techniques to prove Lemma~\ref{lemma_Sigmah}, together with Lemma~\ref{lemma_eigflp},
%Lemma~14 of \cite{fang2022}, 
%Lemma~4.3 of \cite{Bbosq1}, 
(\ref{K.bd}), (\ref{eigenbd}), the orthonormality of $\{\psi_{jl}\}, \{\phi_{km}\}$ and Condition~\ref{cond_Y_eigen} that
$
\max_{l\in [d], m \in [\tilde d]}\{|I_1|/(l^{2(\alpha+1)} \vee m^{2(\tilde \alpha+1)})\} \lesssim  (\cM_1^W)^2\delta^2 +  (\cM_1^Y)^2 \delta^2$, 
$\max_{l\in [d], m \in [\tilde d]}|I_2| \leq  \cM_{W,Y} \delta$, 
$\max_{l\in [d], m \in [\tilde d]}\{|I_3|/(l^{\alpha+1} \vee m^{\tilde \alpha+1})\} \lesssim  \cM_1^W \cM_{W,Y}\delta^2 + \cM_1^Y \cM_{W,Y}\delta^2$ and $
\max_{l\in [d], m \in [\tilde d]}\{|I_4|/(l^{\alpha+1} \vee m^{\tilde \alpha+1})\} \lesssim   \cM_1^W \delta +  \cM_1^Y \delta$. 
Combining the above results and $\cM_{W,Y}=\cM^W_1 + \cM^{Y}_1 + \cM^{W,Y}_{1,1}$ yields that
$
\max_{l\in [d], m \in [\tilde d]}\{\sum_{i=1}^4 |I_i|/(l^{\alpha+1} \vee m^{\tilde \alpha+1})\} \leq c\cM_{W,Y} \delta + \bar{c}(d^{\alpha+1} \vee \tilde d^{\tilde \alpha+1})(\cM_{W,Y} )^2 \delta^2$. 
Following the same developments to prove (\ref{coneq_cov}), we apply (\ref{coneq_K}), Theorem~2, Lemma~24 of \cite{fang2022} and the union bound of probability, choose suitable positive constants $\tilde{c}, \check{c}, \dot{c}$ and hence obtain that
\begin{equation}
	\label{coneq_covWY}
	\begin{split}
		&\mathbb{P}\bigg\{ \underset{j\in [p], k \in [\tilde p], l \in [d], m \in [\tilde d]}{\max}~ \frac{|\hat \sigma_{h,jklm}^{W,Y} - \sigma_{h,jklm}^{W,Y}|}{l^{\alpha+1} \vee m^{\tilde \alpha+1}} >
		\cM_{W,Y} \delta + \tilde{c} (d^{\alpha+1} \vee \tilde d^{\tilde \alpha+1} ) (\cM_{W,Y})^2 \delta^2 \bigg\} \\
		&~~~~~~~~~~\le \check{c} p \tilde p [\exp\{- \dot{c} n \min(\delta^2,\delta)\} + \exp(- \dot{c} n)]\,.
	\end{split}
\end{equation}
Choosing $\delta=\rho_3\sqrt{n^{-1}\log(p \tilde p)} \leq 1$ and $1+\tilde{c}(d^{\alpha+1} \vee \tilde d^{\tilde \alpha +1}) \cM_{W,Y}\delta \leq \rho_4$ for some positive constants $\rho_3, \rho_4,$ which can be achieved for sufficiently large $n \gtrsim (d^{2\alpha+2} \vee {\tilde d}^{2\tilde \alpha + 2}) (\cM_{W,Y})^2\log(p\tilde p),$ it follows from (\ref{coneq_covWY}) that there exist positive constants $c,\bar{c}$ such that, with probability greater than $1-c(p \tilde p)^{-\bar{c}},$
$$\underset{j\in [p], k \in [\tilde p], l \in [d], m \in [\tilde d]}{\max}~ \frac{|\hat \sigma_{h,jklm}^{W,Y} - \sigma_{h,jklm}^{W,Y}|}{l^{\alpha+1} \vee m^{\tilde \alpha+1}} \leq \rho_3\rho_4\cM_{W,Y} \sqrt{\frac{\log(p \tilde p)}{n}}\,,$$
which completes the proof of Proposition~\ref{prop_score_cross}. $\hfill\Box$

\subsection{Proof of Proposition~\ref{prop_score_cross_mix}}
For each $(h, j, k, l),$ we write
$
\hat \varrho_{h,jkl}^{W,Z} - \varrho_{h,jkl}^{W,Z}=\langle(\hat\psi_{jl}-\psi_{jl}),(\hat\Sigma_{h,jk}^{W,Z}-\Sigma_{h,jk}^{W,Z})\rangle + \langle \psi_{jl}, (\widehat\Sigma_{h,jk}^{W,Z}-\Sigma_{h,jk}^{W,Z})\rangle + \langle(\hat\psi_{jl}-\psi_{jl}),\Sigma_{h,jk}^{W,Z}\rangle =: T_1 + T_2 + T_3$. Let $\Omega_{hjk}^{W,Z}=\{\|\widehat\Sigma_{h,jk}^{W,Z}-\Sigma_{h,jk}^{W,Z}\|_{\cS} \leq \cM_{W,Z} \delta\}.$ On the event $\Lambda_j \cap \Omega_{hjk}^{W,Z},$ it follows from (\ref{K.bd}), (\ref{eigenbd}), the orthonormality of $\{\psi_{jl}\}$ and $\|\Sigma_{h,jk}^{WZ}\| \leq \omega_0^{1/2} \sigma_{0,kk}^Z$ that
$
\max_{l \in [d]}(|T_1|/l^{\alpha+1}) \lesssim  \cM_1^W \delta \cM_{W,Z} \delta$, 
$\max_{l \in [d]}|T_2| \leq  \cM_{W,Z} \delta$ and 
$\max_{l \in [d]}(|T_3|/l^{\alpha+1}) \lesssim  \cM_1^W \delta$.
Combining the above results and $\cM_{W,Z}=\cM^W_1 + \cM^{Z}_1 + \cM^{W,Z}_{1,1}$ implies that
$
\max_{l \in [d]}(\sum_{i=1}^3 |T_i|/l^{\alpha+1}) \leq c\cM_{W,Z} \delta + \bar{c}(\cM_{W,Z} )^2 \delta^2$. 
Following the same developments to prove (\ref{coneq_cov}), we apply (\ref{coneq_K}), Remark~3 and Lemma~28 of \cite{fang2022} and the union bound of probability, choose suitable positive constants $\tilde{c}, \check{c}, \dot{c}$ and hence obtain that
\begin{equation}
	\label{coneq_covWZ}
	\begin{split}
		&\mathbb{P}\bigg\{\underset{j\in [p], k \in [\tilde p], l \in [d]}{\max} ~\frac{|\hat \varrho_{h,jkl}^{W,Z} - \varrho_{h,jkl}^{W,Z}|}{l^{\alpha+1}} >
		\cM_{W,Z} \delta + \tilde{c}  (\cM_{W,Z})^2 \delta^2 \bigg\} \le \check{c} p \tilde p [\exp\{- \dot{c} n \min(\delta^2,\delta)\} +  \exp(- \dot{c} n)]\,.
	\end{split}
\end{equation}
Choosing $\delta=\rho_5\sqrt{n^{-1}\log(p \tilde p)} \leq 1$ and $1+\tilde{c} \cM_{W,Z}\delta \leq \rho_6$ for some positive constants $\rho_5, \rho_6,$ which can be achieved for sufficiently large $n \gtrsim (\cM_{W,Z})^2\log(p\tilde p),$ it follows from (\ref{coneq_covWZ}) that there exist positive constants $c,\bar{c}$ such that, with probability greater than $1-c(p \tilde p)^{-\bar{c}},$
$$\underset{j\in [p], k \in [\tilde p], l \in [d]}{\max} ~\frac{|\hat \varrho_{h,jkl}^{W,Z} - \varrho_{h,jkl}^{W,Z}|}{l^{\alpha+1}} \leq \rho_5\rho_6\cM_{W,Z} \sqrt{\frac{\log(p \tilde p)}{n}}\,,$$
which completes the proof of Proposition~\ref{prop_score_cross_mix}. $\hfill\Box$

\subsection{Proof of Theorem~\ref{thm_BRMD}}
By $\bg(\btheta)=\bG\btheta+ \bg({\bf 0})$ and (\ref{linear.g}), we have $\bg(\widehat\btheta)=\bG\widehat\btheta + \bg({\bf 0}),$ $\bG \btheta_0 + \bg({\bf 0}) +\bR={\bf 0}$ and
$\widehat\bg(\widehat\btheta)=\widehat\bG\widehat\btheta + \widehat\bg({\bf 0}).$ Consider event $A=\{\|\widehat \bG - \bG\|_{\max}^{(d,d)} \vee \|\widehat \bg({\bf 0}) - \bg({\bf 0})\|_{\max}^{(d, \tilde d)} \leq \epsilon_{n1}\} \cap \{\|\widehat \bg (\btheta_0)\|_{\max}^{(d,\tilde d)} \leq \gamma_n\}.$ By the union bound of probability and Conditions~\ref{cond_emc}(i) and (iii), this event occurs with probability at least $1-\delta_{n1}-\delta_{n2}.$ On event $A$, we have
\begin{equation}
	\begin{split}
		\label{G.delta_bd}
		\| \bG(\widehat \btheta - \btheta_0) \|^{(d, \tilde d)}_{\max}
		&\leq \| \bg(\widehat\btheta) \|^{(d, \tilde d)}_{\max} + \| \bR \|^{(d, \tilde d)}_{\max}\\
		&\leq \| \widehat{\bg}(\widehat{\btheta}) - \bg(\widehat{\btheta}) \|^{(d, \tilde d)}_{\max} + \| \widehat{\bg}(\widehat{\btheta}) \|^{(d, \tilde d)}_{\max} + \| \bR \|^{(d, \tilde d)}_{\max}\\
		&\leq \| ( \widehat{\bG} - \bG ) \widehat{\btheta} \|_{\max}^{(d, \tilde d)} + \| \widehat{\bg}({\bf 0}) - \bg({\bf 0}) \|^{(d, \tilde d)}_{\max} + \| \widehat{\bg}(\widehat{\btheta}) \|^{(d, \tilde d)}_{\max} + \| \bR \|^{(d, \tilde d)}_{\max} \\
		&\leq \| \widehat{\bG} - \bG \|^{(d,d)}_{\max}  \| \btheta_0 \|^{(d, \tilde d)}_1 + \| \widehat{\bg}({\bf 0}) - \bg({\bf 0}) \|^{(d, \tilde d)}_{\max} + \| \widehat{\bg}(\widehat{\btheta}) \|^{(d, \tilde d)}_{\max} + \| \bR \|^{(d, \tilde d)}_{\max} \\
		&\leq K \epsilon_{n1} + \epsilon_{n1} + \gamma_n + \epsilon_2\,, %\leq (K+2)\epsilon_{n1} + \epsilon_2,
	\end{split}
\end{equation}
where, in the last two inequalities, we have used facts that
$\|(\widehat\bG-\bG)\widehat\btheta\|_{\max}^{(d,\tilde d)} = \max_{i \in [q]} \sum_{j=1}^p \|(\widehat\bG-\bG)_{ij} \widehat \btheta_j\|_{\tF} \leq \max_{i,j} \|(\widehat\bG-\bG)_{ij}\|_{\tF} \sum_{j} \|\widehat\btheta_j\|_{\tF}=\|\widehat\bG-\bG\|_{\max}^{(d,d)} \|\widehat\btheta\|_{1}^{(d, \tilde d)},$ $\|\widehat\btheta\|_1^{(d,\tilde d)} \leq \|\btheta_0\|_1^{(d,\tilde d)} \leq K$ and  $\|\widehat{\bg}(\widehat{\btheta}) \|^{(d, \tilde d)}_{\max} \leq \gamma_n$ by the definition of the block RMD estimator in (\ref{opt_BRMD1}) and $\| \bR \|^{(d, \tilde d)}_{\max} \leq \epsilon_2$ by Condition~\ref{cond_emc}(ii).

%Let $\bdelta = \btheta - \btheta_0.$ We define an identifiable factor
%\begin{equation}
%\label{def.kappa}
%    \kappa(\btheta_0) = \inf_{\bdelta \in C_T, \| \bdelta \|^{(d, \widetilde d)}_1 > 0} \frac{ \| \bG \bdelta \|^{(d, \widetilde d)}_{\max} }{\| \bdelta \|^{(d, \widetilde d)}_{1}},
%\end{equation}
%where $ C_T= \big\{ \bdelta \in \eR^{pd \times \widetilde d}: \|\bdelta_{T^C}\|^{(d, \widetilde d)}_1 \leq \|\bdelta_{T}\|^{(d, \widetilde d)}_1,T \subset \{1, \dots, p\}, |T| = s\big\}.$
On event $A,$ choosing the set $T=S$ in (\ref{def.kappa}) and applying Lemma~\ref{lemma_delta_S_bd} under Condition~\ref{cond_emc}(iii) yields $\|\widehat \bdelta_{S^{\rm c}}\|^{(d, \tilde d)}_1 \leq \| \widehat \bdelta_{S}\|^{(d, \tilde d)}_1$ and hence $\widehat \bdelta \in C_S.$ Then by (\ref{def.kappa}), (\ref{G.delta_bd}) and Lemma~\ref{lemma_kappa_bd2} under Condition~\ref{cond_svd_bd}, we have $\| \widehat{\btheta} - \btheta_0 \|^{(d, \tilde d)}_{1} \leq \kappa(\btheta_0)^{-1} \cdot \| \bG(\widehat{\btheta} - \btheta_0) \|^{(d, \tilde d)}_{\max} \lesssim s\mu^{-2}\{(K+1)\epsilon_{n1} +\gamma_n + \epsilon_2\}$,  which completes the proof. $\hfill\Box$

\subsection{Proof of Proposition~\ref{prop_BRMD}}
Define $\tilde \kappa(\btheta_0)$ by substituting $\bG$ in (\ref{def.kappa})  by $\widetilde \bG.$ By $\widetilde \bG =\bD_x\bG\bD_y$ with $\bD_x$ and $\bD_y$ being diagonal matrices, we have
$
\| \widehat{\btheta} - \btheta_0 \|^{(d, \tilde d)}_{1} \leq\tilde \kappa(\btheta_0)^{-1} \cdot \| \widetilde\bG(\widehat{\btheta} - \btheta_0) \|^{(d, \tilde d)}_{\max}\leq \tilde \kappa(\btheta_0)^{-1} \cdot  \|\bD_x\|_{\max} \|\bD_y\|_{\max} \|\bG(\widehat{\btheta} - \btheta_0) \|^{(d, \tilde d)}_{\max}$.
Following the same procedure to prove Theorem~\ref{thm_BRMD}, we can obtain (\ref{rate_BRMD1}). $\hfill\Box$

\subsection{Proof of Theorem~\ref{thm_sflr}}
We first verify Condition~\ref{cond_emc}(i) for SFLR. For sufficiently large positive constants $c, \bar{c}$, 
define events
\begin{equation}
	\label{event.I1}
	I_1=\bigg\{ \underset{j,k \in [p], h\in [L], l,m \in [d]}{\max} \big|\hat\sigma_{jklm}^{(h)} - \sigma_{jklm}^{(h)}\big| \leq c d^{\alpha+1}\cM_1^{W} \sqrt{\frac{\log p}{n}}\bigg\}\,,
\end{equation}
$$I_2=\bigg\{\underset{k\in [p], h\in [L],m \in [d]}{\max} \bigg|\frac{1}{n-h}\sum_{t=h+1}^n\hat\eta_{(t-h)km}Y_t - \eE\{\eta_{(t-h)km}Y_t\}\bigg|\leq \bar{c} d^{\alpha+1} \cM_{W,Y} \sqrt{\frac{\log p}{n}}\bigg\}\,.$$
On event $I_1 \cap I_2,$ we have
%It follows from Theorem~\ref{thm_score} that with probability at least $1-c_1(pd)^{-c_2},$
\begin{equation}
	\label{verify.G_sflr}
	\begin{split}
		\|\widehat\bG - \bG\|_{\max}^{(d,d)}&= \underset{j,k \in [p],h \in [L]}{\max}\bigg\|\frac{1}{n-h}\sum_{t=h+1}^n\widehat\bfeta_{(t-h)k}\widehat\bfeta_{tj}^{\T} - \mathbb{E}\{\bfeta_{(t-h)k}\bfeta_{tj}^{\T}\}\bigg\|_{\tF}\leq 
		c d^{\alpha+2} \cM_1^W \sqrt{\frac{\log p}{n}}\,,
	\end{split}
\end{equation}
%Moreover, it follows from Proposition~\ref{prop_score_cross_mix} that with probability at least $1-c_1(pd)^{-c2}$
\begin{equation}
	\label{verify.g0_sflr}
	\begin{split}
		\|\widehat\bg({\bf 0}) - \bg({\bf 0})\|_{\max}^{(d,1)}&= \underset{k\in [p], h\in [L]}{\max}\bigg\|\frac{1}{n-h}\sum_{t=h+1}^n\widehat\bfeta_{(t-h)k}Y_t - \eE(\bfeta_{(t-h)j}Y_t)\bigg\|\leq 
		\bar{c} d^{\alpha+3/2} \cM_{W,Y} \sqrt{\frac{\log p}{n}}\,.
	\end{split}
\end{equation}
By Theorem~\ref{thm_score}, Proposition~\ref{prop_score_cross_mix} and the union bound of probability, $\mathbb{P}(I_1 \cap I_2) \geq 1- \tilde{c} p^{-\check{c}}$ for some positive constants $\tilde{c}, \check{c}$. By (\ref{verify.G_sflr}) and (\ref{verify.g0_sflr}), Condition~\ref{cond_emc}(i) can be verified by choosing $\delta_{n1}=\tilde{c} p^{-\check{c}}$ ($p$ depends on $n$) and
\begin{equation}
	\label{verify.emc_sflr}
	\epsilon_{n1} = (c\vee \bar{c})d^{\alpha+2} \cM_{W,Y} \sqrt{\frac{\log p}{n}}\,.
\end{equation}

We next verify Condition~\ref{cond_emc}(ii) for SFLR.  If follows from $r_t=\sum_{j=1}^p \sum_{l=d+1}^{\infty} \eta_{tjl} \langle \psi_{jl}, \beta_{0j} \rangle,$ orthonormality of $\{\psi_{jl}\},$ Cauchy--Schwartz inequality and Condition~\ref{cond_sflr}(i) that
\begin{equation*}
	\begin{split}
		\big\{\|\bR\|_{\max}^{(d,1)} \big\}^2&= \underset{k\in [p], h\in [L]}{\max}\|\eE\{\bfeta_{(t-h)k}r_t\}\|^2 =
		\underset{k,h}{\max} \sum_{m=1}^d\bigg\{\eE\bigg(\eta_{(t-h)km}\sum_{j=1}^p\sum_{l=d+1}^{\infty}\eta_{tjl}a_{jl}\bigg)\bigg\}^2\\
		& \leq \max_{k,h} \sum_{m=1}^d \bigg[\sum_{j \in S}\sum_{l=d+1}^{\infty}\sqrt{\eE\{\eta_{(t-h)km}^2\}\eE(\eta_{tjl}^2)}a_{jl}\bigg]^2\\
		&\leq s^2 \max_{k,j}\sum_{m=1}^d \bigg(\sum_{l=d+1}^{\infty} \lambda_{km}^{1/2}\lambda_{jl}^{1/2} a_{jl}\bigg)^2 \\
		&\leq s^2  \max_{k}\sum_{m=1}^d \lambda_{km} \max_{j} \bigg\{\sum_{l=d+1}^{\infty} \lambda_{jl} \sum_{l=d+1}^{\infty} a_{jl}^2\bigg\}\lesssim \lambda_0^2 s^2 \sum_{l=d+1}^{\infty} l^{-2\tau} =  O(s^2 d^{-2\tau+1})\,.
	\end{split}
\end{equation*}
where the asymptotic inequality comes from Condition~\ref{cond_sflr}(i) and %Lemma~2 of \cite{fang2022} under Condition~\ref{cond_flp} that 
$\lambda_0=\max_{j}\sum_{l=1}^{\infty}\lambda_{jl} = O(1)$ implied by some calculations based on (\ref{df.K}) and Lemma~\ref{lemma_eigflp}.
Therefore
\begin{equation}
	\label{verify.R_sflr}
	\|\bR\|_{\max}^{(d,1)} \leq \dot{c}s d^{-\tau+1/2}=\epsilon_2\,.
\end{equation}
By the similar technique above and Condition~\ref{cond_sflr}(i),
\begin{equation}
	\label{theta.bd_sflr}
	\|\bbb_0\|_1^{(d,1)}= \sum_{j \in S}\bigg(\sum_{l=1}^d a^2_{jl}\bigg)^{1/2}  \lesssim s \max_{j \in S} \bigg(\sum_{l=1}^d l^{-2\tau}\bigg)^{1/2}=O(s)\,.
\end{equation}

Finally, we verify Condition~\ref{cond_emc}(iii) for SFLR. On event $I_1 \cap I_2,$ combining (\ref{verify.emc_sflr}) (\ref{verify.R_sflr}) and (\ref{theta.bd_sflr}) yields that
\begin{equation*}
	\begin{split}
		% \label{verify.gamma}
		\| \widehat\bg(\bbb_0) \|^{(d, 1)}_{\max} &\leq \| \widehat{\bg}(\bbb_0) - \bg(\bbb_0) \|^{(d,1)}_{\max} + \| \bR \|^{(d, 1)}_{\max}\\
		&\leq \| ( \widehat{\bG} - \bG ) \bbb_0 \|^{(d, 1)} + \| \widehat{\bg}({\bf 0}) - \bg({\bf 0}) \|^{(d, 1)}_{\max} + \| \bR \|^{(d, 1)}_{\max} \\
		&\leq \| \widehat{\bG} - \bG \|^{(d, d)}_{\max}  \| \bbb_0 \|^{(d,1)}_1 + \| \widehat{\bg}({\bf 0}) - \bg({\bf 0}) \|^{(d,1)}_{\max} + \| \bR \|^{(d, 1)}_{\max} \\
		& \leq c s\Big( d^{\alpha+2} \cM_{W,Y} \sqrt{\frac{\log p}{n}} + d^{-\tau+1/2}\Big) =\gamma_n\,.
	\end{split}
\end{equation*}
By Condition~\ref{cond_eigen} with $\max_j \|\bD_j\|_{\max} \leq \max_j\lambda_{jd}^{-1/2} = O(d^{\alpha/2})$ and Proposition~\ref{prop_BRMD} under Condition~\ref{cond_sflr}(ii), we have
\begin{equation}
	\label{b_rate_sflr}
	\|\widehat\bbb -\bbb_0\|_1^{(d,1)} = O_\p\bigg\{ \mu^{-2} s^2 d^{\alpha}\bigg( d^{\alpha+2} \cM_{W,Y} \sqrt{\frac{\log p}{n}} + d^{-\tau+1/2}\bigg) \bigg\}\,.
\end{equation}
For each $j \in [p],$ let $R_{j}(u) = \sum_{l=d+1}^{\infty} a_{jl} \psi_{jl}(u).$ By the orthonormality of $\{\psi_{jl}\}$ and  $\|R_j\|^2 = \| \sum_{l=d+1}^\infty a_{jl} \psi_{jl} \|^2 = \sum_{l=d+1}^{\infty} a_{jl}^2 \lesssim d^{-2\tau+1}$ for $j \in S$ under Condition~\ref{cond_sflr}(i), we have
\begin{align*}
	\|\hat{\beta}_{j} - \beta_{0j}\|
	=\|\widehat{\bpsi}_{j}^{\T} \widehat{\bbb}_{j} - \bpsi_{j}^{\T} \bbb_{0j} - R_{j} \| \leq&~\|(\widehat{\bpsi}_{j} - \bpsi_{j})^{\T}\widehat{\bbb}_{j}\| +
	\|\bpsi_{j}^{\T} \{\widehat{\bbb}_{j} - \bbb_{0j}\}\|  + \|R_{j}\|\\
	\leq &~d^{1/2} \max_{l \in [d]}\|\hat{\psi}_{jl} - \psi_{jl}\| \|\widehat \bbb_{j} \| + \|\widehat\bbb_j-\bbb_{0j}\| + O( d^{-\tau+1/2})\,, %\Big(\sum_{l=d+1}^{\infty} a_{jl}^2\Big)^{1/2}.
\end{align*}
which implies that
$
\|\widehat\bbeta-\bbeta_0\|_1\leq d^{1/2}\max_{j\in [p], l\in [d]} \|\hat{\psi}_{jl} - \psi_{jl}\| \|\widehat\bbb\|_1^{(d,1)} +\|\widehat\bbb-\bbb_0\|_1^{(d,1)}  + O(s d^{-\tau+1/2})$, where the third term above is of a smaller order of the second term due to (\ref{b_rate_sflr}).
By $\|\widehat \bbb \|_1^{(d,1)} \le \|\widehat \bbb - \bbb_0\|_1^{(d,1)} + \| \bbb_0 \|_1^{(d,1)},$ (\ref{theta.bd_sflr}) and Theorem~\ref{thm_eigen}, the first term above is of a smaller order of the second term. Hence, we obtain (\ref{rate_sflr}) from (\ref{b_rate_sflr}), which completes the proof. $\hfill\Box$

\subsection{Proof of Theorem~\ref{thm_fflr}}
We first verify Condition~\ref{cond_emc}(i) for FFLR. In addition to event $I_1$ in (\ref{event.I1}), we define event
%$$I_1=\left\{ \underset{\underset{1 \leq l,m \leq d}{1\le j,k\le p, 1 \le h\le L}}{\max} \big|\widehat\sigma_{jklm}^{(h)} - \sigma_{jklm}^{(h)}\big| \leq c_1 d^{\alpha+1}\cM^{Z} \sqrt{\frac{\log(pd)}{n}}\right\},$$
$$I_3=\bigg\{\underset{k\in [p], h\in [L], m \in [d], l \in [\tilde d]}{\max} \bigg|\frac{1}{(n-h)}\sum_{t=h+1}^n\hat\eta_{(t-h)km}\hat\zeta_{tl} - \eE\{\eta_{(t-h)km}\zeta_{tl}\}\bigg| \leq \bar{c} d^{\alpha \vee \tilde \alpha+1} \cM_{W,Y} \sqrt{\frac{\log p}{n}}\bigg\}$$
for some sufficiently large $\bar{c}$. On event $I_1 \cap I_3,$ %define $\widetilde\bD=\text{diag}(\theta_{1}^{1/2}, \dots, \theta_{\widetilde d}^{1/2})$ with $\sum_m\theta_m<\infty,$
we have
\begin{equation}
	\label{verify.g0_fflr}
	\begin{split}
		\|\widehat\bg({\bf 0}) - \bg({\bf 0})\|_{\max}^{(d,\tilde d)}&= \underset{k\in [p], h \in L}{\max}\bigg\|\frac{1}{n-h}\sum_{t=h+1}^n\widehat\bfeta_{(t-h)k}\widehat\bzeta_t^{\T} - \eE(\bfeta_{(t-h)j}\bzeta_{t}^{\T})\bigg\|_{\tF} \leq
		\bar{c} d^{\alpha \vee \tilde \alpha + 2}\cM_{W,Y} \sqrt{\frac{\log p}{n}}\,.
	\end{split}
\end{equation}
By Theorem~\ref{thm_score}, Proposition~\ref{prop_score_cross} and the union bound probability, $\mathbb{P}(I_1 \cap I_3) \geq 1- \tilde{c} p^{-\check{c}}$ for some positive constants $\tilde{c},\check{c}$. By (\ref{verify.G_sflr}) and (\ref{verify.g0_fflr}), Condition~\ref{cond_emc}(i) can be verified with the choice of
\begin{equation}
	\label{verify.emc_fflr}
	\epsilon_{n1} = (c\vee \bar{c})d^{\alpha \vee \tilde \alpha + 2} \cM_{W,Y} \sqrt{\frac{\log p}{n}}\,.
\end{equation}

We next verify Condition~\ref{cond_emc}(ii) for FFLR.  If follows from $\br_t=(r_{t1}, \dots, r_{t \tilde d})^{\T}$ with each $r_{tm'}=\sum_{j=1}^p\sum_{l=d + 1}^{\infty}\eta_{tjl}\langle\langle \psi_{jl}, \beta_{0j} \rangle, \phi_{m'}\rangle,$ orthonormality of $\{\psi_{jl}\},$ $\{\phi_{m'}\},$ Cauchy--Schwartz inequality and Condition~\ref{cond_coef_fflr} that
\begin{equation*}
	\begin{split}
		\big\{\|\bR\|_{\max}^{(d,\tilde d)} \big\}^2&= \underset{k \in [p], h \in [L]}{\max}\|\eE\{\bfeta_{(t-h)k}\br_t^{\T}\}\|_{\tF}^2 =
		\underset{k,h}{\max} \sum_{m=1}^d\sum_{m'=1}^{\tilde d }\bigg\{\eE\bigg(\eta_{(t-h)km}\sum_{j=1}^p\sum_{l=d+1}^{\infty}\eta_{tjl}a_{jlm'}\bigg)\bigg\}^2\\
		& \leq \max_{k,h} \sum_{m=1}^d \sum_{m'=1}^{\tilde d }\bigg[\sum_{j \in S}\sum_{l=d+1}^{\infty}\sqrt{\eE\{\eta_{(t-h)km}^2\}\eE(\eta_{tjl}^2)}a_{jlm'}\bigg]^2\\
		&\leq s^2 \max_{k,j}\sum_{m=1}^d\sum_{m'=1}^{\tilde d } \bigg(\sum_{l=d+1}^{\infty} \lambda_{km}^{1/2}\lambda_{jl}^{1/2} a_{jlm'}\bigg)^2 \\
		&\leq s^2  \max_{k}\sum_{m=1}^d \lambda_{km} \max_{j} \bigg\{\sum_{l=d+1}^{\infty} \lambda_{jl} \sum_{m'=1}^{\tilde d }\sum_{l=d+1}^{\infty} a_{jlm'}^2\bigg\}\\
		&\lesssim \lambda_0^2 s^2 \sum_{m'=1}^{\tilde d }\sum_{l=d+1}^{\infty} (l+m')^{-2\tau-1} =  O(s^2 d^{-2\tau+1})\,,
	\end{split}
\end{equation*}
which implies that
\begin{equation}
	\label{verify.R_fflr}
	\|\bR\|_{\max}^{(d,\tilde d)} \leq \dot{c} s d^{-\tau+1/2}=\epsilon_2\,.
\end{equation}
By the similar technique above and Condition~\ref{cond_coef_fflr},
\begin{equation}
	\label{theta.bd_fflr}
	\|\bB_0\|_1^{(d,\widetilde d)}= \sum_{j \in S}\bigg(\sum_{l=1}^d \sum_{m=1}^{\widetilde d} a^2_{jlm}\bigg)^{1/2}  \lesssim s \max_{j \in S}\bigg\{\sum_{l=1}^d\sum_{m=1}^{\widetilde d} (l+m)^{-2\tau-1}\bigg\}^{1/2}=O(s)\,.
\end{equation}

Finally, we verify Condition~\ref{cond_emc}(iii) for FFLR. On event $I_1 \cap I_3,$ combining (\ref{verify.emc_fflr}) (\ref{verify.R_fflr}) and (\ref{theta.bd_fflr}) and applying the similar techniques for SFLR, we have
\begin{equation*}
	\begin{split}
		% \label{verify.gamma}
		\| \widehat\bg(\bB_0) \|^{(d, \tilde d)}_{\max} %&\leq \| \widehat{\bg}(\bbb_0) - \bg(\bbb_0) \|^{(d,1)}_{\max}\\
		%&\leq \| ( \widehat{\bG} - \bG ) \bbb_0 \|^{(d, 1)} + \| \widehat{\bg}({\bf 0}) - \bg({\bf 0}) \|^{(d, 1)}_{\max} + \| \bR \|^{(d, 1)}_{\max} \\
		&\leq \| \widehat{\bG} - \bG \|^{(d, d)}_{\max}  \| \bB_0 \|^{(d,\tilde d)}_1 + \| \widehat{\bg}({\bf 0}) - \bg({\bf 0}) \|^{(d,\tilde d)}_{\max} + \| \bR \|^{(d, \tilde d)}_{\max} \\
		& \leq c s\bigg( d^{\alpha\vee \tilde \alpha + 2} \cM_{W,Y} \sqrt{\frac{\log p}{n}} + d^{-\tau+1/2}\bigg) =\gamma_n\,.
	\end{split}
\end{equation*}
%By Condition~\ref{cond_eigen} with $\max_j \|\bD_j\|_{\max} \leq \max_j\lambda_{jd}^{-1/2} = O(d^{\alpha/2})$ and By
By Condition~\ref{cond_eigen} and Proposition~\ref{prop_BRMD} under Condition~\ref{cond_sflr}(ii), we have
\begin{equation}
	\label{b_rate_fflr}
	\|\widehat\bB-\bB_0\|_1^{(d,\tilde d)} = O_\p\bigg\{ \mu^{-2} s^2 d^{\alpha}\bigg( d^{\alpha \vee \tilde \alpha +2} \cM_{W,Y} \sqrt{\frac{\log p}{n}} + d^{-\tau+1/2}\bigg) \bigg\}\,.
\end{equation}
For each $j \in [p],$ let $R_{j}(u,v) = (\sum_{l=1}^{d}\sum_{m=1}^{\tilde d}-\sum_{l,m=1}^{\infty}) a_{jlm} \psi_{jl}(u) \phi_{m}(v)$ and write
\begin{align*}
	\hat{\beta}_{j}(u,v) - \beta_{0j}(u,v)
	=&~\widehat{\bpsi}_{j}(u)^{\T} \widehat{\bB}_j \widehat{\bphi}(v) - \bpsi_{j}(u)^{\T} \bB_{0j} \bphi(v) + R_{j}(u,v) \\
	=&~\widehat{\bpsi}_{j}(u)^{\T} \widehat{\bB}_{j} \{\widehat{\bphi}(v) - \bphi(v)\} + \{\widehat{\bpsi}_{j}(u) - \bpsi_{j}(u)\}^{\T} \widehat{\bB}_{j} \bphi(v) \\
	& +\bpsi_{j}(u)^{\T} \{\widehat{\bB}_{j} - \bB_{0j}\} \bphi(v) + R_{j}(u,v)\,.
\end{align*}
By Lemma~9 of \cite{guo2022}, we bound the first three terms by
\begin{align}
	\label{I13_fflr}
	&\big\|\widehat \bpsi_j^{\T} \widehat \bB_{j} (\widehat \bphi -\bphi )\big\|_{\cS}\leq
	{\tilde d}^{1/2} \max_{m \in [\tilde d]}\|\widehat \phi_{m} - \phi_{m}\| \|\widehat \bB_{j} \|_{\tF}\,,  \notag\\
	&\big\|(\widehat \bpsi_j-\bpsi_j)^{\T} \widehat \bB_{j}\bphi \big\|_{\cS} \leq
	d^{1/2} \max_{l \in [d]}\|\widehat \psi_{jl} - \psi_{jl}\| \|\widehat \bB_{j} \|_{\tF}\,,\\
	&\big\| \bpsi_j^{\T} (\widehat \bB_{j} - \bB_{0j})\bphi\big\|_{\cS}= \|\widehat \bB_{j} - \bB_{0j}\|_{\tF}\,.\notag
\end{align}
We next bound the fourth term. For $j \in S,$ by the orthonormality of $\{\psi_{jl}\}$ and $\{\phi_m\},$
\begin{equation}
	\label{R.bd_fflr}
	\begin{split}
		\|R_j\|_{\cS}^2
		&= \bigg\|\bigg(\sum_{l=1}^{d}\sum_{m=1}^{\tilde d}-\sum_{l,m=1}^{\infty}\bigg) a_{jlm} \psi_{jl}\phi_{m}\bigg\|_{\cS}^2 \\
		&= O(1)\cdot\bigg(\sum_{l=1}^d\sum_{m=\tilde d+1}^\infty a_{jlm}^2 + \sum_{l=1}^\infty\sum_{m=1}^{\tilde d} a_{jlm}^2\bigg)\\
		&= O(1)\cdot \bigg\{\sum_{l=1}^d\sum_{m=\tilde d+1}^\infty (l+m)^{-2\tau-1} + \sum_{l=1}^\infty\sum_{m=1}^{\tilde d} (l+m)^{-2\tau-1}\bigg\}= O(d^{-2\tau+1})\,.
	\end{split}
\end{equation}
Combining (\ref{I13_fflr}) and (\ref{R.bd_fflr}), we obtain
$
\|\widehat{\bbeta}-\bbeta_0\|_{1}
\leq \|\widehat{\bB}\|_1^{(d,\tilde d)}
\{{\tilde d}^{1/2}\max_{m \in [\tilde d]}\|\widehat{\phi}_{m}-\phi_{m}\| + d^{1/2}\max_{j\in [p],l\in [d]}\|\widehat{\psi}_{jl}-\psi_{jl}\| \}+\|\widehat{\bB} - \bB_{0}\|_1^{(d,\tilde d)} +
O(sd^{-\tau+1/2})$, 
where the third term above is of a smaller order of the second term due to (\ref{b_rate_fflr}).
By $\|\widehat \bB \|_1^{(d,\widetilde d)} \le \|\widehat \bB - \bB_0\|_1^{(d,\tilde d)} + \| \bB_0 \|_1^{(d,\tilde d)},$ (\ref{theta.bd_fflr}) and Theorem~\ref{thm_eigen}, % and Theorem~3 of \cite{guo2022}, 
the first term is of a smaller order of the second term. According to (\ref{b_rate_fflr}), we complete the proof. $\hfill\Box$

\subsection{Proof of Theorem~\ref{thm_vfar}}
For each $j \in [p],$ we first verify Condition~\ref{cond_emc}(i) for VFAR. On event $I_1$ in (\ref{event.I1}),
\begin{equation}
	\label{verify.G_vfar}
	\begin{split}
		\|\widehat\bG_j - \bG_j\|_{\max}^{(d,d)}&= \underset{j',k \in [p],h \in [L], h'\in [H]}{\max}\bigg\|\frac{1}{n-H-h}\sum_{t=H+h+1}^n\widehat\bfeta_{(t-H-h)k}\widehat\bfeta_{(t-h')j'}^{\T} - \eE\{\bfeta_{(t-H-h)k}\bfeta_{(t-h')j'}^{\T}\}\bigg\|_{\tF}\\
		&\leq %d \underset{\underset{1 \leq l,m \leq d}{1\le j,k\le p, 1 \le h\le L}}{\max} \left|\widehat\sigma_{jklm}^{(h)} - \sigma_{jklm}^{(h)}\right| \lesssim
		c d^{\alpha+2} \cM_1^W \sqrt{\frac{\log p}{n}}\,,
	\end{split}
\end{equation}
\begin{equation}
	\label{verify.g0_vfar}
	\begin{split}
		\|\widehat\bg_j({\bf 0}) - \bg_j({\bf 0})\|_{\max}^{(d,d)}&= \underset{k \in [p], h \in [L]}{\max}\bigg\|\frac{1}{n-H-h}\sum_{t=H+h+1}^n\widehat\bfeta_{(t-H-h)k}\widehat\bfeta_{tj}^{\T} - \eE\{\bfeta_{(t-H-h)k}\bfeta_{tj}^{\T}\}\bigg\|_{\tF} \\
		&\leq %d^{1/2}\underset{1\le k\le p, 1 \le h\le L}{\max} \left|\frac{1}{(n-h)}\sum_{t=h}^n\widehat\eta_{(t+h)km}Y_t - \eE\{\eta_{(t+h)km}Y_t\}\right|\lesssim
		c d^{\alpha+2} \cM_1^{W} \sqrt{\frac{\log p}{n}}\,.
	\end{split}
\end{equation}
It follows from Theorem~\ref{thm_score} that $\mathbb{P}(I_1) \geq 1- \tilde{c} p^{-\check{c}}$ for some positive constants $\tilde{c}, \check{c}$. By (\ref{verify.G_vfar}) and (\ref{verify.g0_vfar}), Condition~\ref{cond_emc}(i) can be verified by choosing
\begin{equation}
	\label{verify.emc_vfar}
	\epsilon_{n1} = cd^{\alpha+2} \cM_1^{W} \sqrt{\frac{\log p}{n}}\,.
\end{equation}

We next verify Condition~\ref{cond_emc}(ii) for VFAR. It follows from $\br_{tj}=(r_{tj1}, \dots, r_{tjd})^{\T}$ with each $r_{tjm'}=\sum_{h'=1}^H\sum_{j'=1}^p\sum_{l=d+1}^{\infty}\eta_{(t-h')j'l}\langle\langle \psi_{j'l}, A_{0,jj'}^{(h')} \rangle, \psi_{jm'} \rangle,$ orthonormality of $\{\psi_{jl}\},$ Cauchy--Schwartz inequality and Condition~\ref{cond_vfar}(i) that
\begin{align*}
	\big\{\|\bR_j\|_{\max}^{(d,d)} \big\}^2&= \underset{k \in [p], h \in [L]}{\max}\|\eE\{\bfeta_{(t-H-h)k}\br_{tj}^{\T}\}\|_{\tF}^2 \\
	&=
	\underset{k,h}{\max} \sum_{m=1}^d\sum_{m'=1}^{d }\bigg[\eE\bigg\{\eta_{(t-H-h)km}\sum_{h'=1}^H\sum_{j'=1}^p\sum_{l=d+1}^{\infty}\eta_{(t-h')j'l}a_{jj'lm'}^{(h')}\bigg\}\bigg]^2\\
	& \leq \max_{k,h} \sum_{m=1}^d \sum_{m'=1}^{d }\bigg[\sum_{(j',h') \in S_j}\sum_{l=d+1}^{\infty}\sqrt{\eE\{\eta_{(t-H-h)km}^2\}\eE\{\eta_{(t-h')j'l}^2\}}a_{jj'lm'}^{(h')}\bigg]^2\\
	&\leq s_j^2 \max_{k,j',h'}\sum_{m=1}^d\sum_{m'=1}^{d } \bigg\{\sum_{l=d+1}^{\infty} \lambda_{km}^{1/2}\lambda_{j'l}^{1/2} a_{jj'lm'}^{(h')}\bigg\}^2 \\
	&\leq s_j^2  \max_{k}\sum_{m=1}^d \lambda_{km} \max_{j',h'} \bigg[\sum_{l=d+1}^{\infty} \lambda_{j'l} \sum_{m'=1}^{d}\sum_{l=d+1}^{\infty} \{a_{jj'lm'}^{(h')}\}^2\bigg]\\
	&\lesssim \lambda_0^2 s_j^2 \sum_{m'=1}^{d}\sum_{l=d+1}^{\infty} (l+m')^{-2\tau-1} =  O(s_j^2 d^{-2\tau+1})\,,
\end{align*}
which implies that
\begin{equation}
	\label{verify.R_vfar}
	\|\bR_j\|_{\max}^{(d,d)} \leq \dot{c} s_j d^{-\tau+1/2}=\epsilon_2\,.
\end{equation}
By the similar technique above and Condition~\ref{cond_vfar}(i), we have
\begin{equation}
	\label{theta.bd_vfar}
	\|\bOmega_{0j}\|_1^{(d,d)}= \sum_{(j',h') \in S_j}\bigg[\sum_{l=1}^d \sum_{m=1}^{d} \{a_{jj'lm}^{(h')}\}^2\bigg]^{1/2}  \lesssim s_j \max_{(j',h') \in S_j}\bigg\{\sum_{l=1}^d\sum_{m=1}^{d} (l+m)^{-2\tau-1}\bigg\}^{1/2}=O(s_j)\,.
\end{equation}

Finally, we verify Condition~\ref{cond_emc}(iii) for VFAR. On event $I_1,$ combining (\ref{verify.emc_vfar}), (\ref{verify.R_vfar}), (\ref{theta.bd_vfar}) and applying the similar techniques, we have
\begin{equation*}
	\begin{split}
		%\label{verify.gamma}
		\| \widehat\bg_j(\bOmega_{0j}) \|^{(d,d)}_{\max}
		&\leq \| \widehat\bG_j - \bG_j \|^{(d, d)}_{\max}  \| \bOmega_{0j} \|^{(d,d)}_1 + \| \widehat{\bg}_j({\bf 0}) - \bg_j({\bf 0}) \|^{(d,d)}_{\max} + \| \bR_j \|^{(d, d)}_{\max} \\
		& \leq c s_j\bigg( d^{\alpha + 2} \cM_1^{W} \sqrt{\frac{\log p}{n}} + d^{-\tau+1/2}\bigg) =\gamma_{nj}\,.
	\end{split}
\end{equation*}
By Condition~\ref{cond_eigen} and Proposition~\ref{prop_BRMD} under Condition~\ref{cond_vfar}(ii), we have
\begin{equation}
	\label{b_rate_vfar}
	\|\widehat\bOmega_{j}-\bOmega_{0j}\|_{1}^{(d,d)} = O_\p\bigg\{ \mu_j^{-2} s_j^2 d^{\alpha}\bigg( d^{\alpha +2} \cM_1^{W} \sqrt{\frac{\log p}{n}} + d^{-\tau+1/2}\bigg) \bigg\}\,.
\end{equation}
For each $j' \in [p],$ let $R_{jj'}^{(h')}(u,v) = (\sum_{l=1}^{d}\sum_{m=1}^{d}-\sum_{l,m=1}^{\infty}) a_{jj'lm}^{(h')} \psi_{j'm}(u) \psi_{jl}(v) $ and write
\begin{align*}
	\hat A_{jj'}^{(h')}(u,v) - A_{0,jj'}^{(h')}(u,v)=&~\widehat{\bpsi}_{j'}(u)^{\T} \widehat\bOmega_{jj'}^{(h')} \widehat\bpsi_j(v) - \bpsi_{j'}(u)^{\T} \bOmega_{0,jj'}^{(h')} \bpsi_j(v) + R_{jj'}^{(h')}(u,v) \\
	=&~\widehat{\bpsi}_{j'}(u)^{\T} \widehat{\bOmega}_{jj'}^{(h')} \{\widehat{\bpsi}_j(v) - \bpsi_j(v)\} + \{\widehat{\bpsi}_{j'}(u) - \bpsi_{j'}(u)\}^{\T} \widehat{\bOmega}_{jj'}^{(h')} \bpsi_j(v) \\
	& +\bpsi_{j'}(u)^{\T} \{\widehat{\bOmega}_{jj'}^{(h')} - \bOmega_{0,jj'}^{(h')}\} \bpsi_j(v) + R_{jj'}^{(h')}(u,v)\,.
\end{align*}
%By Lemma~9 of \cite{guo2022}, 
By the same techniques to prove (\ref{I13_fflr}), we bound the first three terms
\begin{equation}
	\label{I13_vfar}
	\begin{split}
		&\big\|\widehat \bpsi_{j'}^{\T} \widehat \bOmega_{jj'}^{(h')} (\hat \bpsi_j -\bpsi_j)\big\|_{\cS}\leq
		{d}^{1/2} \max_{l \in [d]}\|\widehat \psi_{jl} - \psi_{jl}\| \|\widehat \bOmega_{jj'}^{(h')} \|_{\tF}\,,  \\
		&\big\|(\widehat \bpsi_{j'}-\bpsi_{j'})^{\T} \widehat \bOmega_{jj'}^{(h')}\bpsi_j \big\|_{\cS} \leq
		d^{1/2} \max_{m \in [d]}\|\hat \psi_{j'm} - \psi_{j'm}\| \|\widehat \bOmega_{jj'}^{(h')} \|_{\tF}\,,\\
		&\big\| \bpsi_{j'}^{\T}\{\widehat \bOmega_{jj'}^{(h')} - \bOmega_{0,jj'}^{(h')}\}\bpsi_j\big\|_{\cS}= \|\widehat \bOmega_{jj'}^{(h')} - \bOmega_{0,jj'}^{(h')}\|_{\tF}\,.
	\end{split}
\end{equation}
We next bound the fourth term. For $(j',h') \in S_j,$ by the orthonormality of $\{\psi_{jl}\},$
\begin{equation}
	\label{R.bd_vfar}
	\begin{split}
		\big\|R_{jj'}^{(h')}\big\|_{\cS}^2
		&= \bigg\|\bigg(\sum_{l=1}^{d}\sum_{m=1}^{d}-\sum_{l,m=1}^{\infty}\bigg) a_{jj'lm}^{(h')} \psi_{jl}\psi_{j'm}\bigg\|_{\cS}^2 \\
		&= O(1)\bigg[\sum_{l=1}^d\sum_{m=d+1}^\infty \{a_{jj'lm}^{(h')}\}^2\bigg]
		= O(1) \bigg\{\sum_{l=1}^d\sum_{m=d+1}^\infty (l+m)^{-2\tau-1}\bigg\}= O(d^{-2\tau+1})\,.
	\end{split}
\end{equation}
Combining (\ref{I13_vfar}) and (\ref{R.bd_vfar}), we obtain
\begin{align*}
	\max_{j \in [p]}\sum_{j'=1}^p\sum_{h'=1}^H\|\hat A_{jj'}^{(h')} - A_{0,jj'}^{(h')}\|_{\cS}\leq&~\max_j\|\widehat{\bOmega}_{j}\|_{1}^{(d,d)}
	\bigg\{{d}^{1/2}\underset{j\in [p],l\in [d]}{\max}\|\hat{\psi}_{jl}-\psi_{jl}\| + d^{1/2}\underset{j'\in [p], m\in [d]}{\max}\|\hat{\psi}_{j'm}-\psi_{j'm}\| \bigg\} \\
	&+\max_j\|\widehat\bOmega_{j}-\bOmega_{0j}\|_{1}^{(d,d)} +
	O(s_jd^{-\tau+1/2})\,,
\end{align*}
where the third term above is of a smaller order of the second term due to (\ref{b_rate_vfar}).
By $\max_j\|\widehat \bOmega_j \|_1^{(d,d)} \le \max_j\|\widehat \bOmega_j - \bOmega_{0j}\|_1^{(d,d)} + \max_j\| \bOmega_{0j} \|_1^{(d,d)},$ (\ref{theta.bd_vfar}) and Theorem~\ref{thm_eigen}, the first term is of a smaller order of the second term. Applying (\ref{b_rate_vfar}) with $\mu=\min_j\mu_j$ and $s=\max_j s_j$ completes our proof. $\hfill\Box$

%\subsection{Additional technical lemmas}
%\label{ap_lemma}

%\subsubsection{Lemma~\ref{lemma_eigflp} and its proof}

\section{List of S\&P 100 component stocks used in Section~\ref{sec.real}}
\label{ap_sim}
%Table~\ref{table_sim} reports numerical summaries of relative errors for VFAR, SFLR and FFLR. 
%Table~\ref{name.tb} presents the list of S\&P 100 component stocks used in Section~\ref{sec.real}.

\begin{table}
	\caption{\label{name.tb} List of S\&P~100 stocks.}
	\vspace{0.5cm}
	\linespread{1.25}
	% \begin{adjustbox}{angle=0}
		\scriptsize
		\begin{tabular}{ll|ll}
			\hline
			Ticker & Company name                         & Ticker & Company name                       \\ \hline
			AAPL   & APPLE INC                            & JPM    & JPMORGAN CHASE \& CO               \\
			ABBV   & ABBVIE INC                           & KHC    & KRAFT HEINZ                        \\
			ABT    & ABBOTT LABORATORIES                  & KMI    & KINDER MORGAN INC                  \\
			ACN    & ACCENTURE PLC CLASS A                & KO     & COCA-COLA                          \\
			AGN    & ALLERGAN                             & LLY    & ELI LILLY                          \\
			AIG    & AMERICAN INTERNATIONAL GROUP INC     & LMT    & LOCKHEED MARTIN CORP               \\
			ALL    & ALLSTATE CORP                        & LOW    & LOWES COMPANIES INC                \\
			AMGN   & AMGEN INC                            & MA     & MASTERCARD INC CLASS A             \\
			AMZN   & AMAZON COM INC                       & MCD    & MCDONALDS CORP                     \\
			AXP    & AMERICAN EXPRESS                     & MDLZ   & MONDELEZ INTERNATIONAL INC CLASS A \\
			BA     & BOEING                               & MDT    & MEDTRONIC PLC                      \\
			BAC    & BANK OF AMERICA CORP                 & MET    & METLIFE INC                        \\
			BIIB   & BIOGEN INC INC                       & MMM    & 3M                                 \\
			BK     & BANK OF NEW YORK MELLON CORP         & MO     & ALTRIA GROUP INC                   \\
			BLK    & BLACKROCK INC                        & MON    & MONSANTO                           \\
			BMY    & BRISTOL MYERS SQUIBB                 & MRK    & MERCK \& CO INC                    \\
			C      & CITIGROUP INC                        & MS     & MORGAN STANLEY                     \\
			CAT    & CATERPILLAR INC                      & MSFT   & MICROSOFT CORP                     \\
			CELG   & CELGENE CORP                         & NEE    & NEXTERA ENERGY INC                 \\
			CHTR   & CHARTER COMMUNICATIONS INC CLASS A   & NKE    & NIKE INC CLASS B                   \\
			CL     & COLGATE-PALMOLIVE                    & ORCL   & ORACLE CORP                        \\
			COF    & CAPITAL ONE FINANCIAL CORP           & OXY    & OCCIDENTAL PETROLEUM CORP          \\
			COP    & CONOCOPHILLIPS                       & PCLN   & THE PRICELINE GROUP INC            \\
			COST   & COSTCO WHOLESALE CORP                & PEP    & PEPSICO INC                        \\
			CSCO   & CISCO SYSTEMS INC                    & PFE    & PFIZER INC                         \\
			CVS    & CVS HEALTH CORP                      & PG     & PROCTER \& GAMBLE                  \\
			CVX    & CHEVRON CORP                         & PM     & PHILIP MORRIS INTERNATIONAL INC    \\
			DHR    & DANAHER CORP                         & PYPL   & PAYPAL HOLDINGS INC                \\
			DIS    & WALT DISNEY                          & QCOM   & QUALCOMM INC                       \\
			DUK    & DUKE ENERGY CORP                     & RTN    & RAYTHEON                           \\
			EMR    & EMERSON ELECTRIC                     & SBUX   & STARBUCKS CORP                     \\
			EXC    & EXELON CORP                          & SLB    & SCHLUMBERGER NV                    \\
			F      & F MOTOR                              & SO     & SOUTHERN                           \\
			FB     & FACEBOOK CLASS A INC                 & SPG    & SIMON PROPERTY GROUP REIT INC      \\
			FDX    & FEDEX CORP                           & T      & AT\&T INC                          \\
			FOX    & TWENTY-FIRST CENTURY FOX INC CLASS B & TGT    & TARGET CORP                        \\
			FOXA   & TWENTY-FIRST CENTURY FOX INC CLASS A & TWX    & TIME WARNER INC                    \\
			GD     & GENERAL DYNAMICS CORP                & TXN    & TEXAS INSTRUMENT INC               \\
			GE     & GENERAL ELECTRIC                     & UNH    & UNITEDHEALTH GROUP INC             \\
			GILD   & GILEAD SCIENCES INC                  & UNP    & UNION PACIFIC CORP                 \\
			GM     & GENERAL MOTORS                       & UPS    & UNITED PARCEL SERVICE INC CLASS B  \\
			GOOG   & ALPHABET INC CLASS C                 & USB    & US BANCORP                         \\
			GS     & GOLDMAN SACHS GROUP INC              & UTX    & UNITED TECHNOLOGIES CORP           \\
			HAL    & HALLIBURTON                          & V      & VISA INC CLASS A                   \\
			HD     & HOME DEPOT INC                       & VZ     & VERIZON COMMUNICATIONS INC         \\
			HON    & HONEYWELL INTERNATIONAL INC          & WBA    & WALGREEN BOOTS ALLIANCE INC        \\
			IBM    & INTERNATIONAL BUSINESS MACHINES CO   & WFC    & WELLS FARGO                        \\
			INTC   & INTEL CORPORATION CORP               & WMT    & WALMART STORES INC                 \\
			JNJ    & JOHNSON \& JOHNSON                   & XOM    & EXXON MOBIL CORP                   \\ \hline
		\end{tabular}
		% \end{adjustbox}
\end{table}

%\linespread{1.05}\selectfont
%\bibliography{paperbib}
%\bibliographystyle{dcu}

% \small
\singlespacing
\small
\bibliographystyle{dcu}
\bibliography{paperbib}

%%%%%%%%%%%%%%%%%%%%%%%%%%%%%%%%%%%%%%%%%%%%%%%%%%%%

%\setcounter{equation}{-1}
%\setcounter{section}{0}
%\renewcommand{\thesection}{\Alph{section}}
%

%\newpage
%\onehalfspacing
%\normalsize

%\begin{center}
%	{\bf  \LARGE
	%		Supplementary Material}  \\
%{\large \bf  (For online publication only)}
%\end{center}

%\setcounter{page}{1}

%\bigskip

%\small
%\bibliographystyle{chicago}
%\bibliography{bib}

\end{document}